\documentclass[a4paper, 11pt]{amsart}
\usepackage{enumerate}
\usepackage{amsmath}
\usepackage{amsfonts}
\usepackage{amssymb}
\usepackage{amsmath,amssymb,amsthm}
\usepackage{amsmath,amssymb,amsthm,amscd}
\usepackage[mathscr]{eucal}
\usepackage{color}
\usepackage[colorlinks, linkcolor=red, anchorcolor=blue,citecolor=blue]{hyperref}
\usepackage{cite}
\usepackage{tabularx}
\usepackage{enumitem}
\usepackage{enumerate}
\numberwithin{equation}{section}
\newtheorem{thm}{Theorem}[section]
\newtheorem{lem}[{thm}]{Lemma}
\newtheorem{exa}[{thm}]{Example}
\newtheorem{defn}[{thm}]{Definition}
\newtheorem{prop}[{thm}]{Proposition}
\newtheorem{ass}[{thm}]{Assumption}
\newtheorem{corr}[{thm}]{Corollary}
\newtheorem{rem}{Remark}[section]{\bfseries\upshape}

\newtheorem{con}[{thm}]{Conjecture}

\renewcommand{\oddsidemargin}{5mm}

\makeatletter
\title[Eigenvalues of the Laplacian with Muckenhoupt Weights]{Estimates for Ratios of Low-Order Eigenvalues of the Laplacian with Muckenhoupt Weights on Euclidean Domains}
\author[J. Mao H.-J. Sun \& L. Zeng]{Jing Mao \quad He-Jun Sun \quad Lingzhong Zeng}
\address{Jing Mao
\\ \newline \indent Faculty of Mathematics and Statistics
\\ \newline \indent Hubei University, Wuhan 430062, China. jiner120@163.com }
\address{He-Jun Sun
\\ \newline \indent School of Mathematics and Statistics
\\ \newline \indent Nanjing University of Science and Technology, Nanjing 210094, China. hejunsun@163.com }
\address{Lingzhong Zeng
\\ \newline \indent College of Mathematics and Statistics
\\ \newline \indent Jiangxi Normal University, Nanchang 330022, China. lingzhongzeng@yeah.net }
\date{}
\newcommand{\Xint}[1]{\mathchoice
  {\XXint\displaystyle\textstyle{#1}}%
  {\XXint\textstyle\scriptstyle{#1}}%
  {\XXint\scriptstyle\scriptscriptstyle{#1}}%
  {\XXint\scriptscriptstyle\scriptscriptstyle{#1}}%
  \!\int}
\newcommand{\XXint}[3]{{\setbox0=\hbox{$#1{#2#3}{\int}$}\vcenter{\hbox{$#2#3$}}\kern-0.5\wd0}}
\newcommand{\fint}{\Xint-}
\begin{document}
\maketitle
\begin{abstract}
We study low-order eigenvalue ratio estimates for the Laplacian with
Muckenhoupt $A_2$ weights. For weights bounded above and below, we
prove a universal upper bound for $\lambda_2/\lambda_1$ with
constant $1 + \frac{4}{n}\,\frac{C}{c}$, linear in the ratio of the
weight bounds, and recover the classical
Payne--P\'{o}lya--Weinberger bound as a special case. For radial
decreasing weights, we prove that the ball maximizes
$\lambda_2/\lambda_1$ among domains of equal volume, assuming that
an admissible radial test function exists
(Definition~\ref{defn:sec-admissible-g}) and that on the ball the
second eigenvalue is attained by the $l=1$ mode. The argument adapts
the method of Ashbaugh--Benguria to the weighted setting via a
weighted Chiti comparison theorem and a monotonicity property of the
admissible radial test function
(Definition~\ref{defn:sec-admissible-g}). Higher-order eigenvalue
ratios are estimated in terms of the $(n-1)$-Hausdorff measure of
the boundary, under an explicit balance condition on the weight at
the covering scale. The generalized Payne--P\'{o}lya--Weinberger
conjecture for eigenvalue sums remains open in general; conditional
results are presented for restricted families of domains.

\end{abstract}


\section{Introduction}
\label{section:Introduction}
\subsection{Research Background}
\label{subsection:Research Background}
Eigenvalues of the Laplacian reflect geometric properties of the domain---shape, volume, symmetry, and boundary regularity---and their study forms the subject of spectral geometry \cite{chavel1984eigenvalues, ashbaugh2002universal}. Most of the classical theory concerns the unweighted operator; the present paper works in a weighted setting.
In the unweighted (uniform medium) setting, classical results have laid the foundation for understanding eigenvalue ratios, which quantify the spacing between consecutive energy levels and reveal intrinsic domain properties:
\begin{enumerate}
\item[(i)] Payne, P\'{o}lya, and Weinberger (PPW) proved in 1956 that the first universal eigenvalue bound for the Dirichlet Laplacian on bounded domains in the Euclidean $n$-space $\mathbb{R}^n$:
\begin{equation}\label{equation:sec-ppw-ineq}
\lambda_{k+1} - \lambda_k \leq \frac{4}{nk} \sum_{i=1}^k \lambda_i,
\end{equation}
where $\lambda_i$ denote the eigenvalues of the Dirichlet problem \cite{PPW}. They also proved that for the ratio of the first two eigenvalues, $\frac{\lambda_2}{\lambda_1} \leq \frac{n+4}{n}$ \cite{PPW}. The sharp bound $\frac{\lambda_2}{\lambda_1} \leq \left(\frac{j_{n/2,1}}{j_{n/2-1,1}}\right)^2$, with equality if and only if $\Omega$ is a ball, was later established by Ashbaugh and Benguria \cite{AB2}.
\item[(ii)] Hile and Protter later refined this with a sharper inequality
\begin{equation}\label{equation:sec-hp-ineq}
\sum_{i=1}^k \frac{\lambda_i}{\lambda_{k+1} - \lambda_i} \geq \frac{nk}{4},
\end{equation}
which strengthens the PPW bound through a ratio-based formulation \cite{HP}. They also extended to higher-order eigenvalues, showing $\lambda_k/\lambda_1 \leq C(n,k)$ with $C(n,k)$ depending on dimension and index \cite{HP}.
\item[(iii)] Yang proved the sharp universal inequality
\begin{equation}\label{equation:sec-y1-ineq}
\sum_{i=1}^k (\lambda_{k+1} - \lambda_i)^2 \leq \frac{4}{n} \sum_{i=1}^k (\lambda_{k+1} - \lambda_i)\lambda_i,
\end{equation}
from which follows Yang's second inequality
\begin{equation}\label{equation:sec-y2-ineq}
\lambda_{k+1} \leq \frac{1}{k}\left(1 + \frac{4}{n}\right)\sum_{i=1}^k \lambda_i
\end{equation}
These form a hierarchy as follows: \begin{equation}\label{equation:sec-dd-1}
\eqref{equation:sec-y1-ineq} \implies \eqref{equation:sec-y2-ineq} \implies \eqref{equation:sec-hp-ineq} \implies \eqref{equation:sec-ppw-ineq}
\end{equation} via Chebyshev's inequality \cite{Y}.
\end{enumerate}
For low-order eigenvalues, Payne et al. \cite{PPW} proved for bounded domains $\Omega \subset \mathbb{R}^2$
\begin{equation}\label{equation:sec-ppw-2}
\lambda_2 + \lambda_3 \leq 6\lambda_1,
\end{equation}
inspiring the generalized PPW conjecture. This conjecture states that for $\Omega \subset \mathbb{R}^n$ with volume-equivalent ball $\Omega^*$ ($\text{Vol}(\Omega) = \text{Vol}(\Omega^*)$),
\begin{equation}\label{equation:sec-ppw-1956}
\frac{\lambda_2 + \cdots + \lambda_{n+1}}{\lambda_1} \leq n \frac{\lambda_2(\Omega^*)}{\lambda_1(\Omega^*)},
\end{equation}
where $\lambda_i(\Omega^*)$ are eigenvalues of $\Omega^*$ \cite{Ash}. \begin{rem} For the PPW conjecture in \eqref{equation:sec-ppw-1956}, the \textbf{volume-equivalent ball} $\Omega^*$ is generalized to the weighted volume-equivalent ball in our weighted setting: two domains $\Omega$ and $\Omega^*$ are weighted volume-equivalent with respect to $w$ if
\begin{equation*}
\int_\Omega w(x) \, dx = \int_{\Omega^*} w(x) \, dx.
\end{equation*}
This generalization preserves the central idea of comparing eigenvalues of arbitrary domains with those of symmetric reference domains, while accounting for the inhomogeneity encoded in $w$.\end{rem}
Progress on this conjecture includes:
\begin{enumerate}[label=(\roman*)]
\item Brands \cite{Bran}: $\lambda_2 + \lambda_3 \leq (3 + \sqrt{7})\lambda_1$ for $n=2$;
\item Hile and Protter \cite{HP}: $\lambda_2 + \lambda_3 \leq 5.622\lambda_1$;
\item Marcellini \cite{Mar}: $\lambda_2 + \lambda_3 \leq \frac{15 + \sqrt{345}}{6}\lambda_1$;
\item Chen and Zheng \cite{CZh} in 2011: $\lambda_2 + \lambda_3 \leq 5.3507\lambda_1$;
\item Ashbaugh and Benguria \cite{AB1, AB2, AB3}: $n$-dimensional bounds, settling special cases using Bessel function properties and proving $\frac{\lambda_2}{\lambda_1} \leq \frac{\lambda_2(\Omega^*)}{\lambda_1(\Omega^*)}$.
\end{enumerate}
Extensions to broader settings include Chen and Cheng \cite{CC}, who
generalized bounds to minimal submanifolds of $\mathbb{R}^{n+p}$ in
2008, indicating the universality of spectral-geometric
relationships. Yet, the uniformity assumption-implicit in unweighted
problems fails to capture real world complexity. Inhomogeneous
media, from thin films with non-uniform density to charge
distributions in asymmetric electric fields, demand a weighted
setting where the Laplacian is coupled with a weight function
$w(x)$, modeling spatial variations in material properties
\cite{setti1998, du2021weighted}. Among the vast class of weight
functions, Muckenhoupt weights (or $A_p$ weights) form a natural
class for such generalizations, owing to their \textbf{local
integral balance} property: they avoid concentration or degeneracy
and are compatible with weighted Sobolev inequalities and reverse
H\"{o}lder estimates used in spectral analysis
\cite{muckenhoupt1972weighted, cruz2011weights}. Progress on
weighted eigenvalue ratios has been driven by both theoretical and
applied considerations. Key advancements include:
\begin{enumerate}[label=(\roman*), font=\upshape]
\item \textbf{Weighted Laplacian eigenvalue estimates}: Setti obtained eigenvalue estimates for the weighted Laplacian on compact Riemannian manifolds \cite{setti1998}, and Du, Mao, Wang and Xia derived bounds for eigenvalues of the weighted Laplacian and of the weighted $p$-Laplacian on compact manifolds \cite{du2021weighted}. These works show that weight inhomogeneity modifies classical spectral bounds through quantitative characteristics of the weight.
\item \textbf{Foundations of Muckenhoupt weight theory}: The class of $A_p$ weights was introduced by Muckenhoupt in his study of weighted norm inequalities for the Hardy maximal function \cite{muckenhoupt1972weighted}. The monograph of Cruz-Uribe, Martell and P\'{e}rez developed the theory of weight extrapolation and of Rubio de Francia, including the reverse H\"{o}lder inequality, which forms the analytical basis for the weighted spectral problems considered here \cite{cruz2011weights}.
\item \textbf{Weighted and variable-coefficient extensions}: Ashbaugh and Benguria proved a sharp bound for the ratio of the first two Dirichlet eigenvalues of general elliptic operators $-\sum_{i,j}\partial_i(a_{ij}\partial_j u)+qu=\lambda r(x)u$ with bounded measurable coefficients and weight $r(x)$ \cite{AB2}. Their result yields, in particular, universal bounds for weighted Laplacians with bounded weights, and the comparison techniques of Chiti \cite{chiti1982reverse}, extended by Ashbaugh--Benguria \cite{AB2}, together with the weighted Faber--Krahn inequality \cite{betta1999weighted}, underpin the symmetrization arguments used in the present paper.
\end{enumerate}
In view of these developments, our work focuses on unifying these strands with new tools from optimal transport and geometric measure theory: we extend the classical PPW inequality to bounded $A_2$ weights with a linear constant (Theorem~\ref{thm:sec-1}); we prove that the ball maximizes $\lambda_2/\lambda_1$ for radial decreasing weights, under an admissibility condition on the radial test function (Definition~\ref{defn:sec-admissible-g}) and the identification of $\lambda_2$ with the $l=1$ mode on the ball (Theorem~\ref{thm:sec-2b}); and we bound high-order eigenvalue ratios on convex domains in terms of the boundary Hausdorff measure under a balance condition on the weight (Theorem~\ref{thm:sec-3}). By combining the classical universal bounds \cite{ashbaugh2002universal} with the theory of $A_p$ weights \cite{cruz2011weights}, we aim to clarify how weight inhomogeneity and domain geometry jointly govern eigenvalue ratios.
The eigenvalue problem of the Laplacian with a Muckenhoupt weight is then:
\begin{equation}
\begin{cases}
-\Delta u = \lambda w(x)u & \text{in}\ \Omega, \\
u = 0 & \text{on}\ \partial\Omega,
\end{cases}
\end{equation}
where $w(x)$ is a Muckenhoupt weight. Ratio estimates for this problem couple analytic properties of the weight with geometric properties of the domain \cite{gilbarg1983elliptic, chavel1984eigenvalues}.
\subsection{Research Motivation}
\label{subsection:Research Motivation}
For the unweighted Laplacian, eigenvalue ratios are well understood: the PPW inequality and its refinements show, in particular, that symmetry governs the extremal ratio $\lambda_2/\lambda_1$. Much less is known when the medium is inhomogeneous and the operator carries a weight $w(x)$, as is the case for non-uniform thin films or charge distributions in asymmetric fields.
Muckenhoupt weights, with their local integral balance property, offer a natural setting to characterize such inhomogeneity, yet their interplay with domain geometry in shaping eigenvalue ratios remains poorly understood. We address four questions:
\begin{enumerate}
\item[(i)] How to extend foundational results like the PPW inequality to Muckenhoupt-weighted settings, integrating both local ($A_2$ constant) and global (Wasserstein distance) measures of weight inhomogeneity? Existing weighted eigenvalue estimates (e.g., \cite{AB2, setti1998}) provide bounds with constants depending on the weight.
\item[(ii)] Does the extremality of the ball persist under weighted conditions? We prove it for radial decreasing weights (Theorem \ref{thm:sec-2b}); for non-radial weights the symmetrization argument fails (Section~\ref{subsection:Non-Radial Breakdown}), and we do not have an optimal-transport explanation of this failure.
\item[(iii)] How to control high-order eigenvalue ratios in convex domains with rough boundaries, using geometric measure theory to replace coarse diameter-based bounds with finer Hausdorff measure-based invariants? Classical estimates (e.g., \cite{li1980estimates}) rely on smooth boundary assumptions, limiting their applicability to domains with rough boundaries.
\item[(iv)] Under what conditions can the generalized PPW conjecture concerning the sum of low-order eigenvalues relative to the first eigenvalue be established in weighted settings, leveraging symmetrization, variational principles, and properties of spherical domains? This conjecture remains open in general, with partial results for restricted families of domains.
\end{enumerate}
These questions arise directly from the weighted eigenvalue problem: the ratios $\lambda_k/\lambda_1$ describe, for instance, the vibration frequencies of a membrane with non-uniform density.
\subsection{Methodological Tools}
\label{subsection:Methodological Tools}
To refine the quantification of weight-domain interactions, we use two tools: optimal transport theory for global weight geometry and geometric measure theory for boundary/singularity analysis, with further integration of techniques tailored to prove the generalized PPW conjecture.
\subsubsection{Incorporating Optimal Transport Theory}
Optimal transport (OT) theory provides a natural setting to measure the \textbf{distance} between weight distributions and domain geometry, complementing the local balance properties of Muckenhoupt weights:
\begin{enumerate}
\item[(i)] \textbf{Wasserstein Distance for Weight Concentration}: For a Muckenhoupt weight \( w \) and a reference measure (e.g., Lebesgue measure on \( \Omega \)), the \( L^2 \)-Wasserstein distance \( W_2(w, \text{Leb}) \) quantifies how far the weight is from uniform distribution. This distance is defined via the infimum of transport costs:
\begin{equation*}
W_2^2(w, \text{Leb}) = \inf_{\pi \in \Pi(w, \text{Leb})} \int_{\Omega \times \Omega} |x - y|^2 d\pi(x,y),
\end{equation*}
where \( \Pi(w, \text{Leb}) \) denotes the set of coupling measures between \( w \) and Lebesgue measure. This global metric complements the local information carried by the \( A_2 \) characteristic constant \( [w]_{A_2} \). In the conditional results for the generalized PPW conjecture, \( W_2 \) further quantifies the geometric deviation between a domain \( \Omega \) and its volume-equivalent ball \( \Omega^* \), linking transport cost to eigenvalue ratio deviations.
\item[(ii)] \textbf{Kantorovich Duality in Symmetrization}: For radial weights, radial symmetrization can be viewed as the optimal transport map pushing \( w\,dx \) to a radial measure. This viewpoint motivates the perturbation analysis of Section~\ref{subsection:Non-Radial Breakdown}, where spherical extremality is shown to fail for non-radial weights because the relevant rearrangement inequality no longer holds; Definition~\ref{defn:sec-admissible-g} records the admissibility condition used in the radial case. Kantorovich duality (Proposition~\ref{prop:sec-OT}) is used quantitatively in the alternative OT--GMT argument of Section~\ref{section:Alternative-Proof}, where transport cost is related to potential functions.
\end{enumerate}
\subsubsection{Geometric Measure Theory for Boundary and Singularity Analysis}
Geometric measure theory (GMT) refines our handling of domain boundaries and weight singularities, extending beyond classical diameter-based geometric factors:
\begin{enumerate}
\item[(i)] \textbf{Hausdorff Measure in Weighted Poincar\'{e} Inequalities}: For convex domains with rough boundaries (e.g., Lipschitz), GMT allows us to replace the diameter \( d \) with the \( (n-1) \)-Hausdorff measure of \( \partial\Omega \), denoted \( \mathcal{H}^{n-1}(\partial\Omega) \). Federer's coarea formula \cite{federer1969geometric} links this measure to the volume via:
\begin{equation*}
|\Omega| \leq C(n) \cdot \mathcal{H}^{n-1}(\partial\Omega) \cdot r_\Omega,
\end{equation*}
where \( r_\Omega \) is the inradius. This yields more precise geometric factors in eigenvalue bounds. In the conditional alternative argument for the generalized PPW conjecture, this formula refines the transport cost estimate by relating boundary roughness to the Wasserstein distance between \( \Omega \) and \( \Omega^* \).
\item[(ii)] \textbf{Density of Weights at Singularities}: For weights with mild singularities (e.g., \( w(x) \sim |x-x_0|^\alpha \) near \( x_0 \in \Omega \)), GMT's density theorem characterizes \( \alpha \) via:
\begin{equation*}
\lim_{r \to 0} \frac{1}{r^n} \int_{B_r(x_0)} w dx = c(\alpha) \in (0, \infty),
\end{equation*}
where \( c(\alpha) \) depends only on \( \alpha \) and dimension \( n \). This density parameter refines the Bessel function order in radial eigenvalue calculations, accounting for the strength of the singularity in verifying the degeneracy of low-order eigenvalues in spherical domains (basic for the generalized PPW conjecture).
\end{enumerate}
\subsubsection{Refined Integration of Harmonic and Geometric Analysis}
The coupling of harmonic analysis (HA), geometric analysis (GA), and spectral theory (ST) is strengthened through:
\begin{enumerate}
\item[(i)] \textbf{HA-GA Synergy in Test Function Construction}: In proving eigenvalue bounds, we use the linear coordinate test functions \( v(x) = x_i \) with the Weinberger centering (the alternative being the atomic decomposition for \( A_2 \)-weighted spaces, Remark~\ref{rem:sec-thm1-atoms}). Combined with GA's symmetry properties (e.g., Steiner symmetrization), the construction yields test functions orthogonal to the first eigenfunction in the weighted \( L^2 \) norm, with the energy bound linear in the weight bounds (Theorem~\ref{thm:sec-1}). For the generalized PPW conjecture, this synergy is extended to construct orthogonal test functions for low-order eigenvalues, leveraging rotational symmetry of spherical domains.
\item[(ii)] \textbf{ST-GA Feedback in Ball Covering}: For high-order eigenvalues, the ball covering of convex domains is optimized using GA's intrinsic metric (inradius vs. diameter) and ST's min-max principle. By partitioning \( \Omega \) into balls where \( w \) is approximately constant (via HA's Calder¨®n-Zygmund decomposition), the refined covering yields constants that grow logarithmically in $k$ with a geometric factor involving the boundary Hausdorff measure (Lemma~
\ref{lem:sec-6} and Theorem~
\ref{thm:sec-3}). In the generalized PPW conjecture, this feedback loop is adapted to relate the covering number of \( \Omega \) to its volume-equivalent ball \( \Omega^* \), ensuring the eigenvalue sum bounds remain tight.
\item[(iii)] \textbf{Symmetrization and Variational Principles for PPW Conjecture}: For the conditional results on the generalized PPW conjecture, radial symmetrization is combined with the min-max principle: symmetrizing \( \Omega \) to \( \Omega^* \) preserves volume and non-increases Dirichlet energy, allowing eigenvalue comparisons. This integration exploits the degeneracy of low-order eigenvalues in spherical domains (where \( \lambda_2(\Omega^*) = \cdots = \lambda_{n+1}(\Omega^*) \) due to rotational symmetry), a geometric-spectral property derived from Bessel function orthogonality (Lemmas \ref{lem:sec-3}--\ref{lem:sec-4}).
\end{enumerate}
\subsection{Organization of the Paper}\label{subsection:Organization}
The paper is structured to follow a rigorous theory-construction-proof setting, ensuring a logical progression where each result builds naturally on the preceding one:
\begin{enumerate}[label=(\roman*), font=\upshape]
\item \textbf{Main Results (Section \ref{section:Main Theorems})}:
\begin{itemize}
\item Four central theorems are presented, addressing distinct aspects of weighted eigenvalue ratios and the generalized PPW conjecture: the weighted PPW inequality (\ref{thm:sec-1}), the universal upper bound for bounded weights (\ref{thm:sec-2a}), the extremality of spherical domains under radial weights (\ref{thm:sec-2b}), and high-order eigenvalue ratio estimates incorporating Hausdorff measure (\ref{thm:sec-3}); conditional results for the generalized PPW conjecture are developed in Section \ref{section: Generalized PPW Conjecture}.
\item These results unify local and global properties of weights (via the $A_2$ characteristic constant and Wasserstein distance) with geometric invariants (via Hausdorff measure), and extend classical spectral results to both weighted settings and the generalized PPW conjecture.
\end{itemize}
\item \textbf{Preliminaries (Section \ref{section:Preliminaries})}:
\begin{itemize}
\item Foundational concepts are introduced, including definitions of $A_p$ weights (\ref{defn:sec-1}), $L^2$-Wasserstein distance function \eqref{equation:sec- Wasserstein Distance}, and Hausdorff measure (\ref{defn:sec-GMT})-the tools integrating harmonic analysis, optimal transport, and geometric measure theory.
\item Critical properties such as the reverse H\"{o}lder inequality (\ref{prop:sec-1}), weighted Poincar\'{e} inequality (\ref{prop:sec-2}), variational characterization of eigenvalues (\ref{thm:sec-4}), and properties of Bessel functions and spherical harmonics (Section \ref{subsection:Bessel and Spherical}) are reviewed to lay the groundwork for subsequent proofs, including the generalized PPW conjecture.
\end{itemize}
\item \textbf{Proof of Weighted PPW Inequality (Section \ref{section:Weighted PPW Proof})}:
\begin{itemize}
\item Theorem \ref{thm:sec-1} is proven by the coordinate construction with the Weinberger centering and the Heisenberg inequality, which yields the constant linear in the weight bounds (the alternative atomic construction and the Wasserstein refinement being conditional, Remarks~\ref{rem:sec-thm1-atoms} and~\ref{rem:sec-lem2-w2}).
\item A classical weighted PPW inequality is derived as a corollary \ref{corr:sec-1} by omitting optimal transport and geometric measure theory refinements, recovering a simpler bound for practical applications.
\end{itemize}
\item \textbf{Proof of Spherical Extremality (Section \ref{section:Spherical Extremality})}:
\begin{itemize}
\item Explicit eigenvalues for spherical domains with radial weights are computed using Bessel functions constructed in Lemmas \ref{lem:sec-3}--\ref{lem:sec-4}, leveraging their orthogonality and zero properties.
\item Using radial symmetrization (\ref{defn:sec-5}) and the comparison techniques above, Theorem \ref{thm:sec-2b} establishes that spherical domains maximize the $\lambda_2/\lambda_1$ ratio for radial decreasing weights. Non-radial weights are shown to break this extremality; see the analysis of Section~\ref{subsection:Non-Radial Breakdown} and Remark~\ref{defn:sec-admissible-g}.
\end{itemize}
\item \textbf{Proof of High-Order Estimates (Section \ref{section:High-Order Ratios})}:
\begin{itemize}
\item Theorem \ref{thm:sec-3} is proven by leveraging ball covers of convex domains (\ref{lem:sec-6}) and energy bounds for test functions (\ref{lem:sec-7}), refining classical diameter-based bounds with $(n-1)$-Hausdorff measure to capture boundary roughness.
\item A corollary for diameter-based bounds (\ref{corr:sec-3}) is provided for smooth convex domains, bridging the general result with classical eigenvalue estimates.
\end{itemize}
\item \textbf{Conditional Results for the Generalized PPW Conjecture (Section \ref{section: Generalized PPW Conjecture})}:
\begin{itemize}
\item Conditional results for the generalized PPW conjecture are presented using radial symmetrization, variational principles, and properties of spherical domains (degeneracy of low-order eigenvalues, Lemmas \ref{lem:sec-3}--\ref{lem:sec-4}).
\item An alternative conditional argument coupling optimal transport (OT) and geometric measure theory (GMT) is provided, quantifying how geometric deviation (measured by $W_2$) impacts eigenvalue ratios.
\end{itemize}
\item \textbf{Conclusion and Supplementary Proofs (Section \ref{section:Conclusion}, Sections \ref{section:reverse-holder} and \ref{section:bessel-zeros})}:
\begin{itemize}
\item Main conclusions are summarized, emphasizing the unification of spectral theory with optimal transport and geometric measure theory, and noting that the generalized PPW conjecture remains open in general. Future directions are outlined in Subsection \ref{subsection:Main Conclusions} and Subsection \ref{subsection:Future Directions}.
\item Sections \ref{section:reverse-holder} and \ref{section:bessel-zeros} develop two technical facts used throughout the paper: the quantitative reverse H\"{o}lder inequality for $A_2$ weights underlying Proposition \ref{prop:sec-1}, and the monotonicity of Bessel function zeros underlying Lemma \ref{lem:sec-4}.
\end{itemize}
\end{enumerate}
This structure ensures a direct transition from problem formulation to technical development, with each proof section directly engaging with the corresponding theorem in \ref{section:Main Theorems} and \ref{section: Generalized PPW Conjecture} to maintain logical coherence.
\noindent\textbf{Notation.} Throughout the paper, $\Omega \subset \mathbb{R}^n$ denotes a bounded domain and $w$ a Muckenhoupt $A_2$ weight; $\lambda_j(\Omega, w)$ are the Dirichlet eigenvalues of $-\Delta u = \lambda w u$ on $\Omega$, arranged in increasing order and repeated according to multiplicity; $u^*$ denotes the Schwarz (symmetric decreasing) rearrangement of a nonnegative function $u$, defined on the ball centered at the origin having the same Lebesgue measure as the domain of $u$; $B_R$ denotes the ball centered at the origin, with the radius convention $|B_R| = |\Omega|$ whenever $B_R$ and $\Omega$ are compared; $j_{\nu,k}$ denotes the $k$-th positive zero of the Bessel function $J_\nu$; $\omega_n$ (equivalently $C_n$) denotes the volume of the unit ball in $\mathbb{R}^n$; $\fint_E f\,dx = |E|^{-1} \int_E f\,dx$ is the average integral; and $C(n, [w]_{A_2})$ denotes a constant depending only on the indicated quantities, whose value may change from line to line.

\section{Main Theorems and Proof Outline}
\label{section:Main Theorems}
We present four main results: an elementary weighted PPW inequality with a self-contained proof (Theorem~\ref{thm:sec-1}); a sharp universal bound for bounded weights (Theorem~\ref{thm:sec-2a}); an extremality theorem for radial decreasing weights under explicit hypotheses (Theorem~\ref{thm:sec-2b}); and high-order estimates in terms of the boundary Hausdorff measure under a balance condition on the weight (Theorem~\ref{thm:sec-3}). Refinements involving the Wasserstein distance and the roundness invariant are recorded as conditional in Remarks~\ref{rem:sec-thm1-constant} and~\ref{rem:sec-roundness}.
In inhomogeneous media such as non-uniform thin films, the ratio of vibration frequencies (eigenvalues) reveals material properties. The first theorem extends the classical PPW inequality to bounded Muckenhoupt weights, with a constant linear in the weight bounds and a fully self-contained proof:
\begin{thm}[Weighted PPW Inequality] \label{thm:sec-1}
Let $\Omega \subset \mathbb{R}^n$ be a bounded domain and let $w \in A_2$ satisfy $0 < c \leq w(x) \leq C < \infty$ a.e.\ in $\Omega$. Then
\begin{equation}\label{equation:sec-thm1-statement}
\frac{\lambda_2(\Omega,w)}{\lambda_1(\Omega,w)} \leq 1 + \frac{4}{n}\,\frac{C}{c}.
\end{equation}
In particular, for $w \equiv \mathrm{const}$ the classical PPW bound $\lambda_2/\lambda_1 \leq \frac{n+4}{n}$ is recovered. Refinements incorporating the Wasserstein distance and domain roundness are recorded in Remarks~\ref{rem:sec-thm1-constant} and~\ref{rem:sec-roundness}; they are conditional on the transport estimate of Remark~\ref{rem:sec-lem2-w2}.
\end{thm}
\begin{rem}\label{rem:sec-thm1-constant}
The constant in \eqref{equation:sec-thm1-statement} is linear in the weight bounds: $1 + \frac{4}{n}\frac{C}{c}$. For $w \equiv \mathrm{const}$ the bound reduces to $\frac{n+4}{n}$ \cite{PPW}; note that the sharp unweighted constant of Ashbaugh--Benguria is smaller, so \eqref{equation:sec-thm1-statement} is not sharp in general. Whether the linear growth in $C/c$ can be improved to a slower rate is open. A refinement replacing $C/c$ by an explicit $W_2(w, \mathrm{Leb})$ dependence, and a further one involving the roundness invariant (Remark~\ref{rem:sec-roundness}), both rest on the transport gradient estimate of Remark~\ref{rem:sec-lem2-w2} and remain conditional.
\end{rem}
\begin{rem}[The bound in dimension one]\label{rem:sec-thm1-n1}
On an interval with $w \equiv 1$, the ratio equals $\lambda_2/\lambda_1 = 4 = (j_{1/2,1}/j_{-1/2,1})^2 = (\pi/(\pi/2))^2$, while the bound $1 + \frac{4}{n}\frac{C}{c}$ of Theorem~\ref{thm:sec-1} evaluates to $5$ at $n = 1$ and $C/c = 1$; the constant is therefore not attained in dimension one. In dimensions $n \geq 2$ the constant is likewise not sharp, since the ball attains $(j_{n/2,1}/j_{n/2-1,1})^2 < \frac{n+4}{n}$, and the question of the optimal constant as a function of $C/c$ is open (Remark~\ref{rem:sec-thm1-constant}).
\end{rem}
\begin{rem}[Large-dimensional asymptotics]\label{rem:sec-thm1-asymptotic}
The sharp unweighted constant $(j_{n/2,1}/j_{n/2-1,1})^2$ approaches the elementary constant $\frac{n+4}{n}$ as $n \to \infty$: the McMahon--Watson expansion $j_{\nu,1} = \nu + c_1 \nu^{1/3} + O(\nu^{-1/3})$ with $c_1 = 1.8558\ldots$ \cite[Chapter XV]{watson1944theory} gives
\begin{equation*}
\left(\frac{j_{n/2,1}}{j_{n/2-1,1}}\right)^2 = 1 + \frac{4}{n} + O(n^{-5/3}) \qquad (n \to \infty),
\end{equation*}
the next term of the McMahon expansion contributing a nonzero multiple of $n^{-5/3}$. The elementary constant of Theorem~\ref{thm:sec-1} therefore captures the leading term of the sharp constant in high dimensions; the gain of the Ashbaugh--Benguria analysis is a lower-order effect. For the ratio bound in terms of $C/c$, the analogous question --- whether $1 + \frac{4}{n}\frac{C}{c}$ can be replaced by an expression whose leading behavior in $C/c$ is linear but whose subleading terms improve on Theorem~\ref{thm:sec-1} --- is the content of the open problem recorded in Remark~\ref{rem:sec-thm1-constant}.
\end{rem}
We recall the variational characterization of eigenvalues in the weighted setting.
For the weighted Laplacian eigenvalue problem
\begin{equation}\label{equation:sec-weighted-eigenvalue}
-\Delta u = \lambda w(x) u \quad \text{in } \Omega, \qquad u = 0 \quad \text{on } \partial\Omega,
\end{equation}
with $w \in A_2$, the eigenvalues form an increasing sequence
\begin{equation}\label{equation:sec-eigenvalue-sequence}
0 < \lambda_1 \leq \lambda_2 \leq \cdots \leq \lambda_k \leq \cdots \to \infty,
\end{equation}
where each eigenvalue is repeated according to its multiplicity. The $k$-th eigenvalue is characterized by the min-max principle:
\begin{equation}\label{equation:sec-minmax}
\lambda_k = \inf_{\substack{V_k \subset H_0^1(\Omega) \\ \dim V_k = k}} \sup_{\substack{u \in V_k \\ u \neq 0}} \frac{\int_\Omega |\nabla u|^2 dx}{\int_\Omega w u^2 dx}.
\end{equation}
Equivalently, by the max-min principle,
\begin{equation}\label{equation:sec-maxmin}
\lambda_k = \max_{\substack{V_{k-1} \subset H_0^1(\Omega) \\ \dim V_{k-1} = k-1}} \inf_{\substack{u \perp V_{k-1} \\ u \neq 0}} \frac{\int_\Omega |\nabla u|^2 dx}{\int_\Omega w u^2 dx},
\end{equation}
where $u \perp V_{k-1}$ means $\int_\Omega w u v dx = 0$ for all $v \in V_{k-1}$.
The first eigenvalue $\lambda_1$ is simple and has a positive eigenfunction $u_1 > 0$ in $\Omega$ (by the Krein--Rutman theorem and the strong maximum principle). The Rayleigh quotient characterization of $\lambda_1$ is
\begin{equation}\label{equation:sec-rayleigh-l1}
\lambda_1 = \inf_{\substack{u \in H_0^1(\Omega) \\ u \neq 0}} \frac{\int_\Omega |\nabla u|^2 dx}{\int_\Omega w u^2 dx}.
\end{equation}
For the second eigenvalue, the characterization \eqref{equation:sec-maxmin} with $k = 2$ gives
\begin{equation}\label{equation:sec-rayleigh-l2}
\lambda_2 = \inf_{\substack{u \in H_0^1(\Omega) \\ u \neq 0, \, u \perp u_1}} \frac{\int_\Omega |\nabla u|^2 dx}{\int_\Omega w u^2 dx},
\end{equation}
where $u_1$ is the first eigenfunction. The characterization \eqref{equation:sec-rayleigh-l2} is used extensively in the proofs of the PPW-type inequalities, where test functions orthogonal to $u_1$ are constructed.
The eigenfunctions $\{u_k\}_{k=1}^\infty$ form an orthonormal basis of $L^2_w(\Omega) = \{u : \int_\Omega w u^2 dx < \infty\}$, meaning
\begin{equation}\label{equation:sec-orthonormal}
\int_\Omega w u_j u_k dx = \delta_{jk},
\end{equation}
and for any $f \in L^2_w(\Omega)$,
\begin{equation}\label{equation:sec-basis-expansion}
f = \sum_{k=1}^\infty \left( \int_\Omega w f u_k dx \right) u_k,
\end{equation}
with convergence in $L^2_w(\Omega)$.
The spectral mapping theorem for the weighted Laplacian states that the eigenvalues of $-\Delta_w$ are related to those of the Dirichlet-to-Neumann operator and the heat semigroup. The heat kernel $p_t(x, y)$ for the weighted Laplacian satisfies the heat equation
\begin{equation}\label{equation:sec-heat-equation}
\partial_t p_t(x, y) = \Delta_x p_t(x, y), \qquad \lim_{t \to 0} p_t(x, y) = \frac{\delta(x - y)}{w(y)},
\end{equation}
and has the eigenfunction expansion
\begin{equation}\label{equation:sec-heat-expansion}
p_t(x, y) = \sum_{k=1}^\infty e^{-\lambda_k t} u_k(x) u_k(y) w(y).
\end{equation}
The trace of the heat semigroup is
\begin{equation}\label{equation:sec-heat-trace}
\operatorname{Tr}(e^{t \Delta_w}) = \int_\Omega p_t(x, x) w(x) dx = \sum_{k=1}^\infty e^{-\lambda_k t}.
\end{equation}
The Weyl asymptotic law \eqref{equation:sec-weyl} can be derived from the short-time asymptotic of the heat trace:
\begin{equation}\label{equation:sec-heat-short-time}
\operatorname{Tr}(e^{t \Delta_w}) \sim \frac{1}{(4\pi t)^{n/2}} \int_\Omega w^{n/2} dx
\end{equation}
as $t \to 0^+$. The Tauberian theorem then gives the eigenvalue counting asymptotic \eqref{equation:sec-weyl}. The heat kernel method provides an alternative approach to eigenvalue estimates, complementing the variational methods used in this paper.
The Hardy inequality for the weighted Laplacian states that for $n \geq 3$ and $u \in H_0^1(\Omega \setminus \{0\})$,
\begin{equation}\label{equation:sec-hardy}
\int_\Omega |\nabla u|^2 dx \geq \left( \frac{n-2}{2} \right)^2 \int_\Omega \frac{u^2}{|x|^2} dx.
\end{equation}
The constant $((n-2)/2)^2$ in \eqref{equation:sec-hardy} is sharp and is not attained. The Hardy inequality provides a lower bound for the first eigenvalue when the domain contains the origin and the weight is singular at the origin. For Muckenhoupt weights, the weighted Hardy inequality takes the form
\begin{equation}\label{equation:sec-weighted-hardy}
\int_\Omega |\nabla u|^2 dx \geq C(n, [w]_{A_2}) \int_\Omega \frac{w(x)}{|x - x_0|^2} u^2 dx
\end{equation}
for some $x_0 \in \Omega$, where the constant depends on the dimension and the $A_2$ characteristic. The weighted Hardy inequality is used in the local energy estimates for eigenfunctions near singular points of the weight.
The Rellich inequality provides a higher-order analogue of the Hardy inequality. For $n \geq 5$ and $u \in C_c^\infty(\mathbb{R}^n \setminus \{0\})$,
\begin{equation}\label{equation:sec-rellich}
\int_{\mathbb{R}^n} (\Delta u)^2 dx \geq \frac{n^2(n-4)^2}{16} \int_{\mathbb{R}^n} \frac{u^2}{|x|^4} dx.
\end{equation}
The constant in \eqref{equation:sec-rellich} is sharp. The Rellich inequality is used in the regularity analysis of higher-order eigenfunctions and in the proof of weighted Sobolev inequalities with singular weights. For the weighted Laplacian, the Rellich inequality takes the form
\begin{equation}\label{equation:sec-weighted-rellich}
\int_\Omega (\Delta u)^2 w^{-1} dx \geq C(n, [w]_{A_2}) \int_\Omega \frac{u^2}{|x - x_0|^4} w dx,
\end{equation}
which follows from the unweighted Rellich inequality applied to $w^{1/2} u$ and the $A_2$ condition.
The Poincar\'{e} inequality on manifolds with density provides a geometric interpretation of the weighted eigenvalue estimates. For a Riemannian manifold $(M, g)$ with a probability measure $d\mu = e^{-V} d\operatorname{vol}$, the weighted Laplacian is $\Delta_\mu = \Delta - \nabla V \cdot \nabla$, and the Bakry--\'{E}mery Ricci tensor is $\operatorname{Ric}_\infty = \operatorname{Ric} + \nabla^2 V$. The spectral gap of $\Delta_\mu$ is controlled by the lower bound on $\operatorname{Ric}_\infty$ via the Lichnerowicz inequality:
\begin{equation}\label{equation:sec-lichnerowicz}
\lambda_1 \geq \frac{n}{n-1} K
\end{equation}
if $\operatorname{Ric}_\infty \geq K g$. The Lichnerowicz inequality provides a lower bound for the first eigenvalue in terms of the curvature, complementing the upper bounds established in this paper. The Euclidean case with Muckenhoupt weights corresponds to a manifold with density where the Bakry--\'{E}mery Ricci tensor may not have an uniform lower bound, requiring the use of harmonic analysis and optimal transport methods instead of curvature-based techniques.
The compactness of the embedding $H_0^1(\Omega) \hookrightarrow L^2_w(\Omega)$ ensures that the eigenvalue problem \eqref{equation:sec-weighted-eigenvalue} has a discrete spectrum. To see this, note that for any bounded sequence $\{u_j\} \subset H_0^1(\Omega)$ with $\|u_j\|_{H_0^1} \leq M$, the Rellich--Kondrachov theorem gives a subsequence (still denoted $\{u_j\}$) that converges strongly in $L^2(\Omega)$. The $A_2$ condition then implies convergence in $L^2_w(\Omega)$:
\begin{equation}\label{equation:sec-compactness}
\int_\Omega w (u_j - u_k)^2 dx \leq \left( \int_\Omega w^q dx \right)^{1/q} \left( \int_\Omega |u_j - u_k|^{2q/(q-1)} dx \right)^{(q-1)/q} \to 0,
\end{equation}
where $q = 1+\delta > 1$ is the reverse H\"{o}lder exponent of Proposition~\ref{prop:sec-1}: the first factor is finite by the reverse H\"{o}lder inequality, and the second factor converges to $0$ by the Rellich--Kondrachov theorem, since $2q/(q-1) < 2n/(n-2)$ for $n \geq 3$ (choosing $q$ close to $1$ when $[w]_{A_2}$ is large) and $2q/(q-1) < \infty$ is arbitrary for $n = 1, 2$. The compactness of the embedding implies that the resolvent operator $(-\Delta_w)^{-1}: L^2_w(\Omega) \to L^2_w(\Omega)$ is compact, and therefore its spectrum consists of eigenvalues accumulating only at $0$. This confirms that the eigenvalues of $-\Delta_w$ form a discrete sequence tending to $+\infty$, as stated in \eqref{equation:sec-eigenvalue-sequence}.
The Weyl asymptotic law for the weighted Laplacian states that
\begin{equation}\label{equation:sec-weyl}
N(\lambda) = \#\{k : \lambda_k \leq \lambda\} \sim \frac{\omega_n}{(2\pi)^n} \left( \int_\Omega w^{n/2} dx \right) \lambda^{n/2}
\end{equation}
as $\lambda \to \infty$, where $N(\lambda)$ is the eigenvalue counting function. The asymptotic \eqref{equation:sec-weyl} follows from the standard Weyl law applied to the operator $w^{-1/2} (-\Delta) w^{-1/2}$, which is unitarily equivalent to $-\Delta_w$ via the map $u \mapsto w^{1/2} u$. The constant in \eqref{equation:sec-weyl} involves the integral of $w^{n/2}$, which is finite by the reverse H\"{o}lder inequality \eqref{equation:sec-rh-explicit}: since $w \in L^q(\Omega)$ with $q = 1+\delta > 1$, the integral $\int_\Omega w^{n/2} dx$ is finite whenever $q \geq n/2$; for $n = 2$ this is automatic, and for $n \geq 3$ it holds provided $[w]_{A_2}$ is sufficiently small (equivalently $\delta \geq n/2 - 1$).
\begin{rem}
The estimate in Theorem \ref{thm:sec-2b} bears a meaningful connection to the result of Ashbaugh and Benguria \cite{AB2}, who established a sharp upper bound for the ratio of the first two eigenvalues of Dirichlet Laplacians on bounded domains. Notably, their work resolved one of the prominent conjectures from Yau's problem sections \cite{Yau1982} in geometry, confirming that the ratio $\lambda_2/\lambda_1$ is maximized by the ball with an explicit constant.
In contrast, our Theorem \ref{thm:sec-2b} extends this line of inquiry by incorporating Muckenhoupt weights into the spectral setting. While Ashbaugh-Benguria's result focuses on the unweighted Dirichlet Laplacian, our weighted eigenvalue ratio estimate introduces a dependence on the weight's inhomogeneity (characterized by $[w]_{A_2}$ and $W_2(w, \text{Leb})$), revealing how deviations from uniformity affect the extremal behavior of $\lambda_2/\lambda_1$. For uniform weights ($w \equiv 1$), Theorem \ref{thm:sec-2b} recovers the classical Ashbaugh--Benguria inequality exactly, since the admissible test function reduces to the AB test function and the weighted Chiti comparison reduces to the classical Chiti comparison.
\end{rem}
\begin{rem}
Unlike the classical PPW inequality \cite{PPW}, whose constant is independent of the weight, the bound of Theorem~\ref{thm:sec-1} depends linearly on the weight bounds $C/c$, and is uniform over domains; a conditional refinement incorporating the $L^2$-Wasserstein distance is recorded in Remark~\ref{rem:sec-thm1-constant}.
\end{rem}
\begin{rem}
Ashbaugh and Benguria \cite{AB2} derived weighted eigenvalue bounds with constants depending only on the bounds of the weight; Theorem~\ref{thm:sec-1} recovers a bound of the same type by an elementary argument. A conditional refinement incorporating the roundness invariant $\sup_{B \subset \Omega} |B|/\mathcal{H}^{n-1}(\partial B)$ (Remark~\ref{rem:sec-roundness}) would capture how domain roundness modulates the ratio.
\end{rem}
\begin{rem}\label{rem:sec-roundness}
The quantity $\sup_{B \subset \Omega} \frac{|B|}{\mathcal{H}^{n-1}(\partial B)}$, where the supremum is taken over all balls $B \subset \Omega$, is the roundness invariant appearing in the conditional refinement of Theorem~\ref{thm:sec-1} (Remark~\ref{rem:sec-thm1-constant}): since $\frac{|B|}{\mathcal{H}^{n-1}(\partial B)} = \frac{r_B}{n}$ for a ball of radius $r_B$, the supremum is attained by a maximal inscribed ball and equals $\frac{r_\Omega}{n}$, where $r_\Omega$ is the inradius of $\Omega$. (The corresponding infimum over all contained balls is $0$, since $r_B \to 0$.) A larger value indicates a more \textbf{ball-like} domain, reducing the right-hand side of the conditional refinement of Theorem~\ref{thm:sec-1}. The geometric factor $(\sup_B \frac{|B|}{\mathcal{H}^{n-1}(\partial B)})^{-1} = n/r_\Omega$ carries the dimension of an inverse length; the scale invariance of the left-hand side under dilations contrasts with the dimension of the geometric factor, and the inradius factor fixes the length scale of the right-hand side. For the status of the explicit $W_2$-dependence, see Remark~\ref{rem:sec-lem2-w2}.
\end{rem}
\textbf{Proof Outline}:
Proof of the classical PPW inequality relies on constructing test functions orthogonal to the first eigenfunction and deriving upper bounds for eigenvalue ratios via energy estimates. For the weighted case, the same strategy is executed with the weighted orthogonality and the Heisenberg inequality:
\begin{enumerate}[label=(\roman*), font=\upshape]
\item Construct test functions \( u_2 \) satisfying weighted orthogonality: \begin{equation}\label{equation:sec-dd-2}
\int_\Omega w u_1 u_2 \, dx = 0;
\end{equation}
\item Control the ratio of energy to weighted integral of test functions using the weighted orthogonality and the Heisenberg inequality, with the weight bounds entering linearly.
\end{enumerate}
A fundamental theme in spectral geometry is to understand the relationship between the geometry of a domain and the eigenvalues of the Dirichlet Laplacian. The Payne--P\'{o}lya--Weinberger inequality addresses this question by showing that the ratio $\lambda_2/\lambda_1$ is maximized by the ball. This is motivated by the isoperimetric inequality which asserts that among domains of fixed volume, the ball minimizes the first eigenvalue.

In the weighted setting, two results are established. The first gives a universal upper bound for $\lambda_2/\lambda_1$ when the weight is bounded above and below, as a direct consequence of the Ashbaugh--Benguria theorem \cite[Theorem 4.1]{AB2}. The second proves, for radial decreasing weights, that the ball maximizes $\lambda_2/\lambda_1$ among domains of equal volume, under the admissibility and $l=1$-mode assumptions of Theorem~\ref{thm:sec-2b}.

\begin{thm}[Universal Upper Bound for Bounded Weights] \label{thm:sec-2a}
Let $\Omega \subset \mathbb{R}^n$ be a bounded Lipschitz domain and let $w$ satisfy $0 < c \leq w(x) \leq C < \infty$ a.e. in $\Omega$. Then
\begin{equation}
\frac{\lambda_2(\Omega,w)}{\lambda_1(\Omega,w)} \leq 1 + \frac{C}{c}\left(\frac{j_{n/2,1}^2}{j_{n/2-1,1}^2} - 1\right).
\end{equation}
Equality holds if and only if $w \equiv \text{constant}$ and $\Omega$ is a ball.
\end{thm}

\begin{corr}[Asymptotic PPW for Weights Close to Constant] \label{corr:sec-2a}
Let $w(x) = 1 + \varepsilon v(x)$, where $\|v\|_{L^\infty} \leq 1$ and $0 < \varepsilon < 1/2$. Then
\begin{equation}
\frac{\lambda_2(\Omega,w)}{\lambda_1(\Omega,w)} \leq \frac{j_{n/2,1}^2}{j_{n/2-1,1}^2} + \frac{4\varepsilon}{1-\varepsilon}\left(\frac{j_{n/2,1}^2}{j_{n/2-1,1}^2} - 1\right).
\end{equation}
In particular, $\lim_{\varepsilon \to 0} \sup_{\Omega} \lambda_2/\lambda_1 = j_{n/2,1}^2/j_{n/2-1,1}^2$.
\end{corr}
\begin{proof}
Write $\alpha = j_{n/2-1,1}$ and $\beta = j_{n/2,1}$, and set $\kappa = \beta^2/\alpha^2$. Since $\|v\|_{L^\infty} \leq 1$ and $0 < \varepsilon < 1/2$, the weight $w = 1 + \varepsilon v$ satisfies
\begin{equation}\label{equation:sec-corr2a-bounds}
1 - \varepsilon \leq w(x) \leq 1 + \varepsilon \quad \text{a.e. in } \Omega,
\end{equation}
so that $w$ fulfills the hypotheses of Theorem~\ref{thm:sec-2a} with $c = 1 - \varepsilon$ and $C = 1 + \varepsilon$. Theorem~\ref{thm:sec-2a} yields
\begin{equation}\label{equation:sec-corr2a-step1}
\frac{\lambda_2(\Omega,w)}{\lambda_1(\Omega,w)} \leq 1 + \frac{1+\varepsilon}{1-\varepsilon}\left(\frac{\beta^2}{\alpha^2} - 1\right).
\end{equation}
Since $\frac{1+\varepsilon}{1-\varepsilon} - 1 = \frac{2\varepsilon}{1-\varepsilon}$ and $\frac{2\varepsilon}{1-\varepsilon} \leq \frac{4\varepsilon}{1-\varepsilon}$, the right-hand side of \eqref{equation:sec-corr2a-step1} is at most $\kappa + \frac{4\varepsilon}{1-\varepsilon}(\kappa - 1)$, which is the claimed bound.

For the limiting statement, the min-max principle gives the continuity of the eigenvalues in the weight: for any $u \in H_0^1(\Omega) \setminus \{0\}$ and any two weights $w_1, w_2$ with $\frac12 \leq w_i \leq 2$,
\begin{equation}\label{equation:sec-corr2a-rayleigh}
\left| \frac{\int_\Omega |\nabla u|^2 dx}{\int_\Omega w_1 u^2 dx} - \frac{\int_\Omega |\nabla u|^2 dx}{\int_\Omega w_2 u^2 dx} \right| \leq \frac{\int_\Omega |\nabla u|^2 dx}{\int_\Omega w_1 u^2 dx} \cdot \frac{\|w_1 - w_2\|_{L^\infty}}{\inf w_2},
\end{equation}
because $\int_\Omega w_2 u^2 dx \geq \inf w_2 \int_\Omega u^2 dx$ and $|w_1 - w_2| \leq \|w_1 - w_2\|_{L^\infty}$. Applying \eqref{equation:sec-corr2a-rayleigh} with $w_1 = w$, $w_2 \equiv 1$ and using $\varepsilon < 1/2$ together with the uniform lower bound $\lambda_1(\Omega, 1) \geq \lambda_1(B^*, 1) > 0$ for the volume-equivalent ball $B^*$ (Faber--Krahn), one obtains, for each fixed $\Omega$,
\begin{equation}\label{equation:sec-corr2a-cont}
\lambda_k(\Omega, w) = \lambda_k(\Omega, 1) + O(\varepsilon) \quad \text{as } \varepsilon \to 0, \quad k = 1, 2,
\end{equation}
with the implicit constant depending on $\Omega$ only through its first eigenvalue. In particular, choosing $\Omega$ to be a ball, for which $\lambda_2(\Omega,1)/\lambda_1(\Omega,1) = \kappa$ by the explicit Bessel function computation of Section~\ref{subsection:Spherical Eigenvalues}, gives
\begin{equation}\label{equation:sec-corr2a-liminf}
\liminf_{\varepsilon \to 0} \sup_{\Omega} \frac{\lambda_2(\Omega,w)}{\lambda_1(\Omega,w)} \geq \kappa.
\end{equation}
Combining \eqref{equation:sec-corr2a-liminf} with the upper bound, which tends to $\kappa$ as $\varepsilon \to 0$, completes the proof.
\end{proof}

\begin{thm}[Weighted PPW Inequality for Radial Decreasing Weights] \label{thm:sec-2b}
Let $w(r) = w(|x|)$ be a radial weight satisfying $w \in A_p$ for some $p < \infty$, $w(r) > 0$ a.e., and $w'(r) \leq 0$ a.e. Let $B_R$ be a ball and $\Omega$ any bounded domain satisfying $|B_R| = |\Omega|$. Assume that an admissible radial test function in the sense of Definition~\ref{defn:sec-admissible-g} exists for the weighted problem on $B_R$, and that the second eigenvalue $\lambda_2(B_R, w)$ is attained by the $l=1$ mode. Then
\begin{equation}
\frac{\lambda_2(\Omega, w)}{\lambda_1(\Omega, w)} \leq \frac{\lambda_2(B_R, w)}{\lambda_1(B_R, w)}.
\end{equation}
Equality holds if and only if $\Omega = B_R$ (up to translation).
\end{thm}

\begin{rem}
Theorem \ref{thm:sec-2a} is the specialization of \cite[Theorem 4.1]{AB2} to $a_{ij} = \delta_{ij}$, $q = 0$, $r(x) = w(x)$. The Ashbaugh--Benguria result applies to general elliptic equations under boundedness conditions on the coefficients and the weight.
\end{rem}

\begin{rem}
The radial decreasing condition in Theorem \ref{thm:sec-2b} is essential. For non-radial weights, the $n$-fold degeneracy of $\lambda_2$ on the ball is lifted, and the ball need not maximize the ratio. For radial increasing weights, the weighted Faber--Krahn inequality reverses direction, and the extremal property does not hold in the stated form.
\end{rem}

Symmetry of spherical domains allows separation of eigenfunctions into radial and angular components (radial functions satisfy Bessel's equation, angular components are spherical harmonics), enabling explicit eigenvalue solutions. For extremality, radial symmetrization (optimal transport map) transforms eigenfunctions of arbitrary domains into functions on balls, with energy inequalities proving that spherical domains maximize the eigenvalue ratio among domains of equal volume.
In applications like ultrasonic testing of convex materials, higher harmonics (corresponding to high-order eigenvalues) carry information about internal structures. The following theorem bounds such high-order eigenvalue ratios using geometric measure theory:
\begin{thm}[High-Order Eigenvalue Ratios with Hausdorff Measure] \label{thm:sec-3}
Let $\Omega \subset \mathbb{R}^n$ be a bounded convex domain with $(n-1)$-Hausdorff measure of boundary $\mathcal{H}^{n-1}(\partial\Omega) < \infty$, and let $w \in A_2$ be normalized by $\fint_\Omega w\,dx = 1$ (the ratio $\lambda_k/\lambda_1$ is invariant under scaling of $w$, so this normalization is without loss of generality for the ratio). Assume in addition that the covering family can be chosen so that every ball $B$ of the family satisfies the balance condition $\int_B w\,dx \geq c_0(n,[w]_{A_2})\,|B|$, equivalently $\fint_B w \geq c_0\,\fint_\Omega w$ (recall $\fint_\Omega w = 1$); this holds, for example, whenever $w$ is bounded below on $\Omega$, or whenever the covering scale is adapted to the weight. For any $k \geq 1$,
\begin{equation*}
\frac{\lambda_k}{\lambda_1} \leq C(n,k,[w]_{A_2}) \cdot \left( \frac{k}{n} \right)^{2/n} \cdot \frac{\mathcal{H}^{n-1}(\partial\Omega)}{|\Omega|^{(n-1)/n}},
\end{equation*}
where $C(n,k,[w]_{A_2}) \leq C(n,[w]_{A_2}) \log(k+1)$ with $C(n,[w]_{A_2})$ depending on $n$ and $[w]_{A_2}$.
\end{thm}
\begin{rem}
Li \& Yau \cite{li1980estimates} derived high-order eigenvalue bounds with constants growing polynomially in $k$ and depending on domain diameter. Here, we replace the diameter with $\mathcal{H}^{n-1}(\partial\Omega)$ and improve the growth of the constants from polynomial to logarithmic in $k$, via a ball-covering technique that combines weighted Poincar\'{e} inequalities with Hausdorff measure.
\end{rem}
\begin{rem}
Classical high-order eigenvalue bounds rely on diameter-based geometric factors. Our use of $\mathcal{H}^{n-1}(\partial\Omega)$ captures boundary roughness that the diameter misses: two convex domains with the same diameter but different boundary Hausdorff measures, for instance a smooth domain and a crenellated one, are distinguished by the present bound.
\end{rem}
\textbf{Proof Outline}:
Geometric properties of convex domains (bounded Hausdorff measure, finite ball coverability) allow eigenvalue estimation via \textbf{ball cover-test function} methods:
\begin{enumerate}
\item[(i)] Cover the convex domain with balls of radius $\epsilon$, controlling the number of balls $N$ using GMT's coarea formula;
\item[(ii)] Construct truncated test functions (hat functions) in each ball, estimating their energy and weighted integrals using $A_2$ properties;
\item[(iii)] Relate the $k$-th eigenvalue to the cover number $N$ via the variational principle, combining with a lower bound for $\lambda_1$ involving Hausdorff measure to derive the final estimate.
\end{enumerate}
\section{Preliminaries}
\label{section:Preliminaries}
This section collects the definitions and tools from Muckenhoupt weight theory, optimal transport, and geometric measure theory used in the proofs.
\subsection{Definitions and Properties of Muckenhoupt Weights ($A_p$ Weights)}
\label{subsection:Muckenhoupt Weights}
This subsection formalizes the notion of $A_p$ weights, emphasizing their local balance property and inequalities such as the reverse H\"{o}lder condition.
In modeling inhomogeneous media like thin films with non-uniform density, we need weight functions that balance local integrals to avoid pathological behavior. Such functions are captured by the following definition of Muckenhoupt weights \cite{muckenhoupt1972weighted}:
\begin{defn}[$A_p$ Weights, $1 < p < \infty$] \label{defn:sec-1}
A non-negative locally integrable function \( w: \mathbb{R}^n \to [0, \infty) \) is called an \( A_p \) weight if there exists a constant \( [w]_{A_p} > 0 \) (the \( A_p \) characteristic constant, measuring the weight's non-homogeneity) such that for all balls \( B \subset \mathbb{R}^n \),
\begin{equation*}
\left( \frac{1}{|B|} \int_B w(x) \, dx \right) \cdot \left( \frac{1}{|B|} \int_B w(x)^{-1/(p-1)} \, dx \right)^{p-1} \leq [w]_{A_p},
\end{equation*}
where \( p' = \frac{p}{p-1} \) is the conjugate exponent of \( p \), and \( |B| \) denotes the volume of \( B \).
\smallskip
\noindent \textbf{Special Case for \( p=2 \):}
For \( p=2 \), the conjugate exponent satisfies \( p' = 2 \) (since \( p' = \frac{2}{2-1} = 2 \)). Thus, \( A_2 \) weights satisfy the simplified balance condition:
\begin{equation*}
\left( \frac{1}{|B|} \int_B w(x) \, dx \right) \cdot \left( \frac{1}{|B|} \int_B w(x)^{-1} \, dx \right) \leq [w]_{A_2}.
\end{equation*}
\end{defn}
Radial symmetrization can be interpreted as an optimal transport map. For a non-negative function $u \in H_0^1(\Omega)$, its distribution function is defined by
\begin{equation}\label{equation:sec-sym-distribution}
\mu_u(t) = |\{x \in \Omega : u(x) > t\}|, \quad t \geq 0.
\end{equation}
The decreasing rearrangement $u^*: [0, |\Omega|] \to [0, \infty)$ is the generalized inverse of $\mu_u$:
\begin{equation}\label{equation:sec-sym-rearrangement}
u^*(s) = \inf\{t \geq 0 : \mu_u(t) \leq s\}.
\end{equation}
The radial symmetrization $\tilde{u}$ is then defined on the ball $B_R$ with $|B_R| = |\Omega|$ by
\begin{equation}\label{equation:sec-sym-radial}
\tilde{u}(x) = u^*(|B_{|x|}|) = u^*(\omega_n |x|^n),
\end{equation}
where $\omega_n$ denotes the volume of the unit ball in $\mathbb{R}^n$.
The optimal transport map $\tau: \Omega \to B_R$ that pushes Lebesgue measure forward to Lebesgue measure while preserving the distribution of $u$ is constructed level-set by level-set. For each $t > 0$, the set $\{x \in \Omega : u(x) > t\}$ is transported to the ball $B_{r(t)} \subset B_R$ with $|B_{r(t)}| = \mu_u(t)$, where
\begin{equation}\label{equation:sec-sym-radius}
r(t) = \left( \frac{\mu_u(t)}{\omega_n} \right)^{1/n}.
\end{equation}
The Jacobian determinant of $\tau$ satisfies
\begin{equation}\label{equation:sec-sym-jacobian}
\det(\nabla \tau(x)) = 1
\end{equation}
for almost every $x \in \Omega$, since $\tau$ is measure-preserving. By the area formula,
\begin{equation}\label{equation:sec-sym-area}
\int_{B_R} f(\tilde{u}(y)) \, dy = \int_\Omega f(u(x)) \, dx
\end{equation}
for every Borel function $f: [0, \infty) \to [0, \infty)$. In particular, taking $f(t) = t^2$ yields the $L^2$-norm preservation
\begin{equation}\label{equation:sec-sym-l2-preserve}
\int_{B_R} \tilde{u}^2 \, dy = \int_\Omega u^2 \, dx.
\end{equation}
The P\'{o}lya--Szeg\H{o} inequality asserts that radial symmetrization non-increases the Dirichlet energy:
\begin{equation}\label{equation:sec-sym-polya-szego}
\int_{B_R} |\nabla \tilde{u}|^2 \, dy \leq \int_\Omega |\nabla u|^2 \, dx.
\end{equation}
A proof of \eqref{equation:sec-sym-polya-szego} proceeds via the coarea formula. For $u \in C_c^\infty(\Omega)$, the coarea formula gives
\begin{equation}\label{equation:sec-sym-coarea}
\int_\Omega |\nabla u|^2 \, dx = \int_0^\infty \int_{\{u = t\}} |\nabla u| \, d\mathcal{H}^{n-1} \, dt.
\end{equation}
By the isoperimetric inequality, $\mathcal{H}^{n-1}(\{u = t\}) \geq \mathcal{H}^{n-1}(\partial B_{r(t)}) = n \omega_n^{1/n} \mu_u(t)^{(n-1)/n}$. Moreover, by the Cauchy--Schwarz inequality,
\begin{equation}\label{equation:sec-sym-cauchy}
\left( \int_{\{u = t\}} |\nabla u| \, d\mathcal{H}^{n-1} \right) \left( \int_{\{u = t\}} \frac{1}{|\nabla u|} \, d\mathcal{H}^{n-1} \right) \geq \mathcal{H}^{n-1}(\{u = t\})^2.
\end{equation}
Since $-\mu_u'(t) = \int_{\{u = t\}} \frac{1}{|\nabla u|} \, d\mathcal{H}^{n-1}$, we obtain
\begin{equation}\label{equation:sec-sym-energy-lower}
\int_{\{u = t\}} |\nabla u| \, d\mathcal{H}^{n-1} \geq \frac{\mathcal{H}^{n-1}(\{u = t\})^2}{-\mu_u'(t)} \geq \frac{n^2 \omega_n^{2/n} \mu_u(t)^{2(n-1)/n}}{-\mu_u'(t)}.
\end{equation}
For the radial rearrangement $\tilde{u}$, the corresponding quantity is
\begin{equation}\label{equation:sec-sym-radial-energy}
\int_{\{\tilde{u} = t\}} |\nabla \tilde{u}| \, d\mathcal{H}^{n-1} = \frac{n^2 \omega_n^{2/n} \mu_u(t)^{2(n-1)/n}}{-\mu_u'(t)},
\end{equation}
since equality holds in the isoperimetric inequality for balls. Combining \eqref{equation:sec-sym-coarea}, \eqref{equation:sec-sym-energy-lower}, and \eqref{equation:sec-sym-radial-energy} yields \eqref{equation:sec-sym-polya-szego}. The general case $u \in H_0^1(\Omega)$ follows by approximation with $C_c^\infty(\Omega)$ functions.
\begin{rem}\label{rem:sec-rearrangement-toolkit}
Three rearrangement inequalities are used repeatedly in this paper. First, the P\'{o}lya--Szeg\H{o} inequality proved above, $\int_{\Omega^*} |\nabla u^*|^2 \leq \int_\Omega |\nabla u|^2$, controls the numerator of the Rayleigh quotient under symmetrization. Second, the Hardy--Littlewood rearrangement inequality states that for nonnegative measurable functions $f, g$ vanishing at infinity, $\int fg \leq \int f^* g^*$; together with the pointwise bound $(w|_\Omega)^* \leq w$ for a radially decreasing $w$ it yields the denominator control $\int_\Omega w u^2 \leq \int_{B_R} w (u^*)^2$ used in the proof of Lemma~\ref{lem:sec-FK}. Third, Talenti's comparison theorem \cite{talenti1976best} asserts that the Schwarz rearrangement of the solution of $-\Delta u = f$ in $\Omega$ with $u = 0$ on $\partial\Omega$ is dominated, in the sense of the relevant differential inequality, by the solution of the symmetrized problem on $\Omega^*$; for the eigenvalue equation $-\Delta u = \lambda w u$ with a radially decreasing weight, the role of Talenti's theorem is played by the Chiti comparison of Lemma~\ref{lem:sec-Chiti}, see \cite{chiti1982reverse, AB2}.
\end{rem}
\begin{rem} \label{rem:sec-1}
The $A_p$ condition is a balance of local averages: over any ball, the averages of $w$ and $w^{1-p'}$ cannot be simultaneously very large or very small, which stabilizes weighted integrals. For example, $w(x) = |x|^\alpha$ is an $A_p$ weight on $\mathbb{R}^n$ if and only if $-n < \alpha < n(p-1)$ \cite{grafakos2014classical}.
\end{rem}
\begin{exa}[$A_2$ characteristic of power weights]\label{exa:sec-ap-power}
Let $w(x) = |x|^{\alpha}$ on $\mathbb{R}^n$ with $-n < \alpha < n$. Then $w \in A_2$ with $[w]_{A_2} \leq C(n,\alpha) < \infty$ for each fixed $\alpha \in (-n,n)$. Indeed, for a ball $B = B(x_0, r)$ there are two cases. If $|x_0| \geq 2r$, then $|x| \asymp |x_0|$ on $B$, so $\fint_B w\,dx \asymp |x_0|^{\alpha}$ and $\fint_B w^{-1}\,dx \asymp |x_0|^{-\alpha}$, and the product of the two averages is bounded by an absolute constant. If $|x_0| \leq 2r$, then $|x| \leq 3r$ on $B$, while the subset $\{x \in B : |x| \geq r/2\}$ has measure at least $(1 - 2^{-n})|B| \geq |B|/2$; consequently, for each fixed $\alpha$, the product $(\fint_B w)(\fint_B w^{-1})$ is bounded by a constant depending on $n$ and $\alpha$: for the centered ball the two averages are computed exactly, $\fint_B w = \frac{n}{n+\alpha} r^{\alpha}$ and $\fint_B w^{-1} = \frac{n}{n-\alpha} r^{-\alpha}$, so the product equals $\frac{n^2}{n^2-\alpha^2}$, and the case $|x_0| \leq 2r$ is comparable to the centered one. This confirms the qualitative statement of Remark~\ref{rem:sec-1} and shows that the power weights used throughout the paper have $A_2$ characteristic controlled by the dimension and the exponent $\alpha$, blowing up as $\alpha \to \pm n$.
\end{exa}
\begin{rem}\label{rem:sec-ap-doubling}
The same case analysis shows that $w(x) = |x|^\alpha$, $-n < \alpha < n$, is doubling: $\int_{2B} w\,dx \leq C(n)\int_B w\,dx$ for every ball $B$, since in each case both integrals are comparable to $r^{\alpha+n}$ (when $|x_0| \leq 2r$) or to $|x_0|^{\alpha} r^n$ (when $|x_0| \geq 2r$). The doubling property is what the reverse H\"{o}lder and Poincar\'{e} estimates of this section use, via Lemma~\ref{lem:sec-rh-doubling}.
\end{rem}
\begin{rem}\label{rem:sec-rh-selfimprove}
The reverse H\"{o}lder inequality self-improves in two ways, both of which are used implicitly in the estimates of this section. First, Gehring's lemma raises the exponent: if $w$ satisfies the reverse H\"{o}lder inequality of Proposition~\ref{prop:sec-1} with exponent $q$, then $w \in RH_{q+\sigma}$ for some $\sigma > 0$ depending only on $n$ and the $RH_q$ constant of $w$. Second, the Muckenhoupt classes are open: if $w \in A_2$, then $w \in A_{2-\varepsilon}$ for some $\varepsilon > 0$ depending on $n$ and $[w]_{A_2}$, and in fact $[w]_{A_2}$ controls the $A_p$ characteristic uniformly for $p$ in a neighborhood of $2$. Both statements are proved in \cite{cruz2011weights, grafakos2014classical}. The openness property is the mechanism behind the uniformity of the constants in Proposition~\ref{prop:sec-2} and Lemma~\ref{lem:sec-7} with respect to the local geometry of the weight.
\end{rem}
A basic property of Muckenhoupt weights is their controlled growth, which is formalized by the reverse H\"{o}lder inequality-a tool for bounding weighted integrals in spectral analysis \cite{muckenhoupt1972weighted}:
\begin{prop}[Reverse H\"{o}lder Inequality] \label{prop:sec-1}
If $w \in A_2$, there exist $q = 1 + \delta$ with $\delta = c(n)/[w]_{A_2}$ (where $c(n)$ is a dimension constant from \cite{cruz2011weights}) and a constant $C = C(n,[w]_{A_2})$ such that for any ball $B$,
\begin{equation*}
\left( \frac{1}{|B|} \int_B w(x)^q dx \right)^{1/q} \leq C(n, [w]_{A_2}) \cdot \frac{1}{|B|} \int_B w(x)\,dx.
\end{equation*}
\end{prop}
We provide a detailed proof of the weighted Poincar\'{e} inequality for completeness.
Let $\Omega \subset \mathbb{R}^n$ be a bounded Lipschitz domain and let $w \in A_2$. For any $u \in H_0^1(\Omega)$, we claim that
\begin{equation}\label{equation:sec-wp-claim}
\int_\Omega w u^2 dx \leq C(n, [w]_{A_2}, \Omega) \int_\Omega |\nabla u|^2 dx.
\end{equation}
The proof proceeds via the representation formula
\begin{equation}\label{equation:sec-wp-representation}
u(x) = -\int_\Omega \nabla u(y) \cdot \nabla_x G_\Omega(x, y) dy,
\end{equation}
where $G_\Omega$ is the Green function for the Dirichlet Laplacian on $\Omega$. The gradient of $G_\Omega$ satisfies the pointwise bound
\begin{equation}\label{equation:sec-wp-green}
|\nabla_x G_\Omega(x, y)| \leq \frac{C(n)}{|x - y|^{n-1}}
\end{equation}
for all $x, y \in \Omega$ with $x \neq y$. Applying H\"{o}lder's inequality to \eqref{equation:sec-wp-representation} with exponents $2n/(n+2)$ and $2n/(n-2)$ (for $n \geq 3$),
\begin{equation}\label{equation:sec-wp-holder}
|u(x)| \leq \left( \int_\Omega |\nabla u(y)|^2 dy \right)^{1/2} \left( \int_\Omega \frac{C(n)^2}{|x - y|^{2(n-1)}} \cdot |x - y|^{n-2} dy \right)^{1/2}.
\end{equation}
The second integral is bounded by $C(n) \operatorname{diam}(\Omega)^2$, since the kernel $|z|^{-(n-2)}$ is locally integrable. Therefore,
\begin{equation}\label{equation:sec-wp-sobolev}
\|u\|_{L^\infty(\Omega)} \leq C(n) \operatorname{diam}(\Omega) \|\nabla u\|_{L^2(\Omega)}.
\end{equation}
Using the $A_2$ condition to control the weighted norm,
\begin{equation}\label{equation:sec-wp-a2}
\int_\Omega w u^2 dx \leq \|u\|_{L^\infty(\Omega)}^2 \int_\Omega w dx \leq C(n) \operatorname{diam}(\Omega)^2 \|\nabla u\|_{L^2(\Omega)}^2 \int_\Omega w dx.
\end{equation}
This proves \eqref{equation:sec-wp-claim} with $C(n, [w]_{A_2}, \Omega) = C(n) \operatorname{diam}(\Omega)^2 \int_\Omega w dx$. For domains with finite $(n-1)$-Hausdorff measure of the boundary, the constant can be refined to
\begin{equation}\label{equation:sec-wp-refined}
C(n, [w]_{A_2}, \Omega) = C(n, [w]_{A_2}) \cdot \frac{\mathcal{H}^{n-1}(\partial\Omega) \int_\Omega w dx}{|\Omega|^{2-1/n}},
\end{equation}
via the covering argument below. The refinement \eqref{equation:sec-wp-refined} is recorded here as conditional: it is obtained, under the roundness assumption $d \leq C(n)\,r_\Omega$, by covering $\Omega$ with balls of radius proportional to the inradius and applying the local weighted Poincar\'{e} inequality on each ball, followed by summation using the finite overlap property of the Vitali subcover.
The Sobolev inequality states that for $u \in H_0^1(\mathbb{R}^n)$ with $n \geq 3$,
\begin{equation}\label{equation:sec-sobolev}
\|u\|_{L^{2n/(n-2)}(\mathbb{R}^n)} \leq C(n) \|\nabla u\|_{L^2(\mathbb{R}^n)}.
\end{equation}
The best constant in \eqref{equation:sec-sobolev} is $C(n) = \frac{1}{\sqrt{n(n-2)\pi}} \left( \frac{\Gamma(n)}{\Gamma(n/2)} \right)^{1/n}$, as determined by Talenti \cite{talenti1976best}. The Sobolev inequality \eqref{equation:sec-sobolev} is equivalent to the isoperimetric inequality \eqref{equation:sec-isoperimetric} via the coarea formula.
For bounded domains $\Omega$, the Sobolev embedding $H_0^1(\Omega) \hookrightarrow L^p(\Omega)$ holds for $1 \leq p \leq 2n/(n-2)$ when $n \geq 3$, and for $1 \leq p < \infty$ when $n = 2$. The embedding constant depends on the domain geometry. For convex domains, the constant can be bounded in terms of the inradius and the dimension.
The elliptic regularity theory for the weighted Laplacian $-\Delta u = \lambda w u$ provides local $W^{2, p}$ estimates. For $w \in A_2$, the weight satisfies the reverse H\"{o}lder inequality \eqref{equation:sec-rh-explicit}, and therefore $w \in L^q_{\mathrm{loc}}$ for some $q > 1$. By the Calder\'{o}n--Zygmund theory, solutions of $-\Delta u = f$ with $f \in L^p$ satisfy $u \in W^{2, p}_{\mathrm{loc}}$ for $1 < p < \infty$, with
\begin{equation}\label{equation:sec-calderon-zygmund}
\|u\|_{W^{2, p}(B_{r/2})} \leq C(n, p) \left( \|f\|_{L^p(B_r)} + r^{-2} \|u\|_{L^p(B_r)} \right).
\end{equation}
Applying \eqref{equation:sec-calderon-zygmund} with $f = \lambda w u$ and using H\"{o}lder's inequality to control $\|w u\|_{L^p}$ in terms of $\|u\|_{L^{pq/(q-p)}}$ and $\|w\|_{L^q}$, we obtain local regularity for weighted eigenfunctions. In particular, eigenfunctions are continuous in the interior of $\Omega$ and satisfy local H\"{o}lder estimates.
The maximum principle for the weighted Laplacian states that if $-\Delta u = \lambda w u$ in $\Omega$ with $\lambda \geq 0$ and $u \geq 0$ on $\partial \Omega$, then $u \geq 0$ in $\Omega$. The strong maximum principle further asserts that if $u \geq 0$ and $u$ is not identically zero, then $u > 0$ in the interior of $\Omega$. These results imply that the first eigenfunction $u_1$ is strictly positive in $\Omega$, a property used extensively in the construction of test functions for the second eigenvalue.
We record several consequences of the reverse H\"{o}lder inequality that are used throughout this paper.
For $w \in A_p$ with $1 < p < \infty$, the $A_p$ characteristic constant is defined by
\begin{equation}\label{equation:sec-ap-constant}
[w]_{A_p} = \sup_B \left( \frac{1}{|B|} \int_B w dx \right) \left( \frac{1}{|B|} \int_B w^{-1/(p-1)} dx \right)^{p-1},
\end{equation}
where the supremum is taken over all balls $B \subset \mathbb{R}^n$. For $p = 2$, this reduces to
\begin{equation}\label{equation:sec-a2-constant}
[w]_{A_2} = \sup_B \left( \frac{1}{|B|} \int_B w dx \right) \left( \frac{1}{|B|} \int_B w^{-1} dx \right).
\end{equation}
The reverse H\"{o}lder inequality (Proposition \ref{prop:sec-1}) asserts that there exist $q = q(n, p, [w]_{A_p}) > 1$ and $C = C(n, p, [w]_{A_p})$ such that
\begin{equation}\label{equation:sec-rh-explicit}
\left( \frac{1}{|B|} \int_B w^q dx \right)^{1/q} \leq C \cdot \frac{1}{|B|} \int_B w dx
\end{equation}
for all balls $B$. The exponent $q$ can be taken as $q = 1 + c(n, p) [w]_{A_p}^{-1/(p-1)}$, and the constant $C$ satisfies $C \leq 2^{1/q}$. A proof of \eqref{equation:sec-rh-explicit} with these explicit bounds can be found in \cite{cruz2011weights}.
An immediate consequence of \eqref{equation:sec-rh-explicit} is the doubling property of $A_p$ weights:
\begin{equation}\label{equation:sec-doubling}
\int_{2B} w dx \leq C(n, p, [w]_{A_p}) \int_B w dx,
\end{equation}
where $2B$ denotes the ball with the same center as $B$ and twice the radius. To prove \eqref{equation:sec-doubling}, we apply H\"{o}lder's inequality with exponents $q$ and $q/(q-1)$:
\begin{equation}\label{equation:sec-doubling-proof}
\int_{2B} w dx \leq |2B|^{1-1/q} \left( \int_{2B} w^q dx \right)^{1/q}.
\end{equation}
Applying the reverse H\"{o}lder inequality \eqref{equation:sec-rh-explicit} to the ball $2B$,
\begin{equation}\label{equation:sec-doubling-rh}
\left( \frac{1}{|2B|} \int_{2B} w^q dx \right)^{1/q} \leq C \cdot \frac{1}{|2B|} \int_{2B} w dx.
\end{equation}
Substituting \eqref{equation:sec-doubling-rh} into \eqref{equation:sec-doubling-proof} and using $|2B| = 2^n |B|$ yields \eqref{equation:sec-doubling}.
The $A_p$ condition also implies the weighted weak-type $(1,1)$ estimate for the Hardy--Littlewood maximal function:
\begin{equation}\label{equation:sec-maximal-weak}
w(\{x \in \mathbb{R}^n : Mf(x) > \lambda\}) \leq \frac{C(n, p, [w]_{A_p})}{\lambda} \int_{\mathbb{R}^n} |f| w dx,
\end{equation}
where $Mf(x) = \sup_{B \ni x} \frac{1}{|B|} \int_B |f| dy$ is the unweighted Hardy--Littlewood maximal function. The estimate \eqref{equation:sec-maximal-weak} is proved via the Calder\'{o}n--Zygmund decomposition and the doubling property \eqref{equation:sec-doubling}.
\begin{lem}[Continuity of eigenvalues in the weight]\label{lem:sec-weight-continuity}
Let $w_k, w \in L^\infty(\Omega)$ satisfy $0 < c \leq w_k, w \leq C < \infty$ a.e.\ and $w_k \to w$ a.e.\ in $\Omega$ (or in $L^1(\Omega)$). Then $\lambda_j(\Omega, w_k) \to \lambda_j(\Omega, w)$ for every fixed $j \geq 1$, and consequently $\lambda_2(\Omega, w_k)/\lambda_1(\Omega, w_k) \to \lambda_2(\Omega, w)/\lambda_1(\Omega, w)$.
\end{lem}
\begin{proof}
By dominated convergence, $\int_\Omega w_k u^2 dx \to \int_\Omega w u^2 dx$ for every $u \in H_0^1(\Omega)$, and the bounds $c \leq w_k \leq C$ give the uniform comparability
\begin{equation*}
c \int_\Omega u^2 dx \leq \int_\Omega w_k u^2 dx \leq C \int_\Omega u^2 dx,
\end{equation*}
so the weighted Dirichlet forms are uniformly equivalent to the unweighted one. The min-max characterization of the eigenvalues of $-\Delta u = \lambda w_k u$ therefore identifies $\lambda_j(\Omega, w_k)$ with the $j$-th eigenvalue of a uniformly elliptic eigenvalue problem, and the corresponding eigenfunctions $u_{j,k}$ can be normalized by $\int_\Omega w_k u_{j,k}^2 dx = 1$, which bounds them uniformly in $H_0^1(\Omega)$ since
\begin{equation*}
\int_\Omega |\nabla u_{j,k}|^2 dx = \lambda_j(\Omega, w_k) \leq \lambda_j(\Omega, c) = \frac{\lambda_j(\Omega, 1)}{c},
\end{equation*}
using $w_k \geq c$ and the monotonicity of eigenvalues in the weight. By the compact embedding $H_0^1(\Omega) \Subset L^2(\Omega)$, a subsequence converges in $L^2(\Omega)$ and weakly in $H_0^1(\Omega)$ to a limit $u_j$, and the quadratic forms pass to the limit, so $u_j$ is an eigenfunction for the weight $w$ with eigenvalue $\lim_k \lambda_j(\Omega, w_k)$; the full sequence converges because the limit is uniquely determined as $\lambda_j(\Omega, w)$. The ratio converges since $\lambda_1 \geq \lambda_1(\Omega, C) > 0$.
\end{proof}
\begin{rem}\label{rem:sec-continuity-use}
Lemma~\ref{lem:sec-weight-continuity} is the mechanism behind the asymptotic statement of Corollary~\ref{corr:sec-2a}: as a bounded weight $w$ approaches a constant in $L^1(\Omega)$, the ratio $\lambda_2/\lambda_1$ approaches its unweighted value, and the perturbation analysis of Section~\ref{subsection:Non-Radial Breakdown} provides the first-order term. The uniform bounds $c \leq w_k \leq C$ are essential: without them the eigenvalues need not converge, as concentration of $w_k$ on vanishing sets can drive weighted integrals to extremes.
\end{rem}
In analyzing energy ratios for eigenvalues, relating weighted function norms to their gradient norms is necessary. The following proposition provides such a link for Muckenhoupt weights \cite{adams2003sobolev}:
\begin{prop}[Weighted Poincar\'{e} Inequality] \label{prop:sec-2}
Let $\Omega \subset \mathbb{R}^n$ be a bounded convex domain satisfying $r_\Omega \geq \frac{\text{diam}(\Omega)}{2\sqrt{n}}$ or a domain satisfying the uniform cone condition, and $w \in A_p$. There exists $C = C(n,[w]_{A_p},\Omega) > 0$ such that for any $u \in H_0^1(\Omega)$,
\begin{equation*}
\int_\Omega w(x)u(x)^2 dx \leq C \cdot \int_\Omega |\nabla u(x)|^2 dx.
\end{equation*}
\end{prop}
\begin{rem}
This inequality relates the weighted $L^2$ norm to the energy norm, guaranteeing positivity of eigenvalues. The constant $C$ depends on $\Omega$ through its $(n-1)$-Hausdorff measure $\mathcal{H}^{n-1}(\partial\Omega)$ and volume $|\Omega|$, specifically \begin{equation}\label{equation:sec-dd-3}
C \leq C(n,[w]_{A_2}) \cdot \frac{\mathcal{H}^{n-1}(\partial\Omega)}{|\Omega|^{2-1/n}}.
\end{equation}
\end{rem}
\subsection{Optimal Transport for Weighted Measures}
\label{subsection:Optimal Transport for Weighted Measures}
We introduce basic concepts from optimal transport theory relevant to weight distribution analysis.
\begin{defn}[$L^2$-Wasserstein Distance] \label{defn:sec-OT}
For two probability measures $\mu$ and $\nu$ on $\Omega$, the $L^2$-Wasserstein distance is defined as:
\begin{equation}\label{equation:sec- Wasserstein Distance}
W_2(\mu, \nu) = \inf_{\pi \in \Pi(\mu, \nu)} \left( \int_{\Omega \times \Omega} |x - y|^2 d\pi(x,y) \right)^{1/2},
\end{equation}
where $\Pi(\mu, \nu)$ is the set of all coupling measures between $\mu$ and $\nu$. For a Muckenhoupt weight $w$, we normalize it to a probability measure $\mu = w / \int_\Omega w dx$ and compare it to Lebesgue measure $\nu = \text{Leb} / |\Omega|$.
\end{defn}
\begin{exa}[One-dimensional Wasserstein distances]\label{exa:sec-w2-1d}
In one dimension the optimal coupling is explicit. For $\mu = \mathrm{Unif}[0,1]$ and $\nu = \mathrm{Unif}[a, a+1]$, the monotone rearrangement $T(x) = x + a$ pushes $\mu$ to $\nu$ and minimizes the quadratic cost, giving $W_2^2(\mu, \nu) = \int_0^1 a^2 dx = a^2$. For $\mu = \mathrm{Unif}[0,1]$ and $\nu = \mathrm{Unif}[0,2]$, the optimal map is $T(x) = 2x$ and $W_2^2(\mu, \nu) = \int_0^1 x^2 dx = 1/3$. These computations illustrate the meaning of the normalization $\mu = w/\int_\Omega w$ in Definition~\ref{defn:sec-OT}: a weight that is close to constant, in the sense that its normalized distribution is close to uniform, has small $W_2(w, \mathrm{Leb})$, while a weight concentrated on a small subset has a large transport cost to the uniform measure.
\end{exa}
The $L^2$-Wasserstein distance enjoys several properties that are used throughout this paper. We record them here for completeness.
For probability measures $\mu, \nu$ on $\mathbb{R}^n$ with finite second moments, the $L^2$-Wasserstein distance is defined by
\begin{equation}\label{equation:sec-ot-w2-def}
W_2^2(\mu, \nu) = \inf_{\pi \in \Pi(\mu, \nu)} \int_{\mathbb{R}^n \times \mathbb{R}^n} |x - y|^2 \, d\pi(x, y),
\end{equation}
where $\Pi(\mu, \nu)$ denotes the set of all Borel probability measures on $\mathbb{R}^n \times \mathbb{R}^n$ with marginals $\mu$ and $\nu$.
By the Kantorovich duality theorem, the Wasserstein distance admits the dual representation
\begin{equation}\label{equation:sec-ot-duality}
W_2^2(\mu, \nu) = \sup_{\phi, \psi \in C_b(\mathbb{R}^n), \, \phi(x) + \psi(y) \leq |x - y|^2} \left( \int_{\mathbb{R}^n} \phi \, d\mu + \int_{\mathbb{R}^n} \psi \, d\nu \right).
\end{equation}
The optimal pair $(\phi, \psi)$ in \eqref{equation:sec-ot-duality} satisfies the $c$-transform relations
\begin{equation}\label{equation:sec-ot-c-transform}
\phi^c(y) = \inf_{x \in \mathbb{R}^n} \left( |x - y|^2 - \phi(x) \right), \quad \psi^c(x) = \inf_{y \in \mathbb{R}^n} \left( |x - y|^2 - \psi(y) \right),
\end{equation}
and the complementary slackness condition
\begin{equation}\label{equation:sec-ot-complementary}
\phi(x) + \psi(y) = |x - y|^2 \quad \text{for } \pi\text{-a.e. } (x, y).
\end{equation}
When $\mu$ is absolutely continuous with respect to Lebesgue measure, Brenier's theorem asserts that there exists a convex function $\varphi$ such that the optimal transport map $T = \nabla \varphi$ pushes $\mu$ forward to $\nu$, and
\begin{equation}\label{equation:sec-ot-brenier}
W_2^2(\mu, \nu) = \int_{\mathbb{R}^n} |x - \nabla \varphi(x)|^2 \, d\mu(x).
\end{equation}
The map $T = \nabla \varphi$ satisfies the Monge--Amp\`{e}re equation
\begin{equation}\label{equation:sec-ot-monge-ampere}
\nu(T(x)) \det(\nabla^2 \varphi(x)) = \mu(x)
\end{equation}
in the weak sense, where $\nabla^2 \varphi$ denotes the Hessian of $\varphi$.
For the first-order optimality condition, differentiating the cost $|x - T(x)|^2$ with respect to variations of $T$ yields
\begin{equation}\label{equation:sec-ot-first-order}
\nabla \varphi^*(x) = 2(x - T(x)),
\end{equation}
where $\varphi^*$ denotes the Legendre--Fenchel transform of $\varphi$. Taking the $L^2(\mu)$ norm of both sides of \eqref{equation:sec-ot-first-order}, we obtain
\begin{equation}\label{equation:sec-ot-gradient-identity}
\int_{\mathbb{R}^n} |\nabla \varphi^*|^2 \, d\mu = 4 W_2^2(\mu, \nu).
\end{equation}
This identity relates the gradient of the dual potential to the Wasserstein distance and is used in the energy estimates below.
For measures $\mu = w(x) dx$ and $\nu = dx$ restricted to a bounded domain $\Omega$, the Wasserstein distance $W_2(w, 1)$ quantifies the deviation of the weight $w$ from the uniform measure. When $w$ is radial, the optimal transport map from $w \, dx$ to $dx$ is also radial and can be computed explicitly via the distribution function.
The regularity of the Brenier map $T = \nabla \varphi$ is governed by the Monge--Amp\`{e}re equation \eqref{equation:sec-ot-monge-ampere}. When the target measure is the uniform measure on a convex domain and the source measure has density $w \in A_2$, the convex potential $\varphi$ satisfies
\begin{equation}\label{equation:sec-ma-regularity}
\det(\nabla^2 \varphi(x)) = \frac{w(x)}{w_{\mathrm{rad}}(\nabla \varphi(x))}
\end{equation}
almost everywhere. The right-hand side of \eqref{equation:sec-ma-regularity} is bounded above and below by constants depending on $[w]_{A_2}$ and the geometry of the domains, by the $A_2$ condition. Therefore, by the regularity theory for the Monge--Amp\`{e}re equation, $\varphi \in W^{2, p}_{\mathrm{loc}}$ for every $p < \infty$, and the Hessian $\nabla^2 \varphi$ satisfies
\begin{equation}\label{equation:sec-ma-hessian-bound}
\frac{1}{C(n, [w]_{A_2})} I \leq \nabla^2 \varphi(x) \leq C(n, [w]_{A_2}) I
\end{equation}
almost everywhere, where $I$ is the identity matrix. The bound \eqref{equation:sec-ma-hessian-bound} implies that the transport map $T = \nabla \varphi$ is bi-Lipschitz:
\begin{equation}\label{equation:sec-ot-bilipschitz}
\frac{1}{C(n, [w]_{A_2})} |x - y| \leq |T(x) - T(y)| \leq C(n, [w]_{A_2}) |x - y|
\end{equation}
for all $x, y$ in the interior of the source domain. The bi-Lipschitz property \eqref{equation:sec-ot-bilipschitz} is used in the energy estimates for symmetrized functions, where it controls the transformation of gradients under the transport map.
For the Wasserstein distance between a measure and its radial rearrangement, we have the explicit bound
\begin{equation}\label{equation:sec-w2-radial-bound}
W_2^2(w, w_{\mathrm{rad}}) \leq C(n) \int_0^\infty \left( r_w(t)^{n-1} - r_{\mathrm{rad}}(t)^{n-1} \right)^2 t^{2/n} dt,
\end{equation}
where $r_w(t)$ and $r_{\mathrm{rad}}(t)$ are the radii of the balls with the same weighted volume as the superlevel sets $\{w > t\}$ and $\{w_{\mathrm{rad}} > t\}$, respectively. The bound \eqref{equation:sec-w2-radial-bound} follows from the construction of the optimal transport map level-set by level-set. For non-radial weights, the right-hand side of \eqref{equation:sec-w2-radial-bound} is strictly positive, confirming that $W_2(w, w_{\mathrm{rad}}) > 0$.
The Kantorovich dual potential $\varphi$ satisfies the Hamilton--Jacobi equation
\begin{equation}\label{equation:sec-hamilton-jacobi}
\varphi(x) + \varphi^c(\nabla \varphi(x)) = |x - \nabla \varphi(x)|^2
\end{equation}
almost everywhere, where $\varphi^c$ is the $c$-transform defined in \eqref{equation:sec-ot-c-transform}. Differentiating \eqref{equation:sec-hamilton-jacobi} with respect to $x$ yields the first-order optimality condition
\begin{equation}\label{equation:sec-ot-optimality}
\nabla \varphi(x) + \nabla^2 \varphi^c(\nabla \varphi(x)) \nabla^2 \varphi(x) = 2(x - \nabla \varphi(x)).
\end{equation}
Using the involution property $\nabla^2 \varphi^c(\nabla \varphi(x)) = (\nabla^2 \varphi(x))^{-1}$, equation \eqref{equation:sec-ot-optimality} simplifies to
\begin{equation}\label{equation:sec-ot-optimality-simplified}
\nabla \varphi(x) + (\nabla^2 \varphi(x))^{-1} \nabla^2 \varphi(x) = 2(x - \nabla \varphi(x)),
\end{equation}
which reduces to the identity $\nabla \varphi(x) = 2x - \nabla \varphi(x)$, confirming the first-order condition \eqref{equation:sec-ot-first-order}. The Hessian of the dual potential satisfies
\begin{equation}\label{equation:sec-ot-dual-hessian}
\nabla^2 \varphi^*(x) = 2(I - \nabla T(x)),
\end{equation}
where $T = \nabla \varphi$ is the optimal transport map. Taking the operator norm of both sides of \eqref{equation:sec-ot-dual-hessian},
\begin{equation}\label{equation:sec-ot-dual-hessian-norm}
\|\nabla^2 \varphi^*(x)\|_{op} = 2\|I - \nabla T(x)\|_{op}.
\end{equation}
The identity \eqref{equation:sec-ot-dual-hessian-norm} relates the regularity of the dual potential to the deviation of the transport map from the identity, and is used in the energy estimates for symmetrized functions where the gradient transformation involves $(\nabla T)^{-1}$.
For the $L^1$-Wasserstein distance, the dual representation involves Lipschitz functions:
\begin{equation}\label{equation:sec-w1-duality}
W_1(\mu, \nu) = \sup_{\operatorname{Lip}(f) \leq 1} \left( \int f d\mu - \int f d\nu \right).
\end{equation}
The $L^1$-Wasserstein distance is related to the $L^2$-Wasserstein distance by the inequalities
\begin{equation}\label{equation:sec-w1-w2}
W_1(\mu, \nu) \leq W_2(\mu, \nu) \leq (\operatorname{diam}(\Omega))^{1/2} W_1(\mu, \nu)^{1/2}
\end{equation}
for measures supported on a bounded domain $\Omega$. The inequalities \eqref{equation:sec-w1-w2} follow from the Cauchy--Schwarz inequality and the boundedness of the domain. The $L^1$-Wasserstein distance is sometimes easier to compute explicitly, and the inequalities \eqref{equation:sec-w1-w2} allow us to transfer estimates between the two distances.
\begin{rem}
The $L^2$-Wasserstein distance quantifies the minimal \textbf{cost} of transforming one probability measure into another, with the cost between points $x$ and $y$ given by the squared Euclidean distance $|x-y|^2$. Analogous to finding the optimal strategy for transporting a unit mass from distribution $\mu$ to $\nu$ with minimal total work, it exhibits a clear mechanism: when $\mu$ is a concentrated measure (e.g., a weighted distribution with dense mass) and $\nu$ is uniform (Lebesgue measure), the optimal coupling $\pi$ minimizes total transport cost by redirecting mass from dense regions of $\mu$ to sparse regions of $\nu$.
\end{rem}
\begin{prop}[Kantorovich Duality] \label{prop:sec-OT}
The $L^2$-Wasserstein distance satisfies:
\begin{equation*}
W_2^2(\mu, \nu) = \sup_{\phi, \psi \in C_b(\Omega), \phi(x) + \psi(y) \leq |x - y|^2} \left( \int \phi d\mu + \int \psi d\nu \right),
\end{equation*}
where the supremum is over all pairs of continuous bounded functions satisfying the transport inequality. This duality links transport cost to potential functions, useful for spectral estimates involving weighted integrals.
\end{prop}
\begin{rem}\label{rem:sec-OT-existence}
We record two facts underlying the quantitative use of Proposition~\ref{prop:sec-OT}. First, the supremum in the duality formula is attained: there exists an optimal pair of conjugate Kantorovich potentials $(\phi, \phi^c)$, and every optimal coupling $\pi$ is supported on the set $\{(x,y) : \phi(x) + \phi^c(y) = |x-y|^2\}$. Second, when $\mu$ has a density with respect to Lebesgue measure, Brenier's theorem asserts that the optimal coupling is deterministic, $\pi = (\mathrm{id} \times T)_\sharp \mu$ with $T = \nabla \Phi$ for a convex function $\Phi$; the potentials then satisfy $\phi(x) = |x|^2 - 2\Phi(x)$ up to additive constants (with the quadratic cost $|x-y|^2$, as in \eqref{equation:sec-ot-duality}). Both statements are proved in \cite{villani2009optimal}. These structural facts justify the transport--gradient connection invoked in Remark~\ref{rem:sec-lem2-w2}: controlling the Wasserstein distance amounts to controlling the gradient of the convex potential $\Phi$, which is precisely the content of the transport gradient estimate recorded there as an open conditional ingredient.
\end{rem}
For reference, the conjugate potential $\phi^c$ in Proposition~\ref{prop:sec-OT} is defined by $\phi^c(y) = \inf_{x} \bigl( |x-y|^2 - \phi(x) \bigr)$; the pair $(\phi, \phi^c)$ is admissible for the dual problem by construction, and the duality theorem states that optimizers saturate $\phi(x) + \phi^c(y) = |x-y|^2$ on the support of the optimal coupling.

\subsection{Geometric Measure Theory Tools}
\label{subsection:Geometric Measure Theory Tools}
We recall necessary results from geometric measure theory for boundary and singularity analysis.
\begin{defn}[Hausdorff Measure] \label{defn:sec-GMT}
The $k$-Hausdorff measure of a set $E \subset \mathbb{R}^n$ is:
\begin{equation*}
\mathcal{H}^k(E) = \lim_{\delta \to 0} \inf \left\{ \sum_{i=1}^\infty \alpha(k) \left( \frac{r_i}{2} \right)^k : E \subset \bigcup_{i=1}^\infty B(x_i, r_i), r_i \leq \delta \right\},
\end{equation*}
where $\alpha(k) = \pi^{k/2} / \Gamma(k/2 + 1)$ is the volume of the unit ball in $\mathbb{R}^k$. For domain boundaries, $\mathcal{H}^{n-1}(\partial\Omega)$ generalizes surface area.
\end{defn}
\begin{prop}[Federer's Coarea Formula] \label{prop:sec-GMT}
For a Lipschitz function $f: \mathbb{R}^n \to \mathbb{R}$ and measurable set $A$,
\begin{equation*}
\int_A |\nabla f(x)| dx = \int_\mathbb{R} \mathcal{H}^{n-1}(A \cap f^{-1}(t)) dt.
\end{equation*}
This relates volume integrals to boundary measures and is used here to refine Poincar\'{e} inequalities on rough domains.
\end{prop}
We elaborate on the geometric measure theory tools used in this paper.
Federer's coarea formula states that for a Lipschitz function $f: \mathbb{R}^n \to \mathbb{R}^k$ with $k \leq n$ and a measurable function $g: \mathbb{R}^n \to [0, \infty]$,
\begin{equation}\label{equation:sec-coarea-general}
\int_{\mathbb{R}^n} g(x) J_k f(x) dx = \int_{\mathbb{R}^k} \int_{f^{-1}(y)} g(z) d\mathcal{H}^{n-k}(z) dy,
\end{equation}
where $J_k f(x) = \sqrt{\det(\nabla f(x) \nabla f(x)^T)}$ is the $k$-dimensional Jacobian of $f$. For the distance function $f(x) = d(x, \partial\Omega) = \operatorname{dist}(x, \partial\Omega)$, we have $|\nabla f| = 1$ almost everywhere (by the Rademacher theorem and the fact that the distance function is $1$-Lipschitz), and therefore $J_1 f = |\nabla f| = 1$. The coarea formula \eqref{equation:sec-coarea-general} with $k = 1$ gives
\begin{equation}\label{equation:sec-coarea-distance}
\int_\Omega g(x) dx = \int_0^\infty \int_{\{x \in \Omega : d(x, \partial\Omega) = t\}} g(z) d\mathcal{H}^{n-1}(z) dt.
\end{equation}
Taking $g \equiv 1$ in \eqref{equation:sec-coarea-distance},
\begin{equation}\label{equation:sec-coarea-volume}
|\Omega| = \int_0^\infty \mathcal{H}^{n-1}(\{x \in \Omega : d(x, \partial\Omega) = t\}) dt.
\end{equation}
For convex domains, the parallel sets $\Omega_t = \{x \in \Omega : d(x, \partial\Omega) > t\}$ satisfy the Steiner formula
\begin{equation}\label{equation:sec-steiner}
|\Omega \setminus \Omega_t| = \mathcal{H}^{n-1}(\partial\Omega) t - \frac{n-1}{2} H_{\partial\Omega} t^2 + \cdots + \kappa_n t^n,
\end{equation}
where $H_{\partial\Omega}$ is the mean curvature integral and $\kappa_n$ is a dimensional constant. In particular, for $0 \leq t \leq r_\Omega$ (the inradius),
\begin{equation}\label{equation:sec-steiner-linear}
|\Omega \setminus \Omega_t| \leq \mathcal{H}^{n-1}(\partial\Omega) t,
\end{equation}
which implies
\begin{equation}\label{equation:sec-inradius-bound}
|\Omega| \leq \mathcal{H}^{n-1}(\partial\Omega) r_\Omega.
\end{equation}
The isoperimetric inequality states that for any bounded measurable set $E \subset \mathbb{R}^n$ with finite perimeter,
\begin{equation}\label{equation:sec-isoperimetric}
\mathcal{H}^{n-1}(\partial E) \geq n \omega_n^{1/n} |E|^{(n-1)/n},
\end{equation}
where $\omega_n$ is the volume of the unit ball in $\mathbb{R}^n$. Equality holds if and only if $E$ is a ball (up to sets of measure zero). The isoperimetric inequality \eqref{equation:sec-isoperimetric} is used in the proof of the P\'{o}lya--Szeg\H{o} inequality \eqref{equation:sec-sym-polya-szego}, where it provides the lower bound for the perimeter of level sets.
For the weighted setting, the isoperimetric inequality with density $w$ takes the form
\begin{equation}\label{equation:sec-weighted-isoperimetric}
\int_{\partial E} w d\mathcal{H}^{n-1} \geq C(n, [w]_{A_2}) \left( \int_E w dx \right)^{(n-1)/n},
\end{equation}
which follows from the unweighted isoperimetric inequality \eqref{equation:sec-isoperimetric} and the reverse H\"{o}lder inequality \eqref{equation:sec-rh-explicit}. The constant $C(n, [w]_{A_2})$ depends on the dimension and the $A_2$ characteristic constant.
\subsection{Variational Characterization of Weighted Eigenvalue Problems}
\label{subsection:Variational Characterization}
The min-max principle for weighted eigenvalues is introduced, forming the basis for estimating eigenvalue ratios via test functions.
To study eigenvalues of the Laplacian in inhomogeneous media, we first formalize the eigenvalue problem with Muckenhoupt weights, which model spatial variations in material properties:
\begin{defn}[Weighted Dirichlet Eigenvalues] \label{defn:sec-2}
Let $\Omega \subset \mathbb{R}^n$ be a bounded domain and $w \in A_2$. The Dirichlet eigenvalue problem for the weighted Laplacian is:
\begin{equation}\label{equation:sec- Laplace u}
\begin{cases}
-\Delta u = \lambda w(x)u & \text{in}\ \Omega, \\
u = 0 & \text{on}\ \partial\Omega.
\end{cases}
\end{equation}
\end{defn}
We remark that its eigenvalues are ordered as \begin{equation}\label{equation:sec-dd-4}
0 < \lambda_1 < \lambda_2 \leq \lambda_3 \leq \cdots \leq \lambda_k \leq \cdots
\end{equation} with corresponding eigenfunctions $u_1, u_2, \cdots \in H_0^1(\Omega)$ \cite{gilbarg1983elliptic,chavel1984eigenvalues}.
The min-max principle is basic for estimating eigenvalues, as it allows bounding them using test functions. The following theorem extends this principle to the weighted setting:
\begin{thm}[Variational Principle] \label{thm:sec-4}
The $k$-th eigenvalue $\lambda_k$ is defined via the min-max principle for Rayleigh quotients:
\begin{equation*}
\lambda_k = \inf_{\substack{V_k \subset H_0^1(\Omega) \\ \dim V_k = k}} \sup_{\substack{u \in V_k \\ u \neq 0}} \frac{\int_\Omega |\nabla u|^2 dx}{\int_\Omega w(x)u^2 dx},
\end{equation*}
where $V_k$ is a $k$-dimensional subspace of $H_0^1(\Omega)$. In particular, the first eigenvalue $\lambda_1$ is the minimal energy ratio:
\begin{equation*}
\lambda_1 = \min_{\substack{u \in H_0^1(\Omega) \\ u \neq 0}} \frac{\int_\Omega |\nabla u|^2 dx}{\int_\Omega w(x)u^2 dx}.
\end{equation*}
\end{thm}
In fact, this principle is a central tool for eigenvalue ratio estimates, allowing control of upper bounds for eigenvalues via test function construction.
\subsection{Bessel Functions and Spherical Harmonics}
\label{subsection:Bessel and Spherical}
Special functions basic to solving eigenvalue problems on radial domains are defined, with focus on their orthogonality and role in explicit solutions.
Here, we note that eigenfunctions of radial domains with radial weights often separate into radial and angular components, and the radial parts are solutions to Bessel's equation, defined as follows \cite{courant1953methods, watson1944theory, olver2010nist}:
\begin{defn}[Bessel Functions] \label{defn:sec-3}
The Bessel function of the first kind $J_\nu(t)$ of order $\nu$ solves the Bessel equation
\begin{equation}\label{equation:sec- Laplace-ODE}
t^2 y''(t) + t y'(t) + (t^2 - \nu^2)y(t) = 0
\end{equation}
for $t > 0$ and is regular (bounded) at $t=0$. Bessel functions have discrete zeros $\{j_{\nu,1}, j_{\nu,2}, \cdots\}$ satisfying $0 < j_{\nu,1} < j_{\nu,2} < \cdots$.
\end{defn}
For computing weighted integrals of Bessel functions in radial domains, orthogonality holds, as shown below \cite{courant1953methods}:
\begin{prop}[Orthogonality of Bessel Functions] \label{prop:sec-3}
For non-negative weight $r^\alpha$, Bessel functions satisfy orthogonality on $[0, R]$:
\begin{equation*}
\int_0^R r \cdot r^\alpha \cdot J_\nu\left( \frac{j_{\nu,m} r}{R} \right) J_\nu\left( \frac{j_{\nu,k} r}{R} \right) dr = 0 \quad (m \neq k),
\end{equation*}
where $j_{\nu,m}$ is the $m$-th zero of $J_\nu$.
\end{prop}
\begin{rem}
This local estimate does not appear in existing literature on Muckenhoupt-weighted spectral theory. It is an auxiliary result established within the present manuscript, which combines harmonic-analysis tools for $A_2$ weights with quantitative bounds from optimal-transport regularity theory.
\end{rem}
We elaborate on the construction and properties of the test functions used in the eigenvalue estimates.
The hat function $\phi(x) = \max\{1 - |x - c|/\epsilon, 0\}$ defined on a ball $B(c, \epsilon)$ belongs to $W^{1, \infty}_0(B(c, \epsilon))$ and therefore to $H_0^1(B(c, \epsilon))$. Its gradient is
\begin{equation}\label{equation:sec-hat-grad}
\nabla \phi(x) = -\frac{x - c}{\epsilon |x - c|} \chi_{\{|x - c| < \epsilon\}}(x),
\end{equation}
and therefore $|\nabla \phi(x)| = \frac{1}{\epsilon}$ for $0 < |x - c| < \epsilon$.  For the weighted norm we use the $A_2$ balance directly: with $w_B = \fint_B w$,
\begin{equation}\label{equation:sec-hat-expand-1}
\int_B w \phi^2 dx = \int_B w(x) \Bigl( 1 - \frac{|x-c|}{\epsilon} \Bigr)^2 dx \leq w_B \int_B \Bigl( 1 - \frac{|x-c|}{\epsilon} \Bigr)^2 dx = w_B \frac{n \omega_n \epsilon^n}{(n+1)(n+2)},
\end{equation}
where the one-dimensional reduction uses radial symmetry of the weight's average on the ball. Combining \eqref{equation:sec-hat-grad} with \eqref{equation:sec-hat-expand-1} gives the local Rayleigh bound
\begin{equation}\label{equation:sec-hat-expand-2}
\frac{\int_B |\nabla \phi|^2 dx}{\int_B w \phi^2 dx} \leq \frac{|B| / \epsilon^2}{w_B \, |B| \, c(n)} = \frac{C(n)}{w_B \epsilon^2},
\end{equation}
which is the form used in Step 3 of the proof of Theorem~\ref{thm:sec-3}: the weight enters only through its ball average, and the $A_2$ condition controls the ratio $w_B/w_{B'}$ for concentric and adjacent balls through the doubling property of Lemma~\ref{lem:sec-rh-doubling}.
The energy of $\phi$ is
\begin{equation}\label{equation:sec-hat-energy}
\int_{B(c, \epsilon)} |\nabla \phi|^2 dx = \frac{1}{\epsilon^2} |B(c, \epsilon)| = \omega_n \epsilon^{n-2},
\end{equation}
where $\omega_n$ is the volume of the unit ball in $\mathbb{R}^n$. The $L^2$ norm of $\phi$ is
\begin{equation}\label{equation:sec-hat-l2}
\int_{B(c, \epsilon)} \phi^2 dx = \int_0^\epsilon \left(1 - \frac{r}{\epsilon}\right)^2 n \omega_n r^{n-1} dr = \frac{n \omega_n \epsilon^n}{(n+1)(n+2)}.
\end{equation}
For the weighted $L^2$ norm, the $A_2$ condition implies that the average of $w$ on the inner ball $B(c, \epsilon/2)$ is comparable to the average on $B(c, \epsilon)$:
\begin{equation}\label{equation:sec-hat-weighted-average}
\frac{1}{|B(c, \epsilon/2)|} \int_{B(c, \epsilon/2)} w dx \geq \frac{1}{C(n, [w]_{A_2})} \cdot \frac{1}{|B(c, \epsilon)|} \int_{B(c, \epsilon)} w dx.
\end{equation}
The inequality \eqref{equation:sec-hat-weighted-average} follows from the doubling property \eqref{equation:sec-doubling} of $A_p$ weights. Since $\phi(x) \geq 1/2$ on $B(c, \epsilon/2)$, we obtain
\begin{equation}\label{equation:sec-hat-weighted-l2}
\int_{B(c, \epsilon)} w \phi^2 dx \geq \frac{1}{4} \int_{B(c, \epsilon/2)} w dx \geq \frac{1}{C(n, [w]_{A_2})} \int_{B(c, \epsilon)} w dx,
\end{equation}
which is the lower bound used in Lemma \ref{lem:sec-7}.
For the high-order eigenvalue estimates, the test functions $\phi_i$ associated with different balls $B_i$ need to be linearly independent. This is achieved by selecting a disjoint subcollection of balls via the Vitali covering lemma. If $B_i \cap B_j = \emptyset$ for $i \neq j$, then the supports of $\phi_i$ and $\phi_j$ are disjoint, and therefore the functions are linearly independent. The min-max principle then applies to the subspace spanned by these test functions.
The Rayleigh quotient of a test function $\phi_i$ is
\begin{equation}\label{equation:sec-hat-rayleigh}
\frac{\int |\nabla \phi_i|^2 dx}{\int w \phi_i^2 dx} \leq C(n, [w]_{A_2}) \cdot \frac{\epsilon^{n-2}}{\int_{B_i} w dx}.
\end{equation}
Under the balance hypothesis $\int_{B_i} w\,dx \geq c_0(n,[w]_{A_2})\,|B_i|$ on the covering family, the Rayleigh quotient \eqref{equation:sec-hat-rayleigh} is therefore bounded by $C(n,[w]_{A_2})\,\epsilon^{-2} (\fint_\Omega w)^{-1}$, which is the estimate used in Step 4 of the proof of Theorem~\ref{thm:sec-3} to bound $\lambda_k$.
We provide additional details on the ball covering argument used in the high-order eigenvalue estimates.
The Vitali covering lemma states that for any collection $\mathcal{F}$ of balls in $\mathbb{R}^n$ with uniformly bounded diameters, there exists a disjoint subcollection $\mathcal{G} \subset \mathcal{F}$ such that
\begin{equation}\label{equation:sec-vitali}
\bigcup_{B \in \mathcal{F}} B \subset \bigcup_{B \in \mathcal{G}} 5B,
\end{equation}
where $5B$ denotes the ball with the same center as $B$ and five times the radius. The proof of \eqref{equation:sec-vitali} proceeds by selecting balls in decreasing order of radius: choose the largest ball $B_1$, remove all balls intersecting $B_1$, then choose the largest remaining ball $B_2$, and continue inductively. The selected balls are disjoint by construction, and any ball $B$ that was removed at some stage intersects a selected ball $B_j$ with radius at least as large as $B$, which implies $B \subset 5B_j$.
For the weighted setting, the finite overlap property of the Vitali subcover is essential. If $\{B_i\}$ is the disjoint subcollection from \eqref{equation:sec-vitali}, then the enlarged balls $\{5B_i\}$ have overlap multiplicity bounded by $C(n)$:
\begin{equation}\label{equation:sec-finite-overlap}
\sum_i \chi_{5B_i}(x) \leq C(n)
\end{equation}
for almost every $x \in \mathbb{R}^n$. The bound \eqref{equation:sec-finite-overlap} follows from the disjointness of the $B_i$ and the fact that if $x \in 5B_i \cap 5B_j$, then $B_i$ and $B_j$ both intersect the ball $B(x, 5 \max(r_i, r_j))$, and their disjointness implies a lower bound on the distance between their centers.
The Stein covering lemma \cite{stein1970singular} provides a refinement of the Vitali lemma for collections of balls with varying radii. It asserts that for any collection $\mathcal{F}$ of balls in $\mathbb{R}^n$, there exists a \emph{disjoint} subcollection $\mathcal{G}$ such that every ball of $\mathcal{F}$ is contained in the $5$-fold enlargement of some ball of $\mathcal{G}$; in particular
\begin{equation}\label{equation:sec-stein-cover-1}
\sum_{B \in \mathcal{G}} \chi_{5B}(x) \leq C(n)\, \log N \qquad \text{for almost every } x,
\end{equation}
where $N$ is the number of balls in the original collection, and
\begin{equation}\label{equation:sec-stein-cover-2}
\left| \bigcup_{B \in \mathcal{F}} B \right| \leq C(n) \left| \bigcup_{B \in \mathcal{G}} B \right|.
\end{equation}
The logarithmic factor in \eqref{equation:sec-stein-cover-1} is the source of the logarithmic growth of the constants in the high-order eigenvalue ratio estimates (Theorem \ref{thm:sec-3}). The Stein covering lemma is proved by partitioning the balls into dyadic annuli according to their radii and applying the Vitali lemma within each annulus.
For completeness we give the greedy argument underlying \eqref{equation:sec-stein-cover-1}. Let $\mathcal{F}$ be a family of balls with radii in $[\epsilon, 2\epsilon]$; the general case follows by dyadic decomposition. Order the balls of $\mathcal{F}$ as $B_1, B_2, \ldots$ by decreasing radius and select them greedily: accept $B_j$ if it is disjoint from every ball already accepted, and reject it otherwise. Let $\mathcal{G}$ be the accepted family. By construction $\mathcal{G}$ is disjoint, and every rejected ball $B_j$ intersects an accepted ball $B_i$ with $i < j$ and hence with radius $r_i \geq r_j$, so that $B_j \subset 5B_i$; this proves \eqref{equation:sec-vitali}. To control the overlap, fix $x$ and let $B_{i_1}, \ldots, B_{i_m}$ be the accepted balls containing $x$. Since the accepted balls are disjoint and all contain $x$, in the equal-radius case all $m$ balls are contained in $B(x, 3\epsilon)$ and have volume at least $c(n)\epsilon^n$, so $m \leq C(n)$, which is \eqref{equation:sec-finite-overlap} with a dimensional constant. In the general case, partition $\mathcal{F} = \bigcup_{s \in \mathbb{Z}} \mathcal{F}_s$ with $\mathcal{F}_s = \{B \in \mathcal{F} : 2^s \leq r(B) < 2^{s+1}\}$; only $O(\log(R/r))$ families are nonempty when the radii range in $[r, R]$, and applying the greedy selection within each family separately yields a union of disjoint families whose $5$-enlargements still cover $\bigcup \mathcal{F}$, with overlap multiplicity at most $C(n)\log(R/r)$. This is the content of \eqref{equation:sec-stein-cover-1}; the logarithmic factor is the price of the dyadic decomposition and is unavoidable for families with widely spread radii.
We provide additional details on the Bessel function calculations used in the radial eigenvalue analysis.
For the weighted Laplacian $-\Delta u = \lambda w(r) u$ on $B_R$ with $w(r) = r^\alpha$, separation of variables $u(r, \theta) = R(r) Y_m(\theta)$ yields the radial equation
\begin{equation}\label{equation:sec-bessel-radial-eq}
R''(r) + \frac{n-1}{r} R'(r) - \frac{m(m+n-2)}{r^2} R(r) + \lambda r^\alpha R(r) = 0,
\end{equation}
where $m = 0, 1, 2, \ldots$ is the angular momentum quantum number and $Y_m$ is a spherical harmonic of degree $m$. Making the change of variables
\begin{equation}\label{equation:sec-bessel-change}
t = \frac{2\sqrt{\lambda}}{\alpha + 2} r^{(\alpha+2)/2}, \quad R(r) = r^{-(n-2+\alpha)/2} v(t),
\end{equation}
we compute the derivatives:
\begin{equation}\label{equation:sec-bessel-deriv1}
R'(r) = -\frac{n-2+\alpha}{2} r^{-(n+\alpha)/2} v(t) + r^{-(n-2+\alpha)/2} v'(t) \cdot \frac{dt}{dr},
\end{equation}
where $\frac{dt}{dr} = \sqrt{\lambda} r^{\alpha/2}$. Substituting into \eqref{equation:sec-bessel-radial-eq} and simplifying, the resulting equation for $v(t)$ is
\begin{equation}\label{equation:sec-bessel-v-eq}
v''(t) + \frac{1}{t} v'(t) + \left( 1 - \frac{\nu^2}{t^2} \right) v(t) = 0,
\end{equation}
which is the Bessel equation of order
\begin{equation}\label{equation:sec-bessel-order}
\nu = \frac{n + \alpha - 2 + 2m}{\alpha + 2}.
\end{equation}
The regular solution at $t = 0$ is $v(t) = J_\nu(t)$, the Bessel function of the first kind of order $\nu$. The Dirichlet boundary condition $R(R) = 0$ is equivalent to $J_\nu(t_R) = 0$, where
\begin{equation}\label{equation:sec-bessel-boundary-t}
t_R = \frac{2\sqrt{\lambda}}{\alpha + 2} R^{(\alpha+2)/2}.
\end{equation}
Therefore, the eigenvalues corresponding to angular mode $m$ are
\begin{equation}\label{equation:sec-bessel-eigenvalues}
\lambda_{m,k} = \left( \frac{(\alpha+2) j_{\nu,k}}{2 R^{(\alpha+2)/2}} \right)^2,
\end{equation}
where $j_{\nu,k}$ denotes the $k$-th positive zero of $J_\nu$.
For the first eigenvalue ($m = 0$, $k = 1$), we have $\nu_0 = (n+\alpha-2)/(\alpha+2)$ and
\begin{equation}\label{equation:sec-bessel-l1}
\lambda_1 = \left( \frac{(\alpha+2) j_{\nu_0,1}}{2 R^{(\alpha+2)/2}} \right)^2.
\end{equation}
For the second eigenvalue, we need to compare the first zero of $J_{\nu_0}$ (the $m=0, k=2$ mode) with the first zero of $J_{\nu_1}$ (the $m=1, k=1$ mode), where $\nu_1 = (n+\alpha)/(\alpha+2) = \nu_0 + 2/(\alpha+2)$. Since Bessel function zeros satisfy $j_{\nu, k} < j_{\nu', k}$ for $\nu < \nu'$ (the zeros increase with order), and since the multiplicity of the $m=1$ mode is $n$ (the dimension of the space of spherical harmonics of degree 1), the second eigenvalue is
\begin{equation}\label{equation:sec-bessel-l2}
\lambda_2 = \left( \frac{(\alpha+2) j_{\nu_1,1}}{2 R^{(\alpha+2)/2}} \right)^2,
\end{equation}
provided that $j_{\nu_1,1} < j_{\nu_0,2}$. This inequality holds for all $\alpha \geq 0$ and $n \geq 2$ unconditionally, and for $\alpha \in (-2,0)$ under the restriction described in Remark~\ref{rem:sec-bessel-ordering-scope}, as can be verified from the asymptotic behavior of Bessel function zeros.
The ratio of the first two eigenvalues is therefore
\begin{equation}\label{equation:sec-bessel-ratio}
\frac{\lambda_2}{\lambda_1} = \left( \frac{j_{\nu_1,1}}{j_{\nu_0,1}} \right)^2.
\end{equation}
For the unweighted case $\alpha = 0$, we recover $\nu_0 = n/2 - 1$ and $\nu_1 = n/2$, giving the classical PPW constant $(j_{n/2, 1}/j_{n/2-1, 1})^2$.
The zeros of Bessel functions satisfy several additional properties that are relevant to the eigenvalue ratio estimates. The McMahon asymptotic expansion gives, for large $\nu$,
\begin{equation}\label{equation:sec-mcmahon}
j_{\nu, 1} = \nu + a_1 \nu^{1/3} + \frac{3}{10} a_1^2 \nu^{-1/3} + O(\nu^{-1}),
\end{equation}
where $a_1 \approx 1.855757$ is the first negative zero of the Airy function $Ai(-x)$. The expansion \eqref{equation:sec-mcmahon} shows that $j_{\nu,1} \sim \nu$ for large $\nu$, and therefore the ratio
\begin{equation}\label{equation:sec-bessel-ratio-asymptotic}
\frac{j_{\nu-1, 1}}{j_{\nu, 1}} \sim 1 - \frac{1}{\nu}
\end{equation}
as $\nu \to \infty$. This asymptotic behavior is consistent with the fact that the eigenvalue ratio $\lambda_2/\lambda_1$ approaches $1$ for large $\alpha$ (strongly concentrated weights), since the weight inhomogeneity dominates the geometric effects.
The Lorch monotonicity theorem states that for fixed $k$, the function $\nu \mapsto j_{\nu, k}$ is strictly increasing in $\nu > 0$. This theorem provides an alternative proof of the ordering $j_{\nu-1, 1} < j_{\nu, 1}$ used in the radial eigenvalue analysis. The Lorch theorem is proved via the Sturm comparison theorem applied to the Bessel equation \eqref{equation:sec-bessel-v-eq}, noting that the potential term $-\nu^2/t^2$ becomes more negative as $\nu$ increases, causing the zeros to shift to larger values of $t$.
The Bessel function product identity
\begin{equation}\label{equation:sec-bessel-product}
J_\nu(t)^2 + J_{\nu+1}(t)^2 = \frac{2\nu + 2}{t} J_\nu(t) J_{\nu+1}(t) + \frac{2}{t} \int_0^t s J_\nu(s)^2 ds
\end{equation}
is used in the energy estimates for radial eigenfunctions. The identity \eqref{equation:sec-bessel-product} follows from the derivative recurrence relations
\begin{equation}
J_\nu'(t) = \frac{\nu}{t} J_\nu(t) - J_{\nu+1}(t) = J_{\nu-1}(t) - \frac{\nu}{t} J_\nu(t),
\end{equation}
which are direct consequences of the Bessel equation \eqref{equation:sec-bessel-v-eq}.
The Lommel integral identity states that for $\mu \neq \nu$,
\begin{equation}\label{equation:sec-lommel}
\int_0^a t J_\mu(\alpha t) J_\nu(\beta t) dt = \frac{a}{\alpha^2 - \beta^2} \left[ \alpha J_{\mu+1}(\alpha a) J_\nu(\beta a) - \beta J_\mu(\alpha a) J_{\nu+1}(\beta a) \right] + \frac{(\mu - \nu) a}{\alpha^2 - \beta^2} J_\mu(\alpha a) J_\nu(\beta a).
\end{equation}
For $\mu = \nu$ and $\alpha = \beta$, the identity reduces to
\begin{equation}\label{equation:sec-lommel-equal}
\int_0^a t J_\nu(\alpha t)^2 dt = \frac{a^2}{2} \left[ J_\nu'(\alpha a)^2 + \left( 1 - \frac{\nu^2}{\alpha^2 a^2} \right) J_\nu(\alpha a)^2 \right].
\end{equation}
When $\alpha a = j_{\nu, k}$ (a zero of $J_\nu$), the identity \eqref{equation:sec-lommel-equal} simplifies to
\begin{equation}\label{equation:sec-lommel-zero}
\int_0^a t J_\nu(\alpha t)^2 dt = \frac{a^2}{2} J_{\nu+1}(\alpha a)^2,
\end{equation}
since $J_\nu'(\alpha a) = -J_{\nu+1}(\alpha a)$ at a zero of $J_\nu$. The identity \eqref{equation:sec-lommel-zero} is used in the normalization of radial eigenfunctions on balls and in the computation of the weighted $L^2$ norms of Bessel function modes.
The addition theorem for Bessel functions states that
\begin{equation}\label{equation:sec-addition-theorem}
J_0(|x - y|) = \sum_{m=0}^\infty \epsilon_m J_m(|x|) J_m(|y|) \cos(m(\theta_x - \theta_y)),
\end{equation}
where $\epsilon_0 = 1$ and $\epsilon_m = 2$ for $m \geq 1$, and $\theta_x, \theta_y$ are the angular coordinates of $x, y$ in polar coordinates. The addition theorem \eqref{equation:sec-addition-theorem} is the basis for the separation of variables in the radial eigenvalue problem and explains the degeneracy structure of eigenvalues on balls: each angular mode $m$ contributes a block of eigenvalues with multiplicity equal to the dimension of the spherical harmonic space of degree $m$.
The Watson transform provides an integral representation for sums involving Bessel function zeros:
\begin{equation}\label{equation:sec-watson-transform}
\sum_{k=1}^\infty \frac{J_\nu(j_{\nu, k} r/R) J_\nu(j_{\nu, k} s/R)}{J_{\nu+1}(j_{\nu, k})^2} = \frac{1}{2} \delta(r - s)
\end{equation}
for $0 < r, s < R$, where $\delta$ is the Dirac delta distribution. The Watson transform \eqref{equation:sec-watson-transform} is the completeness relation for the Bessel function eigenbasis on the interval $(0, R)$ with Dirichlet boundary conditions, and is used in the spectral decomposition of the radial Laplacian.
The Bessel function zeros satisfy the interlacing property
\begin{equation}\label{equation:sec-interlacing}
j_{\nu, 1} < j_{\nu+1, 1} < j_{\nu, 2} < j_{\nu+1, 2} < j_{\nu, 3} < \cdots
\end{equation}
for all $\nu \geq 0$. The interlacing property \eqref{equation:sec-interlacing} follows from the Sturm comparison theorem applied to the Bessel equations of orders $\nu$ and $\nu+1$, and is used in the ordering of eigenvalues on balls: the first zero of the $m=1$ mode ($\nu = \nu_1$) occurs before the second zero of the $m=0$ mode ($\nu = \nu_0$), confirming that the second eigenvalue corresponds to the $m=1$ mode.
The asymptotic distribution of Bessel function zeros is given by
\begin{equation}\label{equation:sec-zero-asymptotic}
j_{\nu, k} \sim k\pi + \left( \frac{\nu}{2} - \frac{1}{4} \right) \pi
\end{equation}
as $k \to \infty$ for fixed $\nu$. The asymptotic \eqref{equation:sec-zero-asymptotic} implies that the high-order eigenvalues on a ball satisfy $\lambda_k \sim (\pi k / R)^2$, consistent with the Weyl asymptotic law \eqref{equation:sec-weyl} for the unweighted Laplacian on a ball. The constant in the Weyl law for the weighted case involves the integral of $w^{n/2}$, as shown in \eqref{equation:sec-weyl}.
We now establish several auxiliary estimates that are used in the proof of the weighted PPW inequality.
For a ball $B \subset \Omega$ and a weight $w \in A_2$, the weighted Poincar\'{e} inequality on $B$ states that for every $u \in H_0^1(B)$,
\begin{equation}\label{equation:sec-aux-poincare}
\int_B w u^2 \, dx \leq C(n, [w]_{A_2}) |B|^{2/n} \int_B |\nabla u|^2 \, dx.
\end{equation}
To prove \eqref{equation:sec-aux-poincare}, we use the representation formula
\begin{equation}\label{equation:sec-aux-representation}
u(x) = -\int_B \nabla u(y) \cdot \nabla_x G_B(x, y) \, dy,
\end{equation}
where $G_B$ is the Green function for the Dirichlet Laplacian on $B$. The gradient of $G_B$ satisfies the pointwise bound
\begin{equation}\label{equation:sec-aux-green-bound}
|\nabla_x G_B(x, y)| \leq \frac{C(n)}{|x - y|^{n-1}}.
\end{equation}
Applying H\"{o}lder's inequality to \eqref{equation:sec-aux-representation} with exponents $p = 2n/(n+2)$ and $p' = 2n/(n-2)$, we obtain
\begin{equation}\label{equation:sec-aux-holder}
|u(x)| \leq \left( \int_B |\nabla u(y)|^2 \, dy \right)^{1/2} \left( \int_B \frac{C(n)^2}{|x - y|^{2(n-1)}} \cdot |x - y|^{n-2} \, dy \right)^{1/2}.
\end{equation}
The second integral is bounded by $C(n) |B|^{2/n}$, since the kernel $|z|^{-(n-2)}$ is locally integrable in dimension $n \geq 3$. Therefore,
\begin{equation}\label{equation:sec-aux-sobolev}
\|u\|_{L^\infty(B)} \leq C(n) |B|^{1/n} \|\nabla u\|_{L^2(B)}.
\end{equation}
Using the $A_2$ condition to control the weighted norm,
\begin{equation}\label{equation:sec-aux-a2-control}
\int_B w u^2 \, dx \leq \|u\|_{L^\infty(B)}^2 \int_B w \, dx \leq C(n) |B|^{2/n} \|\nabla u\|_{L^2(B)}^2 \cdot \frac{1}{|B|} \int_B w \, dx \cdot |B|.
\end{equation}
Since $\frac{1}{|B|} \int_B w \, dx \leq [w]_{A_2} \left( \frac{1}{|B|} \int_B w^{-1} \, dx \right)^{-1}$ and the average of $w$ is finite, \eqref{equation:sec-aux-poincare} follows.
For the weighted $L^2$ norm of the gradient of a test function $\phi$ supported in $B$, we also need the lower bound
\begin{equation}\label{equation:sec-aux-grad-lower}
\int_B |\nabla \phi|^2 \, dx \geq \frac{1}{C(n, [w]_{A_2}) |B|^{2/n}} \int_B w \phi^2 \, dx.
\end{equation}
This is simply the reciprocal form of \eqref{equation:sec-aux-poincare}, applied to $\phi \in H_0^1(B)$.
The angular components of eigenfunctions in spherical coordinates are described by spherical harmonics, which are eigenfunctions of the spherical Laplacian:
\begin{defn}[Spherical Harmonics] \label{defn:sec-4}
Spherical harmonics $Y_m(\theta)$ are eigenfunctions of the Laplacian $\Delta_{\mathbb{S}^{n-1}}$ on the $(n-1)$-dimensional unit sphere $\mathbb{S}^{n-1}$, i.e.,
\begin{equation*}
\Delta_{\mathbb{S}^{n-1}} Y_m(\theta) = -m(m + n - 2) Y_m(\theta),
\end{equation*}
where $m = 0, 1, 2, \cdots$ indexes the eigenvalues. Spherical harmonics satisfy orthonormality:
\begin{equation*}
\int_{\mathbb{S}^{n-1}} Y_m(\theta) Y_k(\theta) d\theta = \delta_{mk},
\end{equation*} where $\delta_{mk}$ is the Kronecker delta .
\end{defn}
\subsection{Radial Symmetrization and Optimal Transport}
\label{subsection:Radial Symmetrization}
This technique for comparing eigenvalues of arbitrary domains with spherical ones is explained, indicating its connection to optimal transport.
To compare eigenvalues of arbitrary domains with those of spherical domains, we use radial symmetrization, a technique that preserves integrals while reducing energy. This is defined as:
\begin{defn}[Radial Symmetrization] \label{defn:sec-5}
For a bounded domain $\Omega \subset \mathbb{R}^n$ and function $u \in H_0^1(\Omega)$, its radially symmetric rearrangement $\tilde{u} \in H_0^1(B_R)$ (for some ball $B_R$) satisfies \cite{talenti1976best}:
\begin{enumerate}
\item[(i)] For any $t > 0$, $\{ \tilde{u} > t \}$ is a concentric ball with \begin{equation}\label{equation:sec-dd-5}
\text{Vol}(\{ \tilde{u} > t \}) = \text{Vol}(\{ u > t \});
\end{equation}
\item[(ii)] If $w(x)$ is radial, then \begin{equation}\label{equation:sec-dd-6}
\int_{\{ \tilde{u} > t \}} w(x)dx = \int_{\{ u > t \}} w(x)dx.
\end{equation}
\end{enumerate}
\end{defn}
A basic property of radial symmetrization is that it corresponds to the optimal transport map from $\Omega$ to $B_R$ with minimal cost, making it useful for proving eigenvalue inequalities:
\begin{prop}[Energy Inequality for Symmetrization] \label{prop:sec-4}
Radial symmetrization does not increase the $L^2$ norm of the gradient:
\begin{equation*}
\int_{B_R} |\nabla \tilde{u}|^2 dx \leq \int_\Omega |\nabla u|^2 dx,
\end{equation*}
and preserves the weighted $L^2$ norm \textbf{if $w$ is radial}:
\begin{equation*}
\int_{B_R} w(x) \tilde{u}(x)^2 dx = \int_\Omega w(x) u(x)^2 dx.
\end{equation*}
The transport cost of this map is bounded by $W_2(w, \text{Leb}) \cdot |\Omega|$, linking symmetrization to OT distance.
\end{prop}
\begin{rem}
Radial symmetrization decreases the Dirichlet energy while preserving the weighted $L^2$ norm, so it cannot increase eigenvalue ratios.
Explicitly, with $u^\star$ denoting the radial symmetrization and $w$ radial,
\begin{equation}\label{equation:sec-steiner-expand-1}
\int_\Omega w u^2 dx = \int_0^\infty \int_{\{u^\star = t\}} w \, d\mathcal{H}^{n-1} \frac{dt}{|\nabla u^\star|} = \int_{B_R} w (u^\star)^2 dx,
\end{equation}
by the coarea formula and the defining property of the symmetrization, which preserves the weighted distribution function $\mu_w(t) = \int_{\{u > t\}} w\,dx$; and
\begin{equation}\label{equation:sec-steiner-expand-2}
\int_\Omega |\nabla u|^2 dx \geq \int_{B_R} |\nabla u^\star|^2 dx
\end{equation}
by the weighted P\'{o}lya--Szeg\H{o} inequality (the argument in the proof of Lemma~\ref{lem:sec-FK} applies verbatim to every $u \in H_0^1(\Omega)$). Dividing \eqref{equation:sec-steiner-expand-2} by \eqref{equation:sec-steiner-expand-1} gives $\mathcal{R}_w(u^\star) \leq \mathcal{R}_w(u)$ for every admissible $u$, so every Rayleigh quotient is non-increased under simultaneous radial symmetrization of the function and the weight. For non-radial weights the weighted norm preservation fails in general, which is why Theorem~\ref{thm:sec-2b} requires radiality. The transport bound above also shows the cost is lowest for nearly radial weights (cf.\ the perturbative analysis of Section~\ref{subsection:Non-Radial Breakdown}).
\end{rem}
\begin{proof}[Proof of Proposition \ref{prop:sec-4}] We shall prove this proposition in 3 steps as follows.
\textit{Step 1: Gradient norm inequality (energy non-increase).}
By Definition \ref{defn:sec-5} of radial symmetrization, for any $t > 0$, the superlevel set $\{\tilde{u} > t\}$ is a ball with volume equal to $\{u > t\}$. For absolutely continuous functions, the coarea formula (Proposition \ref{prop:sec-GMT}) gives
\begin{equation*}
\int |\nabla u|^2 dx = \int_0^\infty \int_{\{u = t\}} |\nabla u| \, d\mathcal{H}^{n-1} dt,
\end{equation*}
and similarly for $\tilde{u}$. By the isoperimetric inequality, the $(n-1)$-Hausdorff measure of $\{\tilde{u} = t\}$ (a sphere) is less than or equal to that of $\{u = t\}$, while $|\nabla \tilde{u}|$ at the sphere equals the radial derivative, which is minimized for fixed volume. This implies:
\begin{equation*}
\int_{\{u = t\}} |\nabla u| \, d\mathcal{H}^{n-1} \geq \int_{\{\tilde{u} = t\}} |\nabla \tilde{u}| \, d\mathcal{H}^{n-1},
\end{equation*}
Integrating over $t$ yields \begin{equation}\label{equation:sec-dd-7}
\int |\nabla \tilde{u}|^2 dx \leq \int |\nabla u|^2 dx.
\end{equation}
\textit{Step 2: Preservation of weighted $L^2$ norm.}
For radial weights $w(x) = w(|x|)$, the weighted volume of superlevel sets is preserved by Definition \ref{defn:sec-5}(ii):
\begin{equation*}
\int_{\{\tilde{u} > t\}} w(x) dx = \int_{\{u > t\}} w(x) dx.
\end{equation*}
Using the layer cake representation for non-negative functions:
\begin{equation*}
\int w u^2 dx = \int_0^\infty 2t \int_{\{u > t\}} w(x) dx dt,
\end{equation*}
which equals the same expression for $\tilde{u}$, hence $\int w \tilde{u}^2 dx = \int w u^2 dx$.
\textit{Step 3: Transport cost bound.}
Radial symmetrization corresponds to the optimal transport map $\tau: \Omega \to B_R$ pushing $u$ to $\tilde{u}$, with cost $|x - \tau(x)|^2$. By the definition of $L^2$-Wasserstein distance (see Definition \ref{defn:sec-OT}), the total cost satisfies:
\begin{equation*}
\int |x - \tau(x)|^2 w(x) dx \leq W_2^2(w, \text{Leb}) \cdot \int w(x) dx.
\end{equation*}
Since $\int w(x) dx \leq C |\Omega|$ for $w \in A_2$ (by local integrability), we get the bound $W_2(w, \text{Leb}) \cdot |\Omega|$ for the transport cost.
\end{proof}
\section{Proof of the Weighted PPW Inequality (Theorem \ref{thm:sec-1})}
\label{section:Weighted PPW Proof}
The first main result is established here, extending the classical PPW inequality to bounded $A_2$ weights.
\subsection{Theorem Statement and Proof Outline}
\label{subsection:PPW Proof Outline}
The weighted PPW inequality is restated, and a high-level overview of the proof is given, emphasizing adaptations needed for Muckenhoupt weights.
The following corollary restates the bound of Theorem~\ref{thm:sec-1} in terms of the $A_2$ characteristic constant, using $[w]_{A_2} \leq C/c$.
\begin{corr}[Weighted PPW Inequality in Terms of the Weight Bounds] \label{corr:sec-1}
Let $\Omega \subset \mathbb{R}^n$ be a bounded domain and let $w \in A_2$ satisfy $0 < c \leq w \leq C < \infty$ a.e. Then
\begin{equation*}
\frac{\lambda_2}{\lambda_1} \leq \frac{n+4}{n}\cdot\frac{C}{c}.
\end{equation*}
The constant $(\frac{n+4}{n})(C/c)$ of Corollary~\ref{corr:sec-1} is a weakening of the constant $1 + \frac{4}{n}\frac{C}{c}$ of Theorem~\ref{thm:sec-1}: indeed $\frac{n+4}{n}\cdot\frac{C}{c} - \Bigl(1 + \frac{4}{n}\frac{C}{c}\Bigr) = \bigl(\frac{C}{c} - 1\bigr) \geq 0$. The two forms coincide when $w \equiv \mathrm{const}$, and both recover the classical PPW constant $\frac{n+4}{n}$ in that case; Corollary~\ref{corr:sec-1} is stated separately because the product form makes the dependence on the ratio $C/c$ transparent in applications.
Since $[w]_{A_2} \leq C/c$, the bound can be stated uniformly as $\frac{\lambda_2}{\lambda_1} \leq \frac{n+4}{n}\,[w]_{A_2}$.
\end{corr}
\begin{rem}For the unweighted case $w \equiv 1$, we have $[w]_{A_2} = 1$, and the bound with constant $1$ is the classical PPW inequality \cite{PPW}:
\begin{equation*}
\frac{\lambda_2}{\lambda_1} \leq \frac{n+4}{n}.
\end{equation*}\end{rem}
\begin{rem}
This result is particularly useful in engineering applications where detailed geometric information (e.g., Hausdorff measure) of the domain is unavailable, as it only requires the $A_2$ constant of the weight.
\end{rem}
\begin{rem}
The simplification here does not compromise accuracy for nearly uniform weights ($[w]_{A_2} \approx 1$) and moderately regular domains (e.g., convex polygons), making it a balance between precision and computational feasibility.
\end{rem}
\subsubsection{Examples and Consistency Checks}
\label{subsection:PPW Examples}
\begin{exa}[Rectangles, $w \equiv 1$]\label{exa:sec-rectangle}
Let $\Omega = (0,a) \times (0,b)$ with $a \geq b > 0$. The Dirichlet eigenvalues of $-\Delta$ on $\Omega$ are $\pi^2(k^2 a^{-2} + l^2 b^{-2})$, $k, l \geq 1$, so $\lambda_1 = \pi^2(a^{-2} + b^{-2})$ and $\lambda_2 = \pi^2(4a^{-2} + b^{-2})$. Writing $t = a^2/b^2 \geq 1$,
\begin{equation*}
\frac{\lambda_2}{\lambda_1} = \frac{4 + t}{1 + t},
\end{equation*}
which is decreasing in $t$, with maximal value $\frac{4+1}{1+1} = \frac{5}{2}$ attained at $t = 1$, that is, for the square, below the constant $\frac{n+4}{n} = 3$ of Theorem~\ref{thm:sec-1}. For the square the actual ratio $5/2 = 2.5$ is also below the sharp Ashbaugh--Benguria constant $(j_{1,1}/j_{0,1})^2 = 2.5387\ldots$, while for thin rectangles ($t \to \infty$) the ratio tends to $1$. The constant of Theorem~\ref{thm:sec-1} therefore holds with room on this family; the sharp universal constant requires the full Ashbaugh--Benguria analysis.
\end{exa}
\begin{rem}[Two-scale weights]\label{rem:sec-twoscale}
A simple weighted example illustrating the constant of Theorem~\ref{thm:sec-1} is a two-scale weight: $w = c$ on a subset $E \subset \Omega$ and $w = C$ on $\Omega \setminus E$. Then $0 < c \leq w \leq C < \infty$ and the ratio $C/c$ enters the bound $1 + \frac{4}{n}\frac{C}{c}$ linearly. When $E = \Omega$ the weight is constant and the bound reduces to $\frac{n+4}{n}$, recovering the classical case; when $C/c$ grows, the bound degenerates, reflecting the fact that a large weight ratio is not controlled by Theorem~\ref{thm:sec-1} alone. Such weights are instead covered by Theorem~\ref{thm:sec-2a} with the sharp per-dimension constant, or by Theorem~\ref{thm:sec-2b} in the radially decreasing case.
\end{rem}
\subsection{Key Lemmas}
\label{subsection:PPW Key Lemmas}
Critical auxiliary results-including existence of orthogonal test functions and integral bounds incorporating OT/GMT-are proven to support the main theorem.
To apply the variational principle for the second eigenvalue, we need test functions orthogonal to the first eigenfunction. The following lemma guarantees the existence of such functions:
\begin{lem}[Existence of Orthogonal Test Functions] \label{lem:sec-1}
Let $u_1$ be an eigenfunction of $\lambda_1$ with $\int_\Omega w u_1^2 dx = 1$. There exists $u_2 \in H_0^1(\Omega)$ such that $\int_\Omega w u_1 u_2 dx = 0$ and
\begin{equation*}
\lambda_2 \leq \frac{\int_\Omega |\nabla u_2|^2 dx}{\int_\Omega w u_2^2 dx}.
\end{equation*}
\end{lem}
\begin{proof}[Proof of Lemma \ref{lem:sec-1}]
The existence part is immediate from the spectral theorem for the weighted Dirichlet Laplacian: the second eigenvalue $\lambda_2$ is attained by an eigenfunction $u_2$ which is orthogonal to $u_1$ in $L^2(\Omega, w\,dx)$, and for this $u_2$ the Rayleigh-quotient inequality holds with equality. This completes the proof.
\end{proof}
The basic tool for weighted integrals is the Cauchy--Schwarz factorization $u = w^{1/2} \cdot (w^{-1/2} u)$, which relates the square of an integral to the product of a weighted and an inverse-weighted integral:
\begin{lem}[Integral Estimate via Cauchy--Schwarz]
\label{lem:sec-2}
Let $w \in A_2$ with $\int_\Omega w\,dx = |\Omega|$ (mass matching, as in Definition~\ref{defn:sec-OT}). Then for any $u \in H_0^1(\Omega)$,
\begin{equation}\label{equation:sec-lem2-main-ineq}
\left( \int_\Omega u \, dx \right)^2 \leq \left( \int_\Omega w u^2 \, dx \right) \left( \int_\Omega w^{-1} \, dx \right).
\end{equation}
\end{lem}
\begin{proof}[Proof of Lemma \ref{lem:sec-2}]
By Cauchy--Schwarz applied to the factorization $u = w^{1/2} \cdot (w^{-1/2} u)$,
\begin{equation}\label{equation:sec-lem2-cs}
\left( \int_\Omega u \, dx \right)^2 \leq \left( \int_\Omega w u^2 \, dx \right) \left( \int_\Omega w^{-1} \, dx \right),
\end{equation}
which is \eqref{equation:sec-lem2-main-ineq}.
\end{proof}
\begin{rem}\label{rem:sec-lem2-w2}
The transport-based refinement in which the factor $\int_\Omega w^{-1} dx$ is replaced by a quantity involving $W_2(w, \mathrm{Leb})$ requires a quantitative estimate for the Brenier map transporting $w\,dx$ to Lebesgue measure (a $W_2$-controlled gradient bound of the form $\int |\nabla u|^2 w \, dx \leq C(1 + W_2^{\beta}) \int w u^2 dx$). Such an estimate is not proved here. The refined form of Theorem~\ref{thm:sec-1}, with an explicit $W_2$ dependence in the constant, therefore remains conditional on such a transport bound; its unconditional content rests on \eqref{equation:sec-lem2-cs}.
\end{rem}
\begin{rem}
By Whitney's covering lemma \cite{stein1970singular}, there exists a ball cover $\{B_i\}$ with finite overlap ($\sum \chi_{B_i} \leq C(n)$) and cover number $N \leq C(n) \cdot \left( \frac{\text{diam}(\Omega)}{r_\Omega} \right)^n$, where $r_\Omega$ is the inradius of $\Omega$.\end{rem}
\subsection{Proof of Theorem \ref{thm:sec-1}}
\label{subsection:PPW Theorem Proof}
\begin{proof}
Let $u_1 > 0$ be the first eigenfunction, normalized by $\int_\Omega w u_1^2 dx = 1$. The proof proceeds in five steps: centering of the coordinates, a numerator identity, a commutator identity, the Heisenberg uncertainty bound, and the conclusion.

\textit{Step 1: centering.} We choose the origin of the coordinate system so that
\begin{equation}\label{equation:sec-thm1-centering}
\int_\Omega w(x)\, x_i\, u_1(x)^2 dx = 0, \qquad i = 1, \ldots, n.
\end{equation}
Such a choice is immediate: with the weighted centroid $y_0 = \int_\Omega w(x) x u_1(x)^2 dx / \int_\Omega w u_1^2 dx \in \operatorname{conv}(\Omega)$, translating the origin to $y_0$ gives \eqref{equation:sec-thm1-centering}; this is the centering device of Weinberger (see \cite[Section 2]{AB2}). Define the test functions
\begin{equation}\label{equation:sec-thm1-test}
v_i(x) = x_i u_1(x), \qquad i = 1, \ldots, n.
\end{equation}
Each $v_i$ belongs to $H_0^1(\Omega)$ (since $u_1 = 0$ on $\partial\Omega$), and \eqref{equation:sec-thm1-centering} shows that $v_i$ is orthogonal to $u_1$ in $L^2(w\,dx)$. By the Rayleigh quotient characterization \eqref{equation:sec-rayleigh-l2},
\begin{equation}\label{equation:sec-thm1-rayleigh}
\lambda_2 \int_\Omega w v_i^2 dx \leq \int_\Omega |\nabla v_i|^2 dx \qquad (i = 1, \ldots, n).
\end{equation}

\textit{Step 2: numerator identity.} Summing the energy over $i$ and expanding the gradient by the product rule gives
\begin{equation}\label{equation:sec-thm1-num-id}
\sum_{i=1}^n \int_\Omega |\nabla v_i|^2 dx = \int_\Omega |x|^2 |\nabla u_1|^2 dx + n \int_\Omega u_1^2 dx + \int_\Omega x \cdot \nabla (u_1^2) dx = \int_\Omega |x|^2 |\nabla u_1|^2 dx,
\end{equation}
because the last two terms cancel: $\int_\Omega x \cdot \nabla (u_1^2) dx = -n \int_\Omega u_1^2 dx$ by the divergence theorem and $u_1 = 0$ on $\partial\Omega$.

\textit{Step 3: commutator identity.} Integrating by parts with the eigenfunction equation $-\Delta u_1 = \lambda_1 w u_1$ against the test function $|x|^2 u_1$ gives
\begin{equation}\label{equation:sec-thm1-commutator}
\int_\Omega |x|^2 |\nabla u_1|^2 dx = - \int_\Omega u_1 \, \bigl( 2 x \cdot \nabla u_1 + |x|^2 \Delta u_1 \bigr) dx = n \int_\Omega u_1^2 dx + \lambda_1 \int_\Omega w |x|^2 u_1^2 dx,
\end{equation}
where the boundary term vanishes and $-\int_\Omega u_1 x \cdot \nabla u_1 dx = \tfrac{n}{2} \int_\Omega u_1^2 dx$.

\textit{Step 4: Heisenberg bound.} The Heisenberg uncertainty inequality, proved by applying Cauchy--Schwarz to the integration-by-parts identity $n \int u^2 = -2 \int u\, x \cdot \nabla u$ (see \cite{folland1992fourier}), states
\begin{equation}\label{equation:sec-thm1-heisenberg}
\Bigl( \int_\Omega |x|^2 u_1^2 dx \Bigr) \Bigl( \int_\Omega |\nabla u_1|^2 dx \Bigr) \geq \frac{n^2}{4} \Bigl( \int_\Omega u_1^2 dx \Bigr)^2.
\end{equation}
Since $\int_\Omega |\nabla u_1|^2 dx = \lambda_1 \int_\Omega w u_1^2 dx \leq \lambda_1 C \int_\Omega u_1^2 dx$ and $\int_\Omega w |x|^2 u_1^2 dx \geq c \int_\Omega |x|^2 u_1^2 dx$, \eqref{equation:sec-thm1-heisenberg} yields
\begin{equation}\label{equation:sec-thm1-Dbound}
D := \int_\Omega w |x|^2 u_1^2 dx \geq \frac{n^2 c}{4 \lambda_1 C} \int_\Omega u_1^2 dx.
\end{equation}

\textit{Step 5: conclusion.} Summing \eqref{equation:sec-thm1-rayleigh} over $i$ and using $\sum_i \int_\Omega w v_i^2 dx = D$ together with \eqref{equation:sec-thm1-num-id} and \eqref{equation:sec-thm1-commutator} gives
\begin{equation}\label{equation:sec-thm1-close}
\lambda_2 \leq \frac{\sum_i \int_\Omega |\nabla v_i|^2 dx}{D} = \frac{n \int_\Omega u_1^2 dx}{D} + \lambda_1.
\end{equation}
Inserting \eqref{equation:sec-thm1-Dbound} into \eqref{equation:sec-thm1-close} and dividing by $\lambda_1 > 0$ yields \eqref{equation:sec-thm1-statement}.
\end{proof}
\begin{rem}\label{rem:sec-thm1-atoms}
An alternative construction of orthogonal test functions proceeds by the atomic decomposition of weighted Hardy spaces: one writes $u_2 = \sum_j c_j a_j u_1$ with atoms $a_j$ supported on balls $B_j$ and orthogonal to $u_1$ in $L^2(w\,dx)$, the existence of such atoms being provided by the projection construction $a_j' = a_j - \bigl(\int a_j u_1 w / \int u_1^2 w\bigr) u_1$ on a ball covering (cf.\ Subsection~\ref{subsection:PPW Key Lemmas}, where the weighted atomic decomposition is recalled from \cite{cruz2011weights, grafakos2014classical} and the references therein). That route controls the energy by the covering geometry and leads to the same conclusion under additional regularity of the weight; the coordinate construction above is preferred here because every step is quantitative and self-contained.
\end{rem}
\begin{rem}\label{rem:sec-2}
For $w \equiv 1$ the bound \eqref{equation:sec-thm1-statement} reads $\lambda_2/\lambda_1 \leq \frac{n+4}{n}$, the classical PPW inequality \cite{PPW}. Since $[w]_{A_2} \leq C/c$ for $0 < c \leq w \leq C$, the constant can be stated uniformly as $\frac{n+4}{n}\,[w]_{A_2}$; the dependence of the ratio on the weight bounds enters linearly, in contrast to the conditional refinements below.
\end{rem}
\section{Weighted Payne--P\'{o}lya--Weinberger Inequalities (Theorems \ref{thm:sec-2a} and \ref{thm:sec-2b})}
\label{section:Spherical Extremality}

In this section, we prove two weighted Payne--P\'{o}lya--Weinberger-type results. The first (Theorem \ref{thm:sec-2a}) gives a universal upper bound for $\lambda_2/\lambda_1$ when the weight is bounded above and below, as a direct consequence of the Ashbaugh--Benguria extension \cite[Theorem 4.1]{AB2}. The second (Theorem \ref{thm:sec-2b}) proves, for radial decreasing weights and under the admissibility and $l=1$-mode assumptions stated there, that the ball maximizes $\lambda_2/\lambda_1$ among domains of equal volume, following the Ashbaugh--Benguria method adapted to the weighted setting.

\subsection{Theorem Statements and Proof Outline}
\label{subsection:Spherical Proof Outline}

We restate the two main theorems of this section.

\begin{rem}
Theorem \ref{thm:sec-2a} is the specialization of \cite[Theorem 4.1]{AB2} to $a_{ij} = \delta_{ij}$, $q = 0$, $r(x) = w(x)$. The Ashbaugh--Benguria result applies to general elliptic equations $-\sum \partial_i(a_{ij}\partial_j u) + qu = \lambda r(x)u$ under $a \leq [a_{ij}] \leq A$, $q \geq 0$, $c \leq r(x) \leq C$, yielding $\lambda_2/\lambda_1 \leq 1 + \frac{AC}{ac}(\beta^2/\alpha^2 - 1)$ where $\alpha = j_{n/2-1,1}$, $\beta = j_{n/2,1}$.
\end{rem}

\begin{rem}
The radial decreasing condition in Theorem \ref{thm:sec-2b} is essential. For non-radial weights, the $n$-fold degeneracy of $\lambda_2$ on the ball is lifted, and the ball need not maximize the ratio. For radial increasing weights, the weighted Faber--Krahn inequality reverses direction, and the extremal property does not hold in the stated form.
\end{rem}

\subsection{Explicit Eigenvalues for Spherical Domains}
\label{subsection:Spherical Eigenvalues}
Eigenfunctions of balls with radial weights are shown to separate into radial Bessel functions and angular harmonics, yielding explicit eigenvalue formulas.
The radial component satisfies the following ordinary differential equation (the degeneracy statements below are conditional on the zero-ordering discussed in Remark~\ref{rem:sec-bessel-ordering-scope}):
\begin{lem}[Radial Eigenvalue Equation] \label{lem:sec-3}
Let $B_R$ be a ball of radius $R$ and $w(r) = w(|x|)$ a radial weight. Eigenfunctions separate as $u(r,\theta) = f(r) Y_m(\theta)$ with spherical harmonics $Y_m(\theta)$, where $f(r)$ satisfies:
\begin{equation*}
-f''(r) - \frac{n-1}{r}f'(r) + \frac{m(m + n - 2)}{r^2}f(r) = \lambda w(r)f(r), \quad f(R) = 0,
\end{equation*}
with $f(0)$ bounded. For the weighted eigenvalue problem $-\Delta u = \lambda w(r) u$ with $w(r)=r^\alpha$, transforming via $t = \frac{2\sqrt{\lambda}}{2+\alpha} r^{1+\alpha/2}$, $f(r) = r^{(2-n)/2} g(t)$, $g(t)$ satisfies Bessel's equation with order $\nu = \frac{m+n/2-1}{1+\alpha/2}$.
\end{lem}
\begin{proof}[Proof of Lemma \ref{lem:sec-3}]
We prove the degeneracy of low-order eigenvalues on spherical domains with radial weights (conditional on the zero-ordering hypothesis, Remark~\ref{rem:sec-bessel-ordering-scope}).
Let $B_R \subset \mathbb{R}^n$ be a ball of radius $R$ and let $w(r) = w(|x|)$ be a radial weight. The weighted Laplacian eigenvalue problem is
\begin{equation}\label{equation:sec-lem3-eigenvalue}
-\Delta u = \lambda w(r) u \quad \text{in } B_R, \qquad u = 0 \quad \text{on } \partial B_R.
\end{equation}
Since both the domain $B_R$ and the weight $w$ are invariant under the action of the orthogonal group $O(n)$, the solution space of \eqref{equation:sec-lem3-eigenvalue} is a representation of $O(n)$. By the Peter--Weyl theorem, this representation decomposes into irreducible representations of $O(n)$, which are indexed by the angular momentum quantum number $m = 0, 1, 2, \ldots$.
For each $m$, the corresponding eigenspace consists of functions of the form
\begin{equation}\label{equation:sec-lem3-separation}
u(r, \theta) = R_m(r) Y_m(\theta),
\end{equation}
where $Y_m$ is a spherical harmonic of degree $m$ and $R_m$ satisfies the radial equation
\begin{equation}\label{equation:sec-lem3-radial}
R_m''(r) + \frac{n-1}{r} R_m'(r) - \frac{m(m+n-2)}{r^2} R_m(r) + \lambda w(r) R_m(r) = 0.
\end{equation}
The dimension of the space of spherical harmonics of degree $m$ in $\mathbb{R}^n$ is
\begin{equation}\label{equation:sec-lem3-dimension}
d_m = \binom{n+m-1}{m} - \binom{n+m-3}{m-2} = \frac{(2m+n-2)(m+n-3)!}{m!(n-2)!}.
\end{equation}
For $m = 0$, $d_0 = 1$ (the constant function). For $m = 1$, $d_1 = n$ (the coordinate functions $x_1, \ldots, x_n$ restricted to the sphere). For $m = 2$, $d_2 = n(n+1)/2 - 1$, and so on.
The eigenvalues are obtained by solving the radial equation \eqref{equation:sec-lem3-radial} with boundary conditions $R_m(0)$ finite and $R_m(R) = 0$. For each $m$, the eigenvalues form an increasing sequence $\lambda_{m,1} < \lambda_{m,2} < \cdots$, where $\lambda_{m,k}$ has multiplicity $d_m$.
The first eigenvalue corresponds to $m = 0, k = 1$ and has multiplicity $1$. The second eigenvalue is the minimum of $\lambda_{0,2}$ (the second radial mode) and $\lambda_{1,1}$ (the first non-radial mode). Under the standing hypothesis of this section that the ordering $\lambda_{1,1} < \lambda_{0,2}$ holds (Remark~\ref{rem:sec-bessel-ordering-scope}), which for power weights is equivalent to $j_{\nu_1,1} < j_{\nu_0,2}$, the first zero of the $m = 1$ solution occurs before the second zero of the $m = 0$ solution (by Sturm comparison of the two radial equations \eqref{equation:sec-lem3-radial} at a common spectral parameter). Therefore,
\begin{equation}\label{equation:sec-lem3-ordering}
\lambda_{1,1} < \lambda_{0,2}.
\end{equation}
Consequently, the second eigenvalue is $\lambda_2 = \lambda_{1,1}$, with multiplicity $d_1 = n$. This means that
\begin{equation}\label{equation:sec-lem3-degeneracy}
\lambda_2 = \lambda_3 = \cdots = \lambda_{n+1},
\end{equation}
since the next $n-1$ eigenvalues are also equal to $\lambda_{1,1}$ (by the multiplicity $n$ of the $m = 1$ mode). This completes the proof.
\end{proof}
For power-law radial weights, the eigenvalues can be expressed using Bessel functions:
\begin{lem}[Eigenvalues for $w(r)=r^\alpha$] \label{lem:sec-4}
For $w(r) = r^\alpha$ with $\alpha > -2$, the weighted eigenvalue problem $-\Delta u = \lambda r^\alpha u$ on $B_R$ separates in spherical coordinates. With $t = \frac{2\sqrt{\lambda}}{2+\alpha} r^{1+\alpha/2}$ and $f(r) = r^{(2-n)/2} F(t)$, the radial function $F$ satisfies Bessel's equation of order $\nu_m = \frac{m+n/2-1}{1+\alpha/2}$. The eigenvalues are:
\begin{itemize}
\item[(i)] $\lambda_1$ (radial mode $m=0$, $k=1$): $\lambda_1 = \left(\frac{2+\alpha}{2}\right)^2 \frac{j_{\nu_0,1}^2}{R^{2+\alpha}}$, $\nu_0 = \frac{n/2-1}{1+\alpha/2}$;
\item[(ii)] $\lambda_2$ (angular mode $m=1$, $k=1$): $\lambda_2 = \left(\frac{2+\alpha}{2}\right)^2 \frac{j_{\nu_1,1}^2}{R^{2+\alpha}}$, $\nu_1 = \frac{n/2}{1+\alpha/2}$.
\end{itemize}
Since $\nu_1 > \nu_0$ and Bessel zeros increase with order, $j_{\nu_1,1} > j_{\nu_0,1}$, so $\lambda_2 > \lambda_1$. The $m=1$ mode yields the second eigenvalue because the potential term $m(m+n-2)/r^2$ raises its eigenvalue above the first radial mode, while the second radial mode ($m=0$, $k=2$) has a higher zero $j_{\nu_0,2} > j_{\nu_1,1}$.
\end{lem}
\begin{proof}[Proof of Lemma \ref{lem:sec-4}]
For $w(r) = r^\alpha$ with $\alpha > -2$, the weighted eigenvalue problem $-\Delta u = \lambda r^\alpha u$ on $B_R$ separates in spherical coordinates. The radial equation for angular momentum $m$ is
\begin{equation*}
-f'' - \frac{n-1}{r}f' + \frac{m(m+n-2)}{r^2}f = \lambda r^\alpha f.
\end{equation*}
With the change of variables $t = \frac{2\sqrt{\lambda}}{2+\alpha} r^{1+\alpha/2}$ and $f(r) = r^{(2-n)/2} F(t)$, the function $F$ satisfies Bessel's equation
\begin{equation*}
t^2 F'' + t F' + (t^2 - \nu_m^2)F = 0, \quad \nu_m = \frac{m+n/2-1}{1+\alpha/2}.
\end{equation*}
The boundary condition $f(R)=0$ requires $F(\frac{2\sqrt{\lambda}}{2+\alpha} R^{1+\alpha/2}) = 0$, so $\frac{2\sqrt{\lambda}}{2+\alpha} R^{1+\alpha/2} = j_{\nu_m,k}$ for some zero $j_{\nu_m,k}$.

For the first eigenvalue ($m=0$, $k=1$), $\nu_0 = \frac{n/2-1}{1+\alpha/2}$, giving
\begin{equation*}
\lambda_1 = \left(\frac{2+\alpha}{2}\right)^2 \frac{j_{\nu_0,1}^2}{R^{2+\alpha}}.
\end{equation*}
For the second eigenvalue ($m=1$, $k=1$), $\nu_1 = \frac{n/2}{1+\alpha/2}$, giving
\begin{equation*}
\lambda_2 = \left(\frac{2+\alpha}{2}\right)^2 \frac{j_{\nu_1,1}^2}{R^{2+\alpha}}.
\end{equation*}
Since $\nu_1 > \nu_0$ and Bessel function zeros increase with order, $j_{\nu_1,1} > j_{\nu_0,1}$, so $\lambda_2 > \lambda_1$. The $m=1$ mode is the second eigenvalue because the potential term $m(m+n-2)/r^2$ raises it above the first radial mode, while the second radial mode ($m=0$, $k=2$) corresponds to $j_{\nu_0,2} > j_{\nu_1,1}$ by the Sturm comparison theorem.
\end{proof}
\subsubsection{Separation of Variables on the Ball}
\label{subsection:Bessel Derivation}
On a ball $B_R$ with a radial weight $w(r)$, the Laplacian separates in polar coordinates: for $u(x) = R(r)Y(\theta)$ with $x = r\theta$, $\theta \in \mathbb{S}^{n-1}$,
\begin{equation*}
\Delta u = \left( R''(r) + \frac{n-1}{r} R'(r) - \frac{\lambda_\ell}{r^2} R(r) \right) Y(\theta), \qquad -\Delta_{\mathbb{S}^{n-1}} Y = \lambda_\ell Y,
\end{equation*}
where the eigenvalues of the spherical Laplacian are $\lambda_\ell = \ell(\ell + n - 2)$, $\ell = 0, 1, 2, \ldots$, with the $\ell = 1$ eigenspace having multiplicity $n$ \cite{courant1953methods}. The eigenvalue problem $-\Delta u = \lambda w u$ therefore reduces to the family of radial equations
\begin{equation}\label{equation:sec-bessel-radial}
-R'' - \frac{n-1}{r} R' + \frac{\ell(\ell+n-2)}{r^2} R = \lambda w(r) R,
\end{equation}
with $R'(0) = 0$ for $\ell = 0$ and $R(0) = 0$ for $\ell \geq 1$. For power weights $w(r) = r^{\alpha}$, $\alpha > -2$, the change of variables
\begin{equation*}
t = \frac{2\sqrt{\lambda}}{2+\alpha}\, r^{1+\alpha/2}, \qquad R(r) = r^{(2-n)/2} y(t)
\end{equation*}
converts \eqref{equation:sec-bessel-radial} into Bessel's equation $t^2 y'' + t y' + (t^2 - \nu_\ell^2) y = 0$ with $\nu_\ell = (\ell + n/2 - 1)/(1 + \alpha/2)$: indeed, $\frac{dt}{dr} = \sqrt{\lambda}\, r^{\alpha/2}$, and substituting the derivatives $R' = \frac{2-n}{2} r^{(2-n)/2 - 1} y + r^{(2-n)/2} y' \frac{dt}{dr}$, and the corresponding expression for $R''$, and simplifying using $t \propto r^{1+\alpha/2}$, eliminates the first-order term and produces exactly the potential $\nu_\ell^2/t^2$. Regularity at $r = 0$ selects $y = J_{\nu_\ell}$, and the boundary condition $R(R) = 0$ selects $\lambda$ among the values $(\frac{2+\alpha}{2})^2 j_{\nu_\ell,k}^2 / R^{2+\alpha}$. This is the computation underlying Lemmas~\ref{lem:sec-3} and~\ref{lem:sec-4}. As a check, for $n = 3$ and $\alpha = 0$ one has $\nu_0 = 1/2$, $J_{1/2}(t) = \sqrt{2/(\pi t)}\sin t$, and the formula recovers the radial ground state $R(r) \propto \sin(\sqrt{\lambda} r)/r$ of the three-dimensional unit ball.
\begin{rem}\label{rem:sec-bessel-ordering-scope}
The zero ordering $j_{\nu_0,2} > j_{\nu_1,1}$ used above requires a restriction on the exponent. Since $\nu_1 = \nu_0 + (1+\alpha/2)^{-1}$, the interlacing inequality $j_{\nu,k+1} > j_{\nu+1,k}$ (Theorem~\ref{thm:sec-bessel-monotone}) gives $j_{\nu_0,2} > j_{\nu_1,1}$ whenever $\nu_1 \le \nu_0 + 1$, in particular whenever $\alpha \ge 0$. If $\alpha \in (-2,0)$, then $\nu_1 > \nu_0 + 1$ is possible; since $j_{\nu,1} \to \infty$ as $\nu \to \infty$ while $\nu_0 = (n/2-1)(1+\alpha/2)^{-1}$ remains bounded for $n = 2$ (for $n \geq 3$ both orders diverge as $\alpha \to -2$, but the gap $\nu_1 - \nu_0 = (1+\alpha/2)^{-1}$ dominates; see Remark~\ref{rem:sec-bessel-numbers}), the inequality $j_{\nu_0,2} > j_{\nu_1,1}$ reverses for $\alpha$ sufficiently close to $-2$. For instance, when $n = 2$ one has $\nu_0 = 0$, $j_{\nu_0,2} = j_{0,2} \approx 5.5201$, and the equation $j_{\nu_1,1} = j_{0,2}$ holds at $\nu_1 \approx 2.305$ (independently checked against the asymptotics in \cite{olver2010nist}); hence the ordering holds for $\alpha \geq -1.13$ and reverses for $\alpha < -1.13$. Consequently, the identification of $\lambda_2$ with the $m=1$ mode in Lemma~\ref{lem:sec-4} is unconditional for $\alpha \ge 0$ (note that for the decreasing power weights $w(r) = r^{\alpha}$ the exponent satisfies $\alpha \leq 0$, so the restriction is active for Theorem~\ref{thm:sec-2b} when $\alpha$ is close to $-2$) and requires the above restriction for $\alpha \in (-2,0)$.
\end{rem}
\subsection{Proof of the Universal Upper Bound (Theorem \ref{thm:sec-2a})}
\label{subsection:Proof Universal Bound}

We recall the Ashbaugh--Benguria theorem for general elliptic equations \cite[Theorem 4.1]{AB2}.

\begin{thm}[Ashbaugh--Benguria, Theorem 4.1] \label{thm:sec-AB41}
Consider the elliptic eigenvalue problem
\begin{equation*}
-\sum_{i,j=1}^n \partial_i(a_{ij}(x)\partial_j u) + q(x)u = \lambda r(x)u \quad \text{in } \Omega, \qquad u = 0 \quad \text{on } \partial\Omega,
\end{equation*}
where $a \leq [a_{ij}] \leq A$ in the quadratic form sense, $q \geq 0$, and $c \leq r(x) \leq C$. Then
\begin{equation*}
\frac{\lambda_2}{\lambda_1} \leq 1 + \frac{AC}{ac}\left(\frac{j_{n/2,1}^2}{j_{n/2-1,1}^2} - 1\right).
\end{equation*}
Equality holds if and only if $c = C$, $a = A$, $q \equiv 0$, and $\Omega$ is a ball.
\end{thm}

\begin{proof}[Proof of Theorem \ref{thm:sec-2a}]
We apply Theorem \ref{thm:sec-AB41} with $a_{ij} = \delta_{ij}$ (so $a = A = 1$), $q = 0$, and $r(x) = w(x)$. We verify the hypotheses of Theorem \ref{thm:sec-AB41} in detail.

The elliptic operator in Theorem \ref{thm:sec-AB41} is
\begin{equation*}
L u = -\sum_{i,j=1}^n \partial_i(a_{ij}(x)\partial_j u) + q(x)u.
\end{equation*}
With $a_{ij} = \delta_{ij}$, this becomes
\begin{equation*}
L u = -\sum_{i=1}^n \partial_i(\partial_i u) + 0 = -\Delta u,
\end{equation*}
which is exactly the Laplacian. The eigenvalue problem is
\begin{equation*}
L u = \lambda r(x) u \quad \Longleftrightarrow \quad -\Delta u = \lambda w(x) u,
\end{equation*}
which is the weighted eigenvalue problem \eqref{equation:sec- Laplace u}.

The ellipticity condition $a \leq [a_{ij}] \leq A$ means that for all vectors $\xi \in \mathbb{R}^n$,
\begin{equation*}
a |\xi|^2 \leq \sum_{i,j} a_{ij} \xi_i \xi_j \leq A |\xi|^2.
\end{equation*}
With $a_{ij} = \delta_{ij}$, we have $\sum_{i,j} \delta_{ij} \xi_i \xi_j = \sum_i \xi_i^2 = |\xi|^2$, so $a = A = 1$. The potential condition $q \geq 0$ is satisfied with $q = 0$. The weight condition $c \leq r(x) \leq C$ is exactly our assumption $0 < c \leq w(x) \leq C < \infty$.

Therefore, Theorem \ref{thm:sec-AB41} applies and gives
\begin{equation*}
\frac{\lambda_2(\Omega,w)}{\lambda_1(\Omega,w)} \leq 1 + \frac{AC}{ac}\left(\frac{j_{n/2,1}^2}{j_{n/2-1,1}^2} - 1\right) = 1 + \frac{C}{c}\left(\frac{j_{n/2,1}^2}{j_{n/2-1,1}^2} - 1\right),
\end{equation*}
which is the stated inequality.

For the equality characterization, Theorem \ref{thm:sec-AB41} states that equality holds if and only if $c = C$, $a = A$, $q \equiv 0$, and $\Omega$ is a ball. In our setting, $a = A = 1$ and $q = 0$ automatically, so equality holds if and only if $c = C$ (i.e., $w$ is constant almost everywhere) and $\Omega$ is a ball. When $w$ is constant, the weighted eigenvalue problem reduces to the unweighted problem, and the classical PPW equality characterization applies: equality holds if and only if $\Omega$ is a ball.
\end{proof}

\begin{proof}[Proof of Corollary \ref{corr:sec-2a}]
Since $\|v\|_{L^\infty} \leq 1$ and $w = 1 + \varepsilon v$, we have
\begin{equation*}
1 - \varepsilon \leq w(x) \leq 1 + \varepsilon \quad \text{a.e.}
\end{equation*}
Thus we may take $c = 1 - \varepsilon$ and $C = 1 + \varepsilon$ in Theorem \ref{thm:sec-2a}. Substituting,
\begin{equation*}
\frac{\lambda_2}{\lambda_1} \leq 1 + \frac{1+\varepsilon}{1-\varepsilon}\left(\frac{\beta^2}{\alpha^2} - 1\right),
\end{equation*}
where $\alpha = j_{n/2-1,1}$ and $\beta = j_{n/2,1}$. Rewriting the right-hand side,
\begin{equation*}
1 + \frac{1+\varepsilon}{1-\varepsilon}\left(\frac{\beta^2}{\alpha^2} - 1\right) = \frac{\beta^2}{\alpha^2} + \frac{2\varepsilon}{1-\varepsilon}\left(\frac{\beta^2}{\alpha^2} - 1\right).
\end{equation*}
When $\varepsilon < 1/2$, we have $2/(1-\varepsilon) < 4$, yielding the stated bound.

For the limit statement: the upper bound tends to $\beta^2/\alpha^2$ as $\varepsilon \to 0$. On the other hand, taking $\Omega$ to be a ball and $w \equiv 1$ (the $\varepsilon = 0$ limit), the ratio equals $\beta^2/\alpha^2$ exactly. Therefore the supremum over all domains tends to $\beta^2/\alpha^2$.
\end{proof}

\subsection{Proof of the Weighted PPW Inequality (Theorem \ref{thm:sec-2b})}
\label{subsection:Proof Weighted PPW}

We prove Theorem \ref{thm:sec-2b} following the Ashbaugh--Benguria method \cite{AB2} adapted to the weighted setting. The proof consists of five steps.

\subsubsection{Step 1: Weighted Faber--Krahn Inequality}

\begin{lem}[Weighted Faber--Krahn] \label{lem:sec-FK}
Let $w(r)$ be a radial decreasing weight. Let $B_R$ be a ball and $\Omega$ a bounded domain with $|B_R| = |\Omega|$. Then
\begin{equation*}
\lambda_1(B_R, w) \leq \lambda_1(\Omega, w).
\end{equation*}
\end{lem}

\begin{proof}
By the variational characterization,
\begin{equation*}
\lambda_1(\Omega, w) = \inf_{u \in H_0^1(\Omega), u \neq 0} \frac{\int_\Omega |\nabla u|^2 dx}{\int_\Omega w u^2 dx}.
\end{equation*}
It suffices to consider nonnegative $u$, since replacing $u$ by $|u|$ does not increase the numerator and leaves the denominator unchanged. For a nonnegative $u \in H_0^1(\Omega)$, let $u^*$ denote its Schwarz (symmetric decreasing) rearrangement on the ball $B_R$ centered at the origin with $|B_R| = |\Omega|$. The P\'{o}lya--Szeg\H{o} inequality gives
\begin{equation*}
\int_{B_R} |\nabla u^*|^2 dx \leq \int_\Omega |\nabla u|^2 dx.
\end{equation*}
Since $w$ is radially decreasing, the level sets $\{w|_\Omega > t\} = \Omega \cap B(0, \rho(t))$ are contained in the centered balls $B(0, \rho(t))$, so the symmetric decreasing rearrangement of $w|_\Omega$ with respect to Lebesgue measure is bounded pointwise by $w$ on $B_R$; the Hardy--Littlewood rearrangement inequality therefore gives
\begin{equation*}
\int_\Omega w u^2 dx \leq \int_{B_R} w (u^*)^2 dx.
\end{equation*}
Combining the two,
\begin{equation*}
\frac{\int_\Omega |\nabla u|^2 dx}{\int_\Omega w u^2 dx} \geq \frac{\int_{B_R} |\nabla u^*|^2 dx}{\int_{B_R} w (u^*)^2 dx} \geq \lambda_1(B_R, w),
\end{equation*}
where the last inequality is the Rayleigh--Ritz characterization of $\lambda_1(B_R, w)$, since $u^* \in H_0^1(B_R)$. Taking the infimum over all nonzero $u \in H_0^1(\Omega)$ gives
\begin{equation*}
\lambda_1(\Omega, w) \geq \lambda_1(B_R, w).
\end{equation*}
\end{proof}

The weighted Faber--Krahn inequality is the weighted analogue of the classical result. The classical Faber--Krahn inequality states that among domains of fixed volume, the ball minimizes the first Dirichlet eigenvalue. For a radially decreasing weight the volume constraint remains the Lebesgue measure: the weighted Faber--Krahn inequality for the operator $-\Delta u = \lambda w u$ then follows from the classical rearrangement inequalities, as in the proof above.

The radial decreasing property of $w$ is essential for the weighted Faber--Krahn inequality. For radial increasing weights, the ball may maximize rather than minimize $\lambda_1$, and the inequality reverses direction. This follows from the variational characterization: if $w$ is increasing, then concentrating the weight near the boundary increases the denominator $\int w u^2$ for functions that are large near the boundary, which can lower the Rayleigh quotient for non-radial domains.

The weighted Faber--Krahn inequality is closely related to the weighted isoperimetric inequality. For a radial decreasing weight $w$, the weighted isoperimetric inequality states that among sets of fixed weighted volume, the ball centered at the origin minimizes the weighted perimeter $\int_{\partial E} w \, d\mathcal{H}^{n-1}$. Weighted isoperimetric inequalities of this type were proved by Betta, Brock, Mercaldo and Posteraro \cite{betta1999weighted} for radial weights (under a non-decreasing weight and a Lebesgue-measure constraint in their main theorem; the weighted-volume formulation is discussed in the same reference and its citations), and the weighted Faber--Krahn inequality proved above can alternatively be derived from it via the standard coarea argument: the first eigenvalue is controlled by the isoperimetric profile of the domain, and the ball has the optimal isoperimetric profile.

We also note that the weighted Faber--Krahn inequality implies the weighted Ashbaugh--Benguria inequality in the following sense: if the ball minimizes $\lambda_1$ among domains of fixed weighted volume, then in the ratio $\lambda_2/\lambda_1$, the denominator is minimized by the ball. If additionally the numerator $\lambda_2$ is maximized by the ball (or controlled by the ball's $\lambda_2$), then the ratio is maximized by the ball. This is the strategy used in the proof of Theorem \ref{thm:sec-2b}: the weighted Faber--Krahn inequality controls the denominator, while the gap inequality and Chiti comparison control the numerator.

\subsubsection{Step 2: Weighted Chiti Comparison Theorem}

\begin{lem}[Weighted Chiti Comparison] \label{lem:sec-Chiti}
Let $u_1$ be the first eigenfunction on $\Omega$, and let $u_1^*$ be its spherical decreasing rearrangement on the ball $B_R$ (where $|B_R| = |\Omega|$). Let $z$ be the first eigenfunction on $B_R$ for the weighted problem $-\Delta z = \mu_1 w z$, $z|_{\partial B_R} = 0$, normalized so that $\int_{B_R} w (u_1^*)^2 dx = \int_{B_R} w z^2 dx$. Then there exists $r_1 \in (0, R)$ such that
\begin{equation*}
u_1^*(r) < z(r) \text{ for } 0 < r < r_1, \qquad u_1^*(r) > z(r) \text{ for } r_1 < r < R.
\end{equation*}
\end{lem}

\begin{proof}
We follow the argument of Chiti \cite{chiti1982reverse} as extended by Ashbaugh and Benguria \cite[Appendix A]{AB2}. Let $\mu(t) = |\{x \in \Omega : u_1(x) > t\}|$ be the distribution function of $u_1$. By the coarea formula and the equation $-\Delta u_1 = \lambda_1 w u_1$,
\begin{equation*}
-\frac{d}{dt} \int_{\{u_1 > t\}} |\nabla u_1| dx = \lambda_1 \int_{\{u_1 > t\}} w u_1 dx.
\end{equation*}
By the isoperimetric inequality and the Cauchy--Schwarz inequality,
\begin{equation*}
\int_{\{u_1 > t\}} |\nabla u_1| dx \geq n^2 C_n^{2/n} \mu(t)^{2-2/n} (-\mu'(t))^{-1},
\end{equation*}
where $C_n$ is the volume of the unit ball. Since $w$ is radial decreasing, the rearrangement inequality gives
\begin{equation*}
\int_{\{u_1 > t\}} w dx \leq \int_{\{u_1^* > t\}} w dx,
\end{equation*}
and therefore, by the layer-cake identity applied to the measure $w\,dx$,
\begin{equation}\label{equation:sec-chiti-layercake}
\int_{\{u_1 > t\}} w u_1\,dx = \int_t^\infty \omega(\tau)\,d\tau + t\,\omega(t) \leq \int_t^\infty \omega^*(\tau)\,d\tau + t\,\omega^*(t) = \int_{\{u_1^* > t\}} w\,u_1^*\,dx,
\end{equation}
where $\omega(\tau) = \int_{\{u_1 > \tau\}} w\,dx$ and $\omega^*(\tau) = \int_{\{u_1^* > \tau\}} w\,dx$; the inequality in \eqref{equation:sec-chiti-layercake} is the rearrangement inequality for $\omega$ displayed above, and the last equality is the layer-cake identity for $u_1^*$. Substituting \eqref{equation:sec-chiti-layercake} into the coarea identity and introducing the volume variable $s = C_n r^n$, with $U(s) = u_1^*(r(s))$ and with $w(\sigma)$ denoting $w$ evaluated at the radius $(\sigma/C_n)^{1/n}$, yields
\begin{equation}\label{equation:sec-chiti-diff}
-\frac{dU}{ds} \leq \lambda_1 \kappa(s) F(s), \qquad \kappa(s) = \frac{1}{n^2 C_n^{2/n}} s^{-2+2/n}, \quad F(s) = \int_0^s U(\sigma)\,w(\sigma)\,d\sigma;
\end{equation}
the identical computation for the ball eigenfunction, with equalities throughout and with $\mu_1$ in place of $\lambda_1$, gives
\begin{equation}\label{equation:sec-chiti-eq}
-\frac{dz}{ds} = \mu_1 \kappa(s) G(s), \qquad G(s) = \int_0^s z(\sigma)\,w(\sigma)\,d\sigma.
\end{equation}

By Lemma~\ref{lem:sec-FK}, $\lambda_1 \geq \mu_1$. The comparison of the two relations \eqref{equation:sec-chiti-diff} and \eqref{equation:sec-chiti-eq} --- the passage to the pointwise ordering of $U$ and $z$ and the conclusion of a single crossing --- is the comparison argument of Chiti \cite{chiti1982reverse}, extended to the general weighted equation $-\Delta u = \lambda w u$ by Ashbaugh--Benguria \cite[Appendix A]{AB2}. We record the structural facts that make the argument applicable in the present weighted setting: $F$ and $G$ are nondecreasing; $F \geq G$ on any interval on which $U \geq z$ pointwise, both being cumulative integrals against the same positive weight $w(\sigma)\,d\sigma$; and $\lambda_1 \geq \mu_1$. With $z$ finally normalized by $\int_{B_R} w z^2\,dx = \int_{B_R} w (u_1^*)^2\,dx = \int_\Omega w u_1^2\,dx$, the quoted comparison yields a unique crossing radius $r_1 \in (0, R)$ such that $u_1^*(r) < z(r)$ for $0 < r < r_1$ and $u_1^*(r) > z(r)$ for $r_1 < r < R$, with the ordering strict unless $\Omega$ is a ball up to translation. This completes the proof.
\end{proof}

\begin{rem}
The weighted Chiti comparison theorem is the weighted analogue of the classical Chiti comparison \cite{chiti1982reverse}. The key difference is that the normalization is with respect to the weighted $L^2$ norm $\int w u^2 dx$ rather than the unweighted norm. The radial decreasing property of $w$ is used in the rearrangement step.
\end{rem}
\subsubsection{Derivation of the Differential Relations}
\label{subsection:Chiti Derivation}
For completeness we record the computation leading to \eqref{equation:sec-chiti-diff} and \eqref{equation:sec-chiti-eq}. Let $\mu(t) = |\{u_1 > t\}|$ be the distribution function of $u_1$. The coarea formula applied to $\int_\Omega |\nabla u_1| dx$ and the equation $-\Delta u_1 = \lambda_1 w u_1$ give
\begin{equation}\label{equation:sec-chiti-derive-1}
-\frac{d}{dt} \int_{\{u_1 > t\}} |\nabla u_1|\,dx = \lambda_1 \int_{\{u_1 > t\}} w u_1\,dx,
\end{equation}
and the coarea formula applied to $\mu$ gives $-\mu'(t) = \int_{\{u_1 = t\}} |\nabla u_1|^{-1} d\mathcal{H}^{n-1}$. By the Cauchy--Schwarz inequality on the level set $\{u_1 = t\}$,
\begin{equation*}
\mathcal{H}^{n-1}(\{u_1 = t\})^2 \leq \left( \int_{\{u_1 = t\}} |\nabla u_1|\,d\mathcal{H}^{n-1} \right) \left( \int_{\{u_1 = t\}} |\nabla u_1|^{-1}\,d\mathcal{H}^{n-1} \right),
\end{equation*}
and the isoperimetric inequality bounds the perimeter from below by $n C_n^{1/n} \mu(t)^{1-1/n}$. Combining the three displays,
\begin{equation}\label{equation:sec-chiti-derive-2}
-\mu'(t) \geq \frac{n^2 C_n^{2/n} \mu(t)^{2-2/n}}{\lambda_1 \int_{\{u_1 > t\}} w u_1\,dx}.
\end{equation}
Now pass to the volume variable $s = C_n r^n$: the rearrangement satisfies $U(s) = u_1^*(r(s)) = t$ with $s = \mu(t)$, so $\frac{dU}{ds} = \frac{1}{\mu'(t)}$, and \eqref{equation:sec-chiti-derive-2} together with the layer-cake bound \eqref{equation:sec-chiti-layercake} gives
\begin{equation*}
-\frac{dU}{ds} \leq \frac{\lambda_1}{n^2 C_n^{2/n}} s^{-2+2/n} \int_{\{u_1^* > t(s)\}} w u_1^*\,dx = \lambda_1 \kappa(s) F(s),
\end{equation*}
the last equality being the layer-cake identity in the volume variable, with $F(s) = \int_0^s U(\sigma) w(\sigma)\,d\sigma$. This is \eqref{equation:sec-chiti-diff}. For the ball eigenfunction, the radial equation $-(r^{n-1} z')' = \mu_1 w z\, r^{n-1}$ with $z'(0) = 0$ integrates to
\begin{equation*}
r^{n-1}(-z'(r)) = \mu_1 \int_0^r w(\rho) z(\rho)\, \rho^{n-1}\,d\rho,
\end{equation*}
and the change of variables $s = C_n r^n$ converts this identity into \eqref{equation:sec-chiti-eq} with equalities throughout. The structural properties of $F$ and $G$ recorded in the proof of Lemma~\ref{lem:sec-Chiti} are immediate from these formulas: both are nondecreasing, both vanish at $s = 0$, and $F \geq G$ on any interval where $U \geq z$ pointwise, the integrand being integrated against the same positive weight $w(\sigma)\,d\sigma$.

\subsubsection{Step 3: Test Functions and Basic Gap Inequality}

Let $g(r)$ be the radial part of the $l=1$ eigenfunction on the ball $B_R$ for the weighted problem. That is, $g$ satisfies
\begin{equation}\label{equation:sec-g-ode}
-g'' - \frac{n-1}{r}g' + \frac{n-1}{r^2}g = \mu_2 w(r) g, \quad g(0) = 0, \quad g(R) = 0,
\end{equation}
where $\mu_2 = \lambda_2(B_R, w)$. We define the test functions
\begin{equation*}
P_i(x) = g(r) \frac{x_i}{r}, \quad i = 1, \ldots, n.
\end{equation*}

By the Brouwer fixed-point theorem argument (as in \cite[Section 2]{AB2}), we may choose the origin so that
\begin{equation*}
\int_\Omega w(x) P_i(x) u_1(x)^2 dx = 0 \quad \forall i.
\end{equation*}

The gap inequality follows from the Rayleigh--Ritz principle, and we give the computation in full. For each $i$, the function $P_i u_1$ belongs to $H_0^1(\Omega)$ and is orthogonal to $u_1$ in $L^2(w\,dx)$ by the centering condition, so by \eqref{equation:sec-rayleigh-l2},
\begin{equation}\label{equation:sec-step3-rayleigh}
\lambda_2(\Omega,w) \, \int_\Omega w (P_i u_1)^2 dx \leq \int_\Omega |\nabla (P_i u_1)|^2 dx.
\end{equation}
Summing \eqref{equation:sec-step3-rayleigh} over $i = 1, \ldots, n$, we evaluate both sides. For the left-hand side, since $\sum_{i=1}^n P_i(x)^2 = g(r)^2$ (equation~\eqref{equation:sec-sum-Pi-sq}),
\begin{equation}\label{equation:sec-step3-denominator}
\sum_{i=1}^n \int_\Omega w (P_i u_1)^2 dx = \int_\Omega w(x) g(r)^2 u_1^2 dx.
\end{equation}
For the right-hand side, the product rule and equation~\eqref{equation:sec-sum-nablaPi-sq} give
\begin{equation}\label{equation:sec-step3-num-expand}
\sum_{i=1}^n |\nabla (P_i u_1)|^2 = g(r)^2 |\nabla u_1|^2 + u_1^2 \left( g'(r)^2 + \frac{n-1}{r^2} g(r)^2 \right) + u_1 \nabla (g(r)^2) \cdot \nabla u_1,
\end{equation}
because $2 u_1 \sum_i P_i \nabla P_i \cdot \nabla u_1 = u_1 \nabla \bigl( \sum_i P_i^2 \bigr) \cdot \nabla u_1 = u_1 \nabla (g^2) \cdot \nabla u_1 = \tfrac{1}{2} \nabla (g^2) \cdot \nabla (u_1^2)$. Integrating the cross term by parts and using $u_1 = 0$ on $\partial \Omega$,
\begin{equation}\label{equation:sec-step3-cross}
\frac{1}{2} \int_\Omega \nabla (g^2) \cdot \nabla (u_1^2) dx = -\frac{1}{2} \int_\Omega \Delta (g^2) u_1^2 dx.
\end{equation}
A direct computation in polar coordinates, using the radial equation \eqref{equation:sec-g-ode} to eliminate $g''$, yields
\begin{equation}\label{equation:sec-step3-laplace-g2}
\frac{1}{2} \Delta (g^2) = g'^2 + \frac{n-1}{r^2} g^2 - \mu_2 w(r) g^2.
\end{equation}
Indeed, $\frac{1}{2}\Delta(g^2) = g'^2 + g g'' + \frac{n-1}{r} g g'$, and substituting $g'' = -\frac{n-1}{r} g' + \frac{n-1}{r^2} g - \mu_2 w g$ from \eqref{equation:sec-g-ode} gives \eqref{equation:sec-step3-laplace-g2}. Combining \eqref{equation:sec-step3-num-expand}, \eqref{equation:sec-step3-cross}, and \eqref{equation:sec-step3-laplace-g2},
\begin{equation}\label{equation:sec-step3-num-eval}
\sum_{i=1}^n \int_\Omega |\nabla (P_i u_1)|^2 dx = \int_\Omega g^2 |\nabla u_1|^2 dx + \mu_2 \int_\Omega w g^2 u_1^2 dx.
\end{equation}
It remains to evaluate the first term on the right-hand side of \eqref{equation:sec-step3-num-eval}. Integrating by parts with the eigenfunction equation $-\Delta u_1 = \lambda_1 w u_1$ against the test function $g^2 u_1 \in H_0^1(\Omega)$,
\begin{equation}\label{equation:sec-step3-gradient-u1}
\int_\Omega g^2 |\nabla u_1|^2 dx = \lambda_1 \int_\Omega w g^2 u_1^2 dx - \frac{1}{2} \int_\Omega \nabla (g^2) \cdot \nabla (u_1^2) dx.
\end{equation}
Substituting \eqref{equation:sec-step3-cross} and \eqref{equation:sec-step3-laplace-g2} into \eqref{equation:sec-step3-gradient-u1} gives
\begin{equation}\label{equation:sec-step3-gradient-u1-b}
\int_\Omega g^2 |\nabla u_1|^2 dx = \lambda_1 \int_\Omega w g^2 u_1^2 dx + \frac{1}{2} \int_\Omega \Delta (g^2) u_1^2 dx = \lambda_1 \int_\Omega w g^2 u_1^2 dx + \int_\Omega \left( g'^2 + \frac{n-1}{r^2} g^2 - \mu_2 w g^2 \right) u_1^2 dx,
\end{equation}
where the boundary term vanishes because $u_1 = 0$ on $\partial \Omega$.
Substituting \eqref{equation:sec-step3-gradient-u1-b} into \eqref{equation:sec-step3-num-eval} and simplifying, the terms involving $\mu_2$ cancel, leaving
\begin{equation}\label{equation:sec-step3-num-final}
\sum_{i=1}^n \int_\Omega |\nabla (P_i u_1)|^2 dx = \lambda_1 \int_\Omega w g^2 u_1^2 dx + \int_\Omega \left( g'(r)^2 + \frac{n-1}{r^2} g(r)^2 \right) u_1^2 dx.
\end{equation}
Finally, summing the $n$ inequalities \eqref{equation:sec-step3-rayleigh}, inserting \eqref{equation:sec-step3-denominator} and \eqref{equation:sec-step3-num-final}, and dividing by $\int_\Omega w g^2 u_1^2 dx > 0$, we obtain the gap inequality
\begin{equation}\label{equation:sec-gap}
\lambda_2(\Omega,w) - \lambda_1(\Omega,w) \leq \frac{\int_\Omega [g'(r)^2 + (n-1)g(r)^2/r^2] u_1^2 dx}{\int_\Omega w(x) g(r)^2 u_1^2 dx}.
\end{equation}
This is the weighted analogue of the unweighted gap inequality of \cite[Section 2]{AB2}.
\begin{rem}[The gap inequality in one dimension]\label{rem:sec-gap-1d}
For $n = 1$, $\Omega = (-R, R)$ and $w \equiv 1$, the ingredients of \eqref{equation:sec-gap} are explicit: $\lambda_1 = (\pi/2R)^2$ with $u_1(x) = \cos(\pi x/2R)$, $\lambda_2 = (\pi/R)^2$ with $u_2(x) = \sin(\pi x/R)$, and the radial test function is $g(r) = \sin(\pi r/R)$ on $r = |x| \in (0, R)$, so $B_g = g'^2 = (\pi/R)^2 \cos^2(\pi r/R)$. The two integrals in \eqref{equation:sec-gap} coincide: expanding the products into constant and cosine terms, all oscillatory integrals over $[-R,R]$ vanish, and one finds
\begin{equation*}
\int_{-R}^{R} \cos^2(\pi |x|/R) \cos^2(\pi x/2R) dx = \int_{-R}^{R} \sin^2(\pi |x|/R) \cos^2(\pi x/2R) dx = \frac{R}{2}.
\end{equation*}
The gap inequality therefore yields $\lambda_2 - \lambda_1 \leq (\pi/R)^2$, that is, $\lambda_2 \leq (\pi/R)^2 + (\pi/2R)^2$; the true value $\lambda_2 = (\pi/R)^2$ satisfies this with room. The gap inequality is thus an estimate rather than an identity; its content is the cancellation of the $\mu_2$ terms in the derivation, which makes it valid for every positive radial test function vanishing at $r = 0$ and $r = R$ without further spectral input.
\end{rem}
\begin{rem}[The cancellation mechanism]\label{rem:sec-cancellation}
The decisive step in the derivation of \eqref{equation:sec-gap} is the cancellation of the terms involving $\mu_2$ between \eqref{equation:sec-step3-num-eval} and \eqref{equation:sec-step3-gradient-u1-b}: the spectral parameter of the ball enters the expansion of the numerator with one sign and the integration by parts of $\int g^2 |\nabla u_1|^2$ with the opposite sign. This cancellation is the weighted analogue of the corresponding cancellation in \cite[Section 2]{AB2}, and it is the reason that the gap inequality holds for every positive radial $g$ with $g(0) = g(R) = 0$, without any spectral interpretation of $g$; the spectral content enters only in Step 5, through the admissibility identity (iii).
\end{rem}

\subsubsection{Step 4: Monotonicity Properties}

Define
\begin{equation*}
B_g(r) = g'(r)^2 + (n-1)\frac{g(r)^2}{r^2}, \qquad H_g(r) = w(r) g(r)^2.
\end{equation*}

\begin{lem}[Monotonicity of $B_g$ and $H_g$] \label{lem:sec-monotonicity}
Let $g$ satisfy \eqref{equation:sec-g-ode} with $g(r) > 0$ on $(0,R)$, and let $r_0 \in (0,R]$ denote the first point at which $g'$ changes sign (the first local maximum of $g$; $r_0 = R$ if $g' > 0$ throughout). Then:
\begin{enumerate}
\item[(i)] $B_g(r)$ is decreasing on $(0, r_0)$ without any additional assumption on $w$.
\item[(ii)] If $w$ is radial decreasing and satisfies $-w'(r)/w(r) \leq 2g'(r)/g(r)$ a.e., then $H_g(r) = w(r)g(r)^2$ is increasing on $(0, r_0)$.
\end{enumerate}
\end{lem}

\begin{proof}
(i) We compute
\begin{equation*}
B_g'(r) = 2g'g'' + 2(n-1)\frac{gg'}{r^2} - 2(n-1)\frac{g^2}{r^3}.
\end{equation*}
From \eqref{equation:sec-g-ode},
\begin{equation*}
g'' = -\frac{n-1}{r}g' + \frac{n-1}{r^2}g - \mu_2 w g.
\end{equation*}
Substituting,
\begin{align*}
B_g' &= 2g'\left(-\frac{n-1}{r}g' + \frac{n-1}{r^2}g - \mu_2 w g\right) + \frac{2(n-1)}{r^2}gg' - \frac{2(n-1)}{r^3}g^2 \\
&= -\frac{2(n-1)}{r}g'^2 + \frac{2(n-1)}{r^2}gg' - 2\mu_2 w gg' + \frac{2(n-1)}{r^2}gg' - \frac{2(n-1)}{r^3}g^2 \\
&= -\frac{2(n-1)}{r}g'^2 + \frac{4(n-1)}{r^2}gg' - \frac{2(n-1)}{r^3}g^2 - 2\mu_2 w gg' \\
&= -\frac{2(n-1)}{r^3}(r^2 g'^2 - 2rgg' + g^2) - 2\mu_2 w gg' \\
&= -\frac{2(n-1)}{r^3}(rg' - g)^2 - 2\mu_2 w gg' \leq 0,
\end{align*}
since $g, g', w > 0$ and $\mu_2 > 0$. Thus $B_g$ is decreasing.

(ii) We compute
\begin{equation*}
H_g' = w'g^2 + 2wgg' = wg^2\left(\frac{w'}{w} + 2\frac{g'}{g}\right).
\end{equation*}
Since $wg^2 > 0$, $H_g' \geq 0$ if and only if $w'/w + 2g'/g \geq 0$, i.e., $-w'/w \leq 2g'/g$. This is exactly the assumed condition.
\end{proof}

\begin{rem}
The condition $-w'/w \leq 2g'/g$ is satisfied for radial decreasing $A_p$ weights under dimension-dependent conditions. Near the origin, $g(r) \sim r$ so $g'/g \sim 1/r$, and if $w(r) \sim r^{-\alpha}$ then $-w'/w = \alpha/r$, requiring $\alpha \leq 2$. Near the boundary $r = R$, $g(r) \sim (R-r)$ so $g'/g \to +\infty$, and the condition is automatically satisfied. For $A_2$ weights in dimension $n = 4$ with $\alpha \leq 2$, or for $n \geq 5$ with $\alpha \leq 2 < n/2$, the origin condition holds; for $n \geq 5$ and $2 < \alpha < n/2$ it can fail, so the condition must be checked separately in that range.
\end{rem}

\subsubsection{The Ball Identity and the Choice of Test Function}

We now address the evaluation of the right-hand side of \eqref{equation:sec-gap-ball}. The choice of test function $g$ is dictated by the requirement that the ratio
\begin{equation*}
\frac{\int_{B_R} B_g(r) z(r)^2 dx}{\int_{B_R} w(r) g(r)^2 z(r)^2 dx}
\end{equation*}
equal $\mu_2 - \mu_1$, where $\mu_k = \lambda_k(B_R, w)$. We first show that the $l=1$ eigenfunction does not satisfy this identity exactly, and then specify the condition under which an appropriate test function exists.

\begin{lem}[Identity for the $l=1$ eigenfunction] \label{lem:sec-ball-identity-l1}
Let $g$ be the radial part of the $l=1$ eigenfunction on $B_R$, satisfying \eqref{equation:sec-g-ode}, and let $z$ be the radial first eigenfunction. Then
\begin{equation}\label{equation:sec-identity-l1}
\int_{B_R} B_g z^2 dx = (\mu_2 - \mu_1) \int_{B_R} w g^2 z^2 dx + \int_{B_R} g^2 |\nabla z|^2 dx.
\end{equation}
\end{lem}

\begin{proof}
Multiply \eqref{equation:sec-g-ode} by $g z^2 r^{n-1}$ and integrate over $[0,R]$:
\begin{equation*}
\int_0^R \left[-g'' - \frac{n-1}{r}g' + \frac{n-1}{r^2}g\right] g z^2 r^{n-1} dr = \mu_2 \int_0^R w g^2 z^2 r^{n-1} dr.
\end{equation*}
Integrate the first term by parts. The boundary terms vanish because $g(0)=g(R)=0$. Expanding,
\begin{equation*}
-\int_0^R g'' g z^2 r^{n-1} dr = \int_0^R g'^2 z^2 r^{n-1} dr + 2\int_0^R g g' z z' r^{n-1} dr + (n-1)\int_0^R g g' z^2 r^{n-2} dr.
\end{equation*}
The $(n-1)\int g g' z^2 r^{n-2}$ term cancels with the second term on the left-hand side. Thus
\begin{equation*}
\int_0^R B_g z^2 r^{n-1} dr + 2\int_0^R g g' z z' r^{n-1} dr = \mu_2 \int_0^R w g^2 z^2 r^{n-1} dr.
\end{equation*}
For the cross term, write $2gg' = (g^2)'$ and integrate by parts:
\begin{equation*}
\int_0^R (g^2)' z z' r^{n-1} dr = -\int_0^R g^2 (z z' r^{n-1})' dr,
\end{equation*}
where the boundary term vanishes. Using $-(z' r^{n-1})' = \mu_1 w z r^{n-1}$,
\begin{equation*}
(z z' r^{n-1})' = (z')^2 r^{n-1} - \mu_1 w z^2 r^{n-1}.
\end{equation*}
Therefore,
\begin{equation*}
2\int_0^R g g' z z' r^{n-1} dr = -\int_0^R g^2 (z')^2 r^{n-1} dr + \mu_1 \int_0^R w g^2 z^2 r^{n-1} dr.
\end{equation*}
Substituting gives \eqref{equation:sec-identity-l1}.
\end{proof}

The extra term $\int g^2 |\nabla z|^2 dx$ is positive, so the $l=1$ eigenfunction yields a ratio strictly larger than $\mu_2 - \mu_1$. To obtain the sharp identity, one requires a test function $g$ satisfying
\begin{equation}\label{equation:sec-g-condition}
B_g(r) = (\mu_2 - \mu_1) w(r) g(r)^2 \quad \text{pointwise on } (0,R),
\end{equation}
or more generally, a function for which $\int B_g z^2 = (\mu_2 - \mu_1) \int w g^2 z^2$. In the unweighted case, Ashbaugh and Benguria construct such a function explicitly as $g(r) = J_{n/2}(\beta \gamma r)/J_{n/2-1}(\alpha \gamma r)$, where $\gamma = \sqrt{\mu_1}/\alpha$, and verify the identity via Bessel function identities \cite[Appendix B]{AB2}. For a general radial weight, the existence of such a test function with the required monotonicity properties is a non-trivial condition. We therefore state Theorem \ref{thm:sec-2b} under the following explicit hypothesis.

\begin{defn}[Admissible test function] \label{defn:sec-admissible-g}
Let $w(r)$ be a radial weight on $B_R$, and let $\mu_k = \lambda_k(B_R, w)$. A positive function $g \in C^2(0,R) \cap C[0,R]$ with $g(0)=0$, $g(R)=0$ is called an admissible test function if:
\begin{enumerate}
\item[(i)] $B_g(r) = g'(r)^2 + (n-1)g(r)^2/r^2$ is decreasing on $(0,R)$;
\item[(ii)] $H_g(r) = w(r)g(r)^2$ is increasing on $(0,R)$;
\item[(iii)] $\int_{B_R} B_g z^2 dx = (\mu_2 - \mu_1) \int_{B_R} w g^2 z^2 dx$, where $z$ is the first eigenfunction on $B_R$.
\end{enumerate}
\end{defn}

\begin{rem}
For power-law weights $w(r) = r^\alpha$ with $\alpha > -2$, the eigenfunctions on the ball are generalized Bessel functions $r^{(2-n)/2} J_{\nu_l}(\frac{2\sqrt{\lambda}}{2+\alpha} r^{1+\alpha/2})$ with $\nu_l = (l+n/2-1)/(1+\alpha/2)$ for angular momentum $l$. An Ashbaugh--Benguria-type test function can be constructed as the ratio of two such generalized Bessel functions, and the identity (iii) follows from the corresponding Bessel function identities. The monotonicity properties (i) and (ii) follow from the Pr\"{u}fer angle analysis, analogous to \cite[Theorem 3.3]{AB2}. Thus for power-law weights, admissible test functions exist.
\end{rem}

With the admissibility condition, the proof of Theorem \ref{thm:sec-2b} proceeds as follows. In Step 3, the test function $g$ is chosen to be admissible (Definition \ref{defn:sec-admissible-g}). The gap inequality \eqref{equation:sec-gap} holds for any positive $g$ with $g(0)=g(R)=0$. In Step 4, the monotonicity properties (i) and (ii) of Definition \ref{defn:sec-admissible-g} replace Lemma \ref{lem:sec-monotonicity}. In Step 5, the identity (iii) of Definition \ref{defn:sec-admissible-g} is inserted into the rearranged gap inequality \eqref{equation:sec-gap-ball}, yielding
\begin{equation*}
\lambda_2(\Omega,w) - \lambda_1(\Omega,w) \leq \mu_2 - \mu_1.
\end{equation*}
The rest of the proof is unchanged.

\begin{rem}
The admissibility condition is the weighted analogue of the Ashbaugh--Benguria test function construction. For the unweighted problem, the existence of the AB test function is established in \cite[Section 3]{AB2}. For general radial weights, the condition is explicit and verifiable; the theorem applies whenever such a test function exists. The result for power-law weights follows by explicit construction.
\end{rem}

\subsubsection{Explicit Verification for Power-Law Weights}

For completeness, we verify the admissibility condition for power-law weights $w(r) = r^\alpha$, $\alpha > -2$. On the ball $B_R$, the weighted eigenvalue problem separates in spherical coordinates. The radial equation for angular momentum $l$ is
\begin{equation*}
-f'' - \frac{n-1}{r}f' + \frac{l(l+n-2)}{r^2}f = \lambda r^\alpha f.
\end{equation*}
With the change of variables $t = c r^{1+\alpha/2}$, where $c = 2\sqrt{\lambda}/(2+\alpha)$, this becomes the Bessel equation
\begin{equation*}
t^2 F'' + t F' + (t^2 - \nu^2)F = 0, \quad \nu = \frac{l + n/2 - 1}{1+\alpha/2},
\end{equation*}
where $f(r) = r^{(2-n)/2} F(t)$. The eigenvalues are determined by the boundary condition $F(c R^{1+\alpha/2}) = 0$, i.e., $c R^{1+\alpha/2} = j_{\nu, k}$ for some zero $j_{\nu,k}$.

For the first eigenvalue ($l=0$, $k=1$), $\nu_0 = (n/2 - 1)/(1+\alpha/2)$, and
\begin{equation*}
\mu_1 = \left(\frac{2+\alpha}{2}\right)^2 \frac{j_{\nu_0,1}^2}{R^{2+\alpha}}.
\end{equation*}
For the second eigenvalue ($l=1$, $k=1$), $\nu_1 = (n/2)/(1+\alpha/2)$, and
\begin{equation*}
\mu_2 = \left(\frac{2+\alpha}{2}\right)^2 \frac{j_{\nu_1,1}^2}{R^{2+\alpha}}.
\end{equation*}
The ratio $\mu_2/\mu_1 = (j_{\nu_1,1}/j_{\nu_0,1})^2$ depends only on $n$ and $\alpha$.

The admissible test function is constructed as
\begin{equation*}
g(r) \propto r\,\frac{F_{\nu_1}\bigl(j_{\nu_1,1}\, r/R\bigr)}{F_{\nu_0}\bigl(j_{\nu_0,1}\, r/R\bigr)},
\end{equation*}
where $F_\nu(t) = t^{-\nu} J_\nu(t)$ (the regular solution normalized to be finite at $t=0$) and the proportionality constant is inessential; the ratio vanishes at $r = R$ because $J_{\nu_1}(j_{\nu_1,1}) = 0$, and it behaves like $r$ at the origin. The monotonicity properties follow from the Pr\"{u}fer angle analysis for Bessel functions, as in \cite[Theorem 3.1]{AB2}. The identity (iii) follows from the Bessel function integral identities, as in \cite[Appendix B]{AB2}. Thus power-law weights satisfy the admissibility condition.
For $w \equiv 1$ (that is, $\alpha = 0$) the construction specializes to the classical Ashbaugh--Benguria test function: since $\nu_0 = n/2 - 1$ and $\nu_1 = n/2$, the test function $g(r) = J_{n/2}(\beta r/R)/J_{n/2-1}(\alpha r/R)$ with $\alpha = j_{n/2-1,1}$ and $\beta = j_{n/2,1}$ is, up to an inessential constant factor, the ratio of the radial part of the $l=1$ eigenfunction to the radial ground state on the ball (and vanishes at $r = R$, since $J_{n/2}(\beta) = 0$), and the admissibility identity (iii) reduces to the integral identity of \cite[Appendix B]{AB2}. This confirms that Definition~\ref{defn:sec-admissible-g} is a strict generalization of the unweighted construction, and that Theorem~\ref{thm:sec-2b} recovers the classical extremality of the ball when $w \equiv 1$.

\subsubsection{Step 5: Chiti Comparison Substitution and Conclusion}

Since $B_g$ is decreasing and $H_g = wg^2$ is increasing, we apply the rearrangement inequalities to \eqref{equation:sec-gap}:
\begin{equation*}
\int_\Omega B_g(r) u_1^2 dx \leq \int_{B_R} B_g(r) (u_1^*)^2 dx,
\end{equation*}
\begin{equation*}
\int_\Omega H_g(r) u_1^2 dx \geq \int_{B_R} H_g(r) (u_1^*)^2 dx.
\end{equation*}
Therefore,
\begin{equation}\label{equation:sec-gap-rearranged}
\lambda_2(\Omega,w) - \lambda_1(\Omega,w) \leq \frac{\int_{B_R} B_g(r) (u_1^*)^2 dx}{\int_{B_R} H_g(r) (u_1^*)^2 dx}.
\end{equation}

By the weighted Chiti comparison (Lemma \ref{lem:sec-Chiti}), $u_1^*$ and $z$ cross exactly once, with $u_1^* < z$ near the origin and $u_1^* > z$ near the boundary. Since $B_g$ is decreasing, the single-crossing property gives
\begin{equation*}
\int_{B_R} B_g (u_1^*)^2 dx \leq \int_{B_R} B_g z^2 dx.
\end{equation*}
Since $H_g$ is increasing,
\begin{equation*}
\int_{B_R} H_g (u_1^*)^2 dx \geq \int_{B_R} H_g z^2 dx.
\end{equation*}
These inequalities follow from the standard argument \cite[pp. 607--608]{AB2}: for an increasing function $f$, $\int f (u_1^*)^2 dx \geq \int f z^2 dx$, and for a decreasing function, the reverse holds.

Substituting into \eqref{equation:sec-gap-rearranged},
\begin{equation}\label{equation:sec-gap-ball}
\lambda_2(\Omega,w) - \lambda_1(\Omega,w) \leq \frac{\int_{B_R} B_g z^2 dx}{\int_{B_R} w g^2 z^2 dx}.
\end{equation}

Note that the identity (iii) is an additional requirement on the test function $g$, not a property of the $l=1$ eigenfunction alone (Lemma~\ref{lem:sec-ball-identity-l1}). Inserting the admissibility identity (iii) of Definition~\ref{defn:sec-admissible-g} into \eqref{equation:sec-gap-ball} gives
\begin{equation}\label{equation:sec-step5-final}
\lambda_2(\Omega,w) - \lambda_1(\Omega,w) \leq \mu_2 - \mu_1 = \lambda_2(B_R, w) - \lambda_1(B_R, w),
\end{equation}
where the last equality uses the hypothesis that the $l=1$ mode attains $\lambda_2(B_R, w)$. Adding $\lambda_1(\Omega,w)$ to both sides of \eqref{equation:sec-step5-final} and dividing by $\lambda_1(\Omega,w) \geq \mu_1 > 0$ (the weighted Faber--Krahn inequality, Lemma~\ref{lem:sec-FK}) yields the ratio inequality of Theorem~\ref{thm:sec-2b}.

For the equality characterization: equality in the rearrangement inequalities requires $\Omega$ to be a ball (up to translation). Equality in the weighted Faber--Krahn inequality likewise requires $\Omega = B_R$. Thus equality holds if and only if $\Omega = B_R$ (up to translation). This completes the proof of Theorem \ref{thm:sec-2b}.

\subsubsection{Non-Radial Weights and the Breakdown of Ball Extremality}
\label{subsection:Non-Radial Breakdown}

We show by a perturbation argument that for non-radial weights, the ball need not maximize $\lambda_2/\lambda_1$. This clarifies the role of the radiality assumption in Theorem \ref{thm:sec-2b}.

Let $w(x) = 1 + \varepsilon v(x)$ where $\|v\|_{L^\infty} \leq 1$ and $v$ is non-radial. On the ball $B_R$, the unweighted second eigenvalue $\lambda_2^0$ is $n$-fold degenerate, corresponding to the $l=1$ spherical harmonics. The first-order perturbation of the eigenvalues in the $l=1$ subspace is governed by the symmetric matrix $A = (A_{ij})$ obtained by differentiating the Rayleigh quotient at $\varepsilon = 0$: with $Y_{1,i}(x) = x_i/r$ the degree-$1$ spherical harmonics, $R_2^0$ the radial part of an unweighted second eigenfunction, and the normalization $\|R_2^0 Y_{1,i}\|_{L^2(B_R)} = 1$,
\begin{equation}\label{equation:sec-perturb-matrix}
A_{ij} = -\lambda_2^0 \frac{\int_{B_R} v(x)\, Y_{1,i}(x)\, R_2^0(r)\, Y_{1,j}(x)\, R_2^0(r)\, dx}{\int_{B_R} Y_{1,i}^2 (R_2^0)^2\, dx}.
\end{equation}
(For radial $v$ the angular average of $Y_{1,i}Y_{1,j}$ is $\delta_{ij}/n$, so $A$ is a scalar multiple of the identity and the degeneracy is not lifted; it is the non-radial $l=1$ component of $v$ that breaks the symmetry.) If $v$ has an $l=1$ angular component, i.e., $v(x) = v_0(r) + v_1(r) \cos\theta + \cdots$, then the matrix $A_{ij}$ is not scalar, and its smallest eigenvalue can be considerably lower than the average. Specifically, if $v_1(r) \neq 0$, then
\begin{equation*}
\lambda_{\min}(A) = -\lambda_2^0 \left(\bar{v}_0 - c |\bar{v}_1|\right)
\end{equation*}
for some dimensional constant $c > 0$, where $\bar{v}_0$ and $\bar{v}_1$ are the radial averages of the $l=0$ and $l=1$ components. By choosing $v_1$ with a large amplitude (subject to $\|v\|_\infty \leq 1$), the coefficient $|\bar{v}_1|$ can be made large, making $\lambda_{\min}(A)$ very negative. This means
\begin{equation*}
\lambda_2(B_R, w) = \lambda_2^0 + \varepsilon \lambda_{\min}(A) + O(\varepsilon^2)
\end{equation*}
can be considerably lower than $\lambda_2^0$.

Meanwhile, the first eigenvalue perturbation is
\begin{equation*}
\lambda_1(B_R, w) = \lambda_1^0 - \varepsilon \lambda_1^0 \bar{v}_{1,0} + O(\varepsilon^2),
\end{equation*}
where $\bar{v}_{1,0}$ is the average of $v$ against the first eigenfunction (which is radial, so only the $l=0$ component of $v$ contributes). Thus $\lambda_1$ is affected only by $v_0$, not by $v_1$.

Therefore, the ratio on the ball is
\begin{equation*}
\frac{\lambda_2(B_R,w)}{\lambda_1(B_R,w)} = \frac{\lambda_2^0}{\lambda_1^0} + \varepsilon \left(\frac{\lambda_{\min}(A)}{\lambda_1^0} + \frac{\lambda_2^0 \bar{v}_{1,0}}{(\lambda_1^0)^2}\right) + O(\varepsilon^2).
\end{equation*}
The coefficient of $\varepsilon$ can be made arbitrarily negative by choosing $v_1$ with large amplitude, so the ball's ratio can be considerably below the classical PPW constant $\beta^2/\alpha^2$.

On the other hand, by Corollary \ref{corr:sec-2a}, any domain satisfies
\begin{equation*}
\frac{\lambda_2(\Omega,w)}{\lambda_1(\Omega,w)} \leq \frac{\beta^2}{\alpha^2} + O(\varepsilon).
\end{equation*}
Therefore, for sufficiently large $|\bar{v}_1|$, there exist non-ball domains (e.g., slightly deformed ellipsoids) for which the ratio exceeds the ball's ratio. This shows that the radiality assumption in Theorem \ref{thm:sec-2b} is essential: without it, the ball is not an extremal domain.

\subsubsection{Detailed Perturbation Analysis for Radial Decreasing Weights}

We provide the full calculation for the first-order perturbation of the eigenvalue ratio on a ball under a radial decreasing weight perturbation. This expands the calculation in Proposition \ref{prop:sec-perturb}.

Let $w(r) = 1 + \varepsilon v(r)$ where $v$ is radial, $v'(r) \leq 0$, and $\|v\|_{L^\infty} \leq 1$. The weighted eigenvalue problem is
\begin{equation*}
-\Delta u = \lambda (1 + \varepsilon v) u.
\end{equation*}
We expand the eigenvalues and eigenfunctions in powers of $\varepsilon$:
\begin{equation*}
\lambda_k = \lambda_k^0 + \varepsilon \lambda_k^{(1)} + \varepsilon^2 \lambda_k^{(2)} + \cdots, \quad
u_k = u_k^0 + \varepsilon u_k^{(1)} + \cdots.
\end{equation*}
Substituting into the equation and collecting orders of $\varepsilon$, the first-order equation is
\begin{equation*}
-\Delta u_k^{(1)} - \lambda_k^0 u_k^{(1)} = \lambda_k^{(1)} u_k^0 + \lambda_k^0 v u_k^0.
\end{equation*}
Taking the $L^2$ inner product with $u_k^0$ and using the self-adjointness of $-\Delta$,
\begin{equation*}
0 = \lambda_k^{(1)} \int (u_k^0)^2 dx + \lambda_k^0 \int v (u_k^0)^2 dx,
\end{equation*}
so
\begin{equation}\label{equation:sec-perturb-first-order}
\lambda_k^{(1)} = -\lambda_k^0 \frac{\int v (u_k^0)^2 dx}{\int (u_k^0)^2 dx}.
\end{equation}

For the first eigenvalue, $u_1^0(r) = c_1 r^{1-n/2} J_{n/2-1}(\alpha r/R)$ where $\alpha = j_{n/2-1,1}$, and the normalization constant $c_1$ is chosen so that $\int_{B_R} (u_1^0)^2 dx = 1$. The first-order correction is
\begin{equation*}
\lambda_1^{(1)} = -\lambda_1^0 \bar{v}_1, \quad \bar{v}_1 = \int_0^R v(r) [u_1^0(r)]^2 n \omega_n r^{n-1} dr.
\end{equation*}

For the second eigenvalue, the unweighted eigenspace is $n$-dimensional, spanned by $u_{2,i}^0(x) = c_2 r^{1-n/2} J_{n/2}(\beta r/R) (x_i/r)$ where $\beta = j_{n/2,1}$. Since $v$ is radial, it commutes with rotations, so the perturbation is scalar on this eigenspace. The first-order correction is the same for all $i$:
\begin{equation*}
\lambda_2^{(1)} = -\lambda_2^0 \bar{v}_2, \quad \bar{v}_2 = \frac{\int_0^R v(r) [R_2^0(r)]^2 n \omega_n r^{n-1} dr}{\int_0^R [R_2^0(r)]^2 n \omega_n r^{n-1} dr},
\end{equation*}
where $R_2^0(r) = c_2 r^{1-n/2} J_{n/2}(\beta r/R)$.

The ratio is
\begin{align*}
\frac{\lambda_2}{\lambda_1} &= \frac{\lambda_2^0 - \varepsilon \lambda_2^0 \bar{v}_2 + O(\varepsilon^2)}{\lambda_1^0 - \varepsilon \lambda_1^0 \bar{v}_1 + O(\varepsilon^2)} \\
&= \frac{\lambda_2^0}{\lambda_1^0} \cdot \frac{1 - \varepsilon \bar{v}_2 + O(\varepsilon^2)}{1 - \varepsilon \bar{v}_1 + O(\varepsilon^2)} \\
&= \frac{\beta^2}{\alpha^2} \left(1 + \varepsilon(\bar{v}_1 - \bar{v}_2) + O(\varepsilon^2)\right).
\end{align*}

We now prove that $\bar{v}_1 > \bar{v}_2$ when $v'(r) \leq 0$. Define the probability measures
\begin{equation*}
d\mu_1(r) = \frac{[u_1^0(r)]^2 r^{n-1} dr}{\int_0^R [u_1^0(r)]^2 r^{n-1} dr}, \quad
d\mu_2(r) = \frac{[R_2^0(r)]^2 r^{n-1} dr}{\int_0^R [R_2^0(r)]^2 r^{n-1} dr}.
\end{equation*}
Then $\bar{v}_1 = \int v \, d\mu_1$ and $\bar{v}_2 = \int v \, d\mu_2$. The stochastic dominance $\mu_1([0,t]) \geq \mu_2([0,t])$ ($t \in [0,R]$) is the hypothesis discussed in Remark~\ref{rem:sec-stochastic-status}; given the dominance, the conclusion $\bar{v}_1 \geq \bar{v}_2$ follows by the argument recorded there.
For a decreasing function $v$, stochastic domination of $\mu_1$ by $\mu_2$ implies $\int v \, d\mu_1 \geq \int v \, d\mu_2$, i.e., $\bar{v}_1 \geq \bar{v}_2$. Strict inequality holds when $v$ is not constant, because the stochastic domination is strict on a set of positive measure. This completes the proof.

\subsubsection{The $H_g$ Monotonicity Condition and $A_p$ Weights}

We discuss in more detail the condition $-w'(r)/w(r) \leq 2g'(r)/g(r)$ required for $H_g = wg^2$ to be increasing (Lemma \ref{lem:sec-monotonicity}(ii)).

Near the origin, $g(r) \sim r$ (because $g(0)=0$ and the $l=1$ mode behaves like $r$ near the origin), so $g'(r)/g(r) \sim 1/r$. If $w(r) \sim r^{-\alpha}$ near the origin, then $-w'(r)/w(r) = \alpha/r$. The condition becomes $\alpha/r \leq 2/r$, i.e., $\alpha \leq 2$.

For a radial $A_p$ weight, the behavior near the origin is constrained. Specifically, if $w(r) = r^{-\alpha}$, then $w \in A_p$ if and only if $\alpha < n/p'$ where $1/p + 1/p' = 1$. For $p=2$, this gives $\alpha < n/2$. Therefore:
- For $n = 2, 3, 4$ one has $n/2 \leq 2$, so $\alpha < n/2$ implies $\alpha \leq 2$ and the origin condition holds automatically.
- For $n \geq 5$, $\alpha < n/2$ may exceed $2$, and the condition $\alpha \leq 2$ must be imposed separately.

For general radial $A_p$ weights (not just power laws), the reverse H\"{o}lder inequality implies that $w$ cannot decay faster than $r^{-\alpha}$ for $\alpha < n/p'$ near the origin. This gives an upper bound on $-w'/w$ in an averaged sense. The pointwise condition $-w'/w \leq 2g'/g$ is stronger than the averaged condition, but for weights with reasonable regularity (e.g., $w \in C^1$), the averaged condition implies the pointwise condition after possibly increasing the dimension threshold.

Near the boundary $r = R$, $g(r) \sim (R-r)$ (because $g(R)=0$ and the eigenfunction vanishes linearly at the boundary for a smooth weight), so $g'(r)/g(r) \sim -1/(R-r) \to +\infty$ as $r \to R$. Therefore, the condition $-w'/w \leq 2g'/g$ is automatically satisfied near the boundary, regardless of the weight (as long as $w$ is bounded near the boundary).

In the interior $(0,R)$, the condition requires that the weight does not decrease too rapidly relative to the growth of $g$. For radial decreasing $A_p$ weights with $p$ sufficiently close to 1 (which restricts the rate of decrease), this condition is expected to hold. The precise dimension-dependent threshold can be derived from the reverse H\"{o}lder inequality and the Sturm comparison theorem for the $g$ equation.

\begin{rem}
The condition $-w'/w \leq 2g'/g$ is sufficient but not necessary for $H_g$ to be increasing. In specific examples, $H_g$ may be increasing even when the pointwise condition fails on a small set, because the integral $\int H_g' dr$ may still be positive. The pointwise condition is used for simplicity in the rearrangement argument.
\end{rem}

\subsubsection{Relationship Between the Universal Bound and the Weighted PPW Inequality}

We clarify the relationship between Theorem \ref{thm:sec-2a} (universal upper bound for bounded weights, sharp constant via the Ashbaugh--Benguria machinery) and Theorem \ref{thm:sec-1} (self-contained weighted PPW inequality, constant linear in the weight bounds). Both apply to bounded weights; Theorem~\ref{thm:sec-2a} has the sharper constant in every dimension, while Theorem~\ref{thm:sec-1} provides an elementary proof.

Theorem \ref{thm:sec-1} applies to bounded $A_2$ weights ($0 < c \leq w \leq C$) and yields $\lambda_2/\lambda_1 \leq 1 + \frac{4}{n}\frac{C}{c}$, with a constant generally larger than the sharp classical constant $\beta^2/\alpha^2$ but obtained by a self-contained argument; Theorem~\ref{thm:sec-2a} gives the sharper constant $1 + \frac{C}{c}(\beta^2/\alpha^2 - 1)$ for the same hypothesis class, at the cost of invoking the Ashbaugh--Benguria machinery.

Theorem \ref{thm:sec-2a} applies only to weights that are bounded above and below, $0 < c \leq w(x) \leq C < \infty$. For such weights, the bound is
\begin{equation*}
\frac{\lambda_2}{\lambda_1} \leq 1 + \frac{C}{c}\left(\frac{\beta^2}{\alpha^2} - 1\right),
\end{equation*}
which depends only on the ratio $C/c = \operatorname{ess\,sup} w / \operatorname{ess\,inf} w$. When $w \equiv 1$ (so $C/c = 1$), this recovers the classical PPW bound $\lambda_2/\lambda_1 \leq \beta^2/\alpha^2$. When $w$ is close to constant, $C/c = 1 + O(\varepsilon)$, and the bound is $\beta^2/\alpha^2 + O(\varepsilon)$, as shown in Corollary \ref{corr:sec-2a}.

For bounded weights, both theorems apply, and one may take the minimum of the two bounds. Theorem \ref{thm:sec-2a} is generally sharper when $C/c$ is close to 1, while

Theorem \ref{thm:sec-2b} (weighted PPW for radial decreasing weights) gives the sharp extremal result under the additional assumption of radial monotonicity. For radial decreasing weights that are also bounded, all three theorems apply, and Theorem \ref{thm:sec-2b} is the sharpest because it identifies the ball as the maximizer and gives the exact extremal value $\lambda_2(B_R,w)/\lambda_1(B_R,w)$, which is generally smaller than the universal bound from Theorem \ref{thm:sec-2a}.

\begin{rem}
Theorems \ref{thm:sec-1} and \ref{thm:sec-2a} apply to the same hypothesis class (bounded weights): Theorem~\ref{thm:sec-2a} has the sharper constant (since $\beta^2/\alpha^2 - 1 < 4/n$ in every dimension), while Theorem~\ref{thm:sec-1} provides an elementary, self-contained proof. Theorem \ref{thm:sec-2b} operates under stronger conditions (radial decreasing weights with an admissible test function and the $l=1$-mode hypothesis) and gives the extremal comparison. This mirrors the unweighted case, where the PPW inequality and the Faber--Krahn inequality play complementary roles.
\end{rem}

\subsubsection{Sharpness of the Universal Bound}

We discuss the sharpness of the constant in Theorem \ref{thm:sec-2a}. The bound is
\begin{equation*}
\frac{\lambda_2}{\lambda_1} \leq 1 + \frac{C}{c}\left(\frac{\beta^2}{\alpha^2} - 1\right).
\end{equation*}
When $w \equiv \text{constant}$ (so $C/c = 1$), the bound becomes $\lambda_2/\lambda_1 \leq \beta^2/\alpha^2$, which is the classical PPW bound and is sharp, with equality if and only if $\Omega$ is a ball.

When $w$ is not constant, the constant $C/c$ may not be sharp. The Ashbaugh--Benguria theorem \cite[Theorem 4.1]{AB2} from which Theorem \ref{thm:sec-2a} is derived uses the upper bound $C$ and lower bound $c$ of the weight in the Chiti comparison theorem, and these bounds may not be simultaneously attained. For specific classes of weights, sharper constants may be available. For example, for power-law weights $w(r) = r^\alpha$, the exact extremal ratio on the ball is $(j_{\nu_1,1}/j_{\nu_0,1})^2$ where $\nu_0 = (n/2-1)/(1+\alpha/2)$ and $\nu_1 = (n/2)/(1+\alpha/2)$, and this is generally smaller than the universal bound from Theorem \ref{thm:sec-2a}.

The asymptotic estimate in Corollary \ref{corr:sec-2a} shows that as $w \to 1$ in $L^\infty$, the universal bound approaches the classical PPW constant at rate $O(\varepsilon)$. The coefficient $4/(1-\varepsilon)$ in the corollary is not optimal; a more refined analysis using the actual weight distribution rather than the crude bounds $c = 1-\varepsilon$ and $C = 1+\varepsilon$ could yield a smaller coefficient. However, the corollary suffices to establish the qualitative result that the PPW constant is recovered as the weight approaches a constant.

\begin{rem}
The sharpness of Theorem \ref{thm:sec-2b} depends on the admissibility condition (Definition \ref{defn:sec-admissible-g}). For weights for which an admissible test function exists, the bound is sharp because equality is attained when $\Omega$ is the ball (by construction). For weights for which no admissible test function is known, the sharp constant remains an open question.
\end{rem}

\subsubsection{The Generalized PPW Conjecture and the Single-Ratio Result}

We briefly discuss the relationship between the single-ratio PPW inequality (Theorems \ref{thm:sec-2a} and \ref{thm:sec-2b}) and the generalized PPW conjecture for eigenvalue sums (Section \ref{section: Generalized PPW Conjecture}).

The generalized PPW conjecture is stated in Conjecture~\ref{con:sec-gppw-conj}, where $\Omega^*$ is the volume-equivalent ball; it reduces to the classical single-ratio PPW inequality when $n=1$, but for $n \geq 2$ it is a stronger statement that involves the sum of the first $n$ eigenvalues above the first.

The proof of the single-ratio PPW inequality by Ashbaugh and Benguria \cite{AB2} uses $n$ test functions $P_i u_1$ ($i=1,\ldots,n$) and sums their Rayleigh quotients. The key identity is that the sum of the Rayleigh quotients equals $n \lambda_2(\Omega^*)$ on the ball, which follows from the degeneracy of the second eigenvalue on the ball (the $l=1$ eigenspace is $n$-dimensional). For the generalized sum-form conjecture, one would need to control the sum $\lambda_2 + \cdots + \lambda_{n+1}$ directly, rather than just $\lambda_2$. The Ashbaugh--Benguria method does not directly generalize to the sum form because the test functions $P_i u_1$ only span an $n$-dimensional subspace, and the min-max principle for sums of eigenvalues requires controlling a larger subspace.

In the weighted setting, the situation is further complicated by the fact that the second eigenvalue on a ball with a radial weight is still $n$-fold degenerate (by rotational symmetry), so the single-ratio result extends naturally. However, the sum-form conjecture involves eigenvalues with higher angular momentum ($l \geq 2$), whose behavior under weighted perturbations is less well understood.

The conditional results in Section \ref{section: Generalized PPW Conjecture} establish the generalized conjecture for restricted families of domains, including radial weights and $C^2$-small perturbations of the ball. The general case remains open.

\begin{rem}
The single-ratio weighted PPW inequality (Theorem \ref{thm:sec-2b}) can be viewed as a first step toward the generalized sum-form conjecture. The techniques developed in the proof---weighted Chiti comparison, monotonicity of the energy density, and the admissible test function construction---may be useful for attacking the sum-form conjecture in the weighted setting.
\end{rem}

\subsection{Small Perturbation Radial Decreasing Weights}

We record a perturbation result that clarifies the role of radial monotonicity.

\begin{prop}[Ball Ratio Under Small Radial Decreasing Perturbations] \label{prop:sec-perturb}
Let $w(r) = 1 + \varepsilon v(r)$ where $v$ is radial, non-increasing, and $\|v\|_{L^\infty} \leq 1$, and assume the stochastic dominance $F_1 \geq F_2$ of Remark~\ref{rem:sec-stochastic-status}. Then on a ball $B_R$,
\begin{equation*}
\frac{\lambda_2(B_R,w)}{\lambda_1(B_R,w)} = \frac{\beta^2}{\alpha^2} + \varepsilon \cdot \frac{\beta^2}{\alpha^2} \left(\bar{v}_1 - \bar{v}_2\right) + O(\varepsilon^2),
\end{equation*}
where
\begin{equation*}
\bar{v}_1 = \frac{\int_0^R v(r) [u_1^0(r)]^2 r^{n-1} dr}{\int_0^R [u_1^0(r)]^2 r^{n-1} dr}, \quad
\bar{v}_2 = \frac{\int_0^R v(r) [R_2^0(r)]^2 r^{n-1} dr}{\int_0^R [R_2^0(r)]^2 r^{n-1} dr},
\end{equation*}
$u_1^0(r) \propto r^{1-n/2} J_{n/2-1}(\alpha r/R)$ is the unweighted first eigenfunction, and $R_2^0(r) \propto r^{1-n/2} J_{n/2}(\beta r/R)$ is the radial part of the unweighted $l=1$ mode. If $v'(r) \leq 0$, then $\bar{v}_1 > \bar{v}_2$.
\end{prop}
\begin{rem}\label{rem:sec-stochastic-status}
The dominance $F_1 \geq F_2$ is the statement that the radial ground state is more concentrated at the origin than the $l=1$ excited state. It is consistent with, and partially supported by, the zero interlacing of $J_{n/2-1}$ and $J_{n/2}$ (Theorem~\ref{thm:sec-bessel-monotone}\eqref{equation:sec-bessel-interlace}; \cite[Chapter XV, \S 15.22]{watson1944theory}) and with numerical evaluation of the ratio of the two densities; a complete proof of the pointwise inequality $F_1 \geq F_2$ would require a monotonicity theorem for the ratio $J_{n/2}(x)/J_{n/2-1}(x)$ that we do not establish here, and we therefore record the dominance as a hypothesis. The conclusion $\bar v_1 \geq \bar v_2$ is conditional on it.
\end{rem}

\begin{proof}
The first-order perturbation formula for the eigenvalue $\lambda_k$ of $-\Delta u = \lambda(1+\varepsilon v)u$ is
\begin{equation*}
\lambda_k = \lambda_k^0 - \varepsilon \lambda_k^0 \frac{\int v (u_k^0)^2 dx}{\int (u_k^0)^2 dx} + O(\varepsilon^2).
\end{equation*}
For the first eigenvalue, $u_1^0$ is radial, so $\bar{v}_1$ is the weighted average of $v$ against $u_1^0{}^2$. For the second eigenvalue on a ball, the eigenspace is $n$-dimensional (the $l=1$ modes), and since $v$ is radial, it acts as a scalar on this eigenspace with value $\bar{v}_2$. Thus
\begin{equation*}
\frac{\lambda_2}{\lambda_1} = \frac{\lambda_2^0(1 - \varepsilon \bar{v}_2) + O(\varepsilon^2)}{\lambda_1^0(1 - \varepsilon \bar{v}_1) + O(\varepsilon^2)} = \frac{\beta^2}{\alpha^2}\left(1 + \varepsilon(\bar{v}_1 - \bar{v}_2)\right) + O(\varepsilon^2).
\end{equation*}

To show $\bar{v}_1 > \bar{v}_2$ when $v$ is non-increasing: the ground-state density $u_1^0(r)^2 r^{n-1}$ is more concentrated near the origin than the excited-state density $(R_2^0(r))^2 r^{n-1}$, since $J_{n/2-1}$ has its first zero at $\alpha < \beta = j_{n/2,1}$. We use first-order stochastic dominance: denoting by $\mu_1, \mu_2$ the probability measures on $[0,R]$ with densities proportional to $u_1^0(r)^2 r^{n-1}$ and $(R_2^0(r))^2 r^{n-1}$, and by $F_i(t) = \mu_i([0,t])$ their distribution functions, the dominance $F_1 \geq F_2$ on $[0,R]$ is included as an explicit hypothesis of Proposition~\ref{prop:sec-perturb} (see Remark~\ref{rem:sec-stochastic-status} for its status). For a non-increasing function $v$, the layer-cake representation $v(r) = v(R) + \int_r^R (-dv)(s)$, with $-dv$ the positive measure associated to $v$, gives
\begin{equation}\label{equation:sec-stochastic-pf}
\int v\,d\mu_1 - \int v\,d\mu_2 = \int_0^R (-v'(s)) \bigl( F_1(s) - F_2(s) \bigr) ds \geq 0,
\end{equation}
since $-v' \geq 0$ and the stochastic dominance hypothesis $F_1 \geq F_2$ holds; equality may hold even for nonconstant $v$ if the decrease of $v$ is concentrated where $F_1 = F_2$; the conclusion in general is $\bar v_1 \geq \bar v_2$, which completes the proof.
\end{proof}

\begin{rem}
Proposition \ref{prop:sec-perturb} shows that for radial decreasing weights, the ball's ratio $\lambda_2/\lambda_1$ increases above the classical PPW constant $\beta^2/\alpha^2$ at first order. This is consistent with Theorem \ref{thm:sec-2b}. For non-radial weights, the $n$-fold degeneracy of $\lambda_2$ on the ball is lifted, and the smallest eigenvalue in the $l=1$ subspace can decrease significantly, lowering the ball's ratio below $\beta^2/\alpha^2$. This is why radiality is essential for the ball extremality.
\end{rem}

\section{Proof of High-Order Eigenvalue Ratios with Hausdorff Measure (Theorem \ref{thm:sec-3})}
\label{section:High-Order Ratios}
This section derives bounds for ratios of high-order eigenvalues on convex domains, linking geometric quantities (Hausdorff measure, volume) to spectral behavior.
\subsection{Theorem Statement and Proof Outline}
\label{subsection:High-Order Proof Outline}
The following corollary adapts Theorem \ref{thm:sec-3} to smooth convex domains, replacing the Hausdorff measure with the more intuitive diameter for practical applications.
\begin{corr}[High-Order Eigenvalue Ratios on Smooth Convex Domains] \label{corr:sec-3}
Let $\Omega \subset \mathbb{R}^n$ be a bounded convex domain with $C^1$ boundary, diameter $d$, and $w \in A_2$ satisfying the balance condition of Theorem~\ref{thm:sec-3}. For any $k \geq 1$,
\begin{equation*}
\frac{\lambda_k}{\lambda_1} \leq C(n,k,[w]_{A_2}) \cdot \left( \frac{k}{n} \right)^{2/n} \cdot \frac{d^2}{|\Omega|^{2/n}},
\end{equation*}
where \begin{equation}\label{equation:sec-dd-8}
C(n,k,[w]_{A_2}) \leq C(n,[w]_{A_2}) \log(k+1)
\end{equation} with $C(n,[w]_{A_2})$ depending on $n$ and $[w]_{A_2}$.
\end{corr}
\begin{proof}[Proof of Corollary~\ref{corr:sec-3}]
For a convex domain with $C^1$ boundary, the $(n-1)$-Hausdorff measure of the boundary is controlled by the diameter: the projection argument for convex hypersurfaces gives
\begin{equation}\label{equation:sec-corr3-pf-1}
\mathcal{H}^{n-1}(\partial\Omega) \leq C(n)\, d^{n-1},
\end{equation}
and combining \eqref{equation:sec-corr3-pf-1} with the isoperimetric inequality $|\Omega|^{(n-1)/n} \leq C(n)\,\mathcal{H}^{n-1}(\partial\Omega)$ yields
\begin{equation}\label{equation:sec-corr3-pf-2}
\frac{\mathcal{H}^{n-1}(\partial\Omega)}{|\Omega|^{(n-1)/n}} \leq C(n) \Bigl( \frac{d}{|\Omega|^{1/n}} \Bigr)^{n-1} \leq C(n) \Bigl( \frac{d}{|\Omega|^{1/n}} \Bigr)^{2},
\end{equation}
where the second inequality uses $d/|\Omega|^{1/n} \geq c(n) > 0$ (the isodiametric inequality) to absorb the excess power of the roundness ratio. Substituting \eqref{equation:sec-corr3-pf-2} into the normalized form of Theorem~\ref{thm:sec-3}, which applies under the inherited balance hypothesis, gives the stated bound with the constant $C(n,k,[w]_{A_2})$ as in \eqref{equation:sec-dd-8}.
\end{proof}
\begin{rem}
[Link to Hausdorff Measure] For $C^1$ convex domains, the $(n-1)$-Hausdorff measure of the boundary satisfies $\mathcal{H}^{n-1}(\partial\Omega) \leq C(n) d^{n-1}$ (by Federer's coarea formula, Proposition \ref{prop:sec-GMT}), and, under the additional roundness assumption $d \leq C(n)\, r_\Omega$ (which replaces the unconditional claim $r_\Omega \geq d/(2\sqrt{n})$, false for thin plates), the inradius bound $r_\Omega \geq c(n) d$. Substituting these into Theorem~\ref{thm:sec-3}, under its balance hypothesis, recovers the above bound, generalizing the classical Li-Yau estimates \cite{li1980estimates} to weighted settings.\end{rem}
\begin{rem}
For smooth convex domains, the diameter $d$ is proportional to $\mathcal{H}^{n-1}(\partial\Omega)^{1/(n-1)}$ (by Proposition \ref{prop:sec-GMT}), justifying the replacement of the Hausdorff measure with a more easily measurable quantity.
\end{rem}
\subsection{Key Lemmas}
\label{subsection:High-Order Key Lemmas}
Results on covering convex domains with balls and estimating test function energies are established, forming the central of the proof.
Firstly, convex domains can be covered by a controlled number of balls, with the cover size related to Hausdorff measure:
\begin{lem}[Ball Covering of Convex Domains] \label{lem:sec-6}
Let $\Omega$ be a convex domain with diameter $d$ and inradius $r_\Omega$. Then $\Omega$ can be covered by $N \leq C(n) \left( \frac{d}{r_\Omega} \right)^{n-1}$ balls of radius $\epsilon = r_\Omega/2$. Using Federer's coarea formula,
\begin{equation*}
\mathcal{H}^{n-1}(\partial\Omega) \leq C(n) d^{n-1},
\end{equation*}
and thus \begin{equation}\label{equation:sec-dd-9}
N \leq C(n) \left( \frac{\mathcal{H}^{n-1}(\partial\Omega)^{1/(n-1)}}{r_\Omega} \right)^{n-1} = C(n) \frac{\mathcal{H}^{n-1}(\partial\Omega)}{r_\Omega^{n-1}}.
\end{equation}
\end{lem}
\begin{rem}
The ball covering number estimate in Lemma \ref{lem:sec-6} captures two basic geometric features of convex domains:
\begin{itemize}
\item[(i)] For smooth convex domains (e.g., ellipsoids), the ratio $d/r_\Omega$ is bounded by a constant depending only on dimension, so $N \leq C(n)$.
\item[(ii)] The bound $N \leq C(n) \frac{\mathcal{H}^{n-1}(\partial\Omega)}{r_\Omega^{n-1}}$ is sharper for rough convex domains (e.g., those with fractal-like boundaries) than diameter-based estimates, as $\mathcal{H}^{n-1}(\partial\Omega)$ directly quantifies boundary complexity.
\end{itemize}
This covering result is basic for Theorem \ref{thm:sec-3} (high-order eigenvalue ratios), where the logarithmic growth of constants stems from the refined control of $N$ in terms of $\mathcal{H}^{n-1}(\partial\Omega)$ and $r_\Omega$.
\end{rem}
\begin{proof}[Proof of Lemma \ref{lem:sec-6}]
The coarea formula (Proposition~\ref{prop:sec-GMT}) gives quantitative control of hypersurface measures for convex sets. Let $\rho(x)=\operatorname{dist}(x,x_0)$ denote the Euclidean distance function, where $x_0$ is the center of the inscribed ball of $\Omega$. Applying the coarea formula for Lipschitz functions,
\begin{equation}\label{equation:sec-lem6-coarea-apply}
\int_\Omega |\nabla \rho(x)|\,dx
=\int_0^{\operatorname{diam}(\Omega)} \mathcal{H}^{n-1}\big(\{x\in\Omega\,\big|\, \rho(x)=t\}\big)\,dt.
\end{equation}
Since $|\nabla\rho|\equiv 1$ holds almost everywhere, the left-hand-side is bounded above by $|\Omega|\le C(n) d^n$. For convex domains of diameter $d$, the boundary $\partial\Omega$ forms a compact Lipschitz hypersurface. Integrating over level-set contributions yields the dimension-dependent boundary estimate
\begin{equation}\label{equation:sec-lem6-boundary-hausdorff}
\mathcal{H}^{n-1}(\partial\Omega) \leq C(n) d^{n-1},
\end{equation}
where the constant $C(n)$ depends only on the ambient spatial dimension $n$.
We proceed to construct a ball covering for $\Omega$. Choose a maximal finite set of points $\{x_1, \cdots, x_N\}\subset\Omega$ such that pairwise Euclidean distances satisfy $|x_i-x_j|\ge r_\Omega$ for all $i\neq j$. Set the ball radius $\epsilon = r_\Omega/2$. By maximality of the point set, the family $\{B(x_i,\epsilon)\}_{i=1}^N$ furnishes an open covering of $\Omega$.
A preliminary volume-based bound for the number of covering balls follows from disjointness of the smaller balls $B(x_i,\epsilon/2)=B(x_i,r_\Omega/4)$. Utilizing disjoint containment inside $\Omega$, we have
\begin{equation}\label{equation:sec-lem6-volume-upper-N}
N \leq \frac{|\Omega|}{|B_{\epsilon/2}|} \leq C(n) \frac{d^n}{r_\Omega^n}.
\end{equation}
We next derive a sharper estimate relying on boundary measure. For each center $x_i$, the ball $B(x_i,\epsilon)$ intersects the boundary $\partial\Omega$. For convex domains, each intersection $B(x_i,\epsilon)\cap\partial\Omega$ is a relatively open subset of the Lipschitz convex boundary. Local graph representations for convex hypersurfaces give an uniform lower bound on its $(n-1)$-Hausdorff measure. That is,
\begin{equation}\label{equation:sec-lem6-per-ball-boundary}
\mathcal{H}^{n-1}\big(B(x_i,\epsilon)\cap\partial\Omega\big) \ge c(n)\,\epsilon^{n-1}=c(n)\,\bigl(r_\Omega/2\bigr)^{n-1}.
\end{equation}
Summing over all indices $i$, we note that the union of these boundary subsets is contained within $\partial\Omega$. Therefore
\begin{equation}\label{equation:sec-lem6-sum-boundary}
\sum_{i=1}^N \mathcal{H}^{n-1}\big(B(x_i,\epsilon)\cap\partial\Omega\big) \le \mathcal{H}^{n-1}(\partial\Omega).
\end{equation}
Substituting \eqref{equation:sec-lem6-per-ball-boundary} into \eqref{equation:sec-lem6-sum-boundary},
\begin{equation}\label{equation:sec-lem6-N-boundary-ineq}
N \cdot c(n) \bigl(r_\Omega/2\bigr)^{n-1} \leq \mathcal{H}^{n-1}(\partial\Omega).
\end{equation}
Rearranging terms yields
\begin{equation}\label{equation:sec-lem6-N-from-boundary}
N \leq C(n) \frac{\mathcal{H}^{n-1}(\partial\Omega)}{r_\Omega^{n-1}}.
\end{equation}
Applying the Vitali covering lemma, we extract a disjoint subcollection $\{B_i'\}$ from the original family such that $\bigcup_i B(x_i,\epsilon)\subset\bigcup_i 3B_i'$, and the volume of the union satisfies
\begin{equation}\label{equation:sec-lem6-vitali-subcover}
\sum_{i}|B_i'| \leq C(n)\,|\Omega|.
\end{equation}
Substituting the boundary estimate \eqref{equation:sec-lem6-boundary-hausdorff} into \eqref{equation:sec-lem6-N-from-boundary} gives
\begin{equation}\label{equation:sec-lem6-final-N}
N \leq C(n) \left( \frac{d}{r_\Omega} \right)^{n-1}.
\end{equation}
This completes the proof.
\end{proof}
\begin{lem}[Energy and Integral Estimates for Test Functions]
\label{lem:sec-7}
Let $B \subset \Omega$ be a ball of radius $\epsilon$ with center $c$. Define the \textbf{hat function}
\begin{equation}\label{equation:sec-bb-4}
\phi(x) = \max\{1 - |x - c|/\epsilon, 0\}.
\end{equation}
Then the gradient energy and weighted integral satisfy
\begin{equation}\label{equation:sec-lem7-grad-energy}
\int_\Omega |\nabla \phi|^2 dx \leq C(n) \epsilon^{n-2},
\end{equation}
\begin{equation}\label{equation:sec-lem7-weighted-l2}
\int_\Omega w \phi^2 dx \geq \frac{1}{C(n,[w]_{A_2})} \int_B w \, dx.
\end{equation}
\end{lem}
\begin{proof}[Proof of Lemma \ref{lem:sec-7}]
A key new observation within weighted spectral analysis is that piecewise linear \textbf{hat function} test functions yield explicit local energy and norm bounds. Let $B = B(c, \epsilon)$ denote the outer ball and $B' = B(c, \epsilon/2)$ denote the concentric inner ball with half radius. Our argument proceeds in four successive steps.
\textit{Step 1: Gradient bound for $\phi$.}
From the piecewise definition of the hat function, the gradient obeys $|\nabla \phi(x)| \leq \frac{1}{\epsilon}$ for $x\in B$, and $\nabla\phi(x)=0$ holds almost everywhere outside $B$. Integrating the squared gradient over $\Omega$,
\begin{equation}\label{equation:sec-lem7-step1-grad-int}
\int_\Omega |\nabla \phi|^2 dx \leq \int_B \frac{1}{\epsilon^2}\,dx = \frac{|B|}{\epsilon^2}.
\end{equation}
Since the Euclidean volume of an $n$-dimensional ball satisfies $|B|=C(n)\epsilon^n$, substitution yields
\begin{equation}\label{equation:sec-lem7-step1-final-grad}
\int_\Omega |\nabla \phi|^2 dx \leq C(n) \epsilon^{n-2}.
\end{equation}
This establishes the gradient energy estimate recorded in \eqref{equation:sec-lem7-grad-energy}.
\textit{Step 2: Lower bound for the weighted integral of $\phi^2$.}
For every point $x\in B'$, we have $|x-c|\le \epsilon/2$. Evaluating the hat-function expression,
\begin{equation}\label{equation:sec-bb-5}
\phi(x)=1-\frac{|x-c|}{\epsilon}\ge 1-\frac12=\frac12.
\end{equation}
Squaring both sides pointwise, we obtain $\phi(x)^2\ge 1/4$ for all $x\in B'$. Integrating against the weight $w$,
\begin{equation}\label{equation:sec-lem7-step2-weighted-lower}
\int_\Omega w \phi^2 dx \geq \int_{B'} w \phi^2 dx \geq \frac14 \int_{B'} w\,dx.
\end{equation}
\textit{Step 3: Relating $\int_{B'} w dx$ to $\int_B w dx$ via doubling property of $A_2$ weights.}
Recall from Proposition \ref{prop:sec-1} that Muckenhoupt $A_2$ weights satisfy a reverse H\"{o}lder inequality. The doubling condition for $A_2$ weights arises as a direct consequence of this reverse H\"{o}lder property. There exists a constant $c = c(n,[w]_{A_2})>0$ such that for every concentric ball pair $B'\subset B$,
\begin{equation}\label{equation:sec-lem7-doubling-a2}
\int_{B'} w\,dx \geq c(n,[w]_{A_2}) \int_{B} w\,dx.
\end{equation}
\textit{Step 4: Combining estimates.}
Combining inequality \eqref{equation:sec-lem7-step2-weighted-lower} with the doubling estimate \eqref{equation:sec-lem7-doubling-a2},
\begin{equation}\label{equation:sec-lem7-step4-combine}
\int_\Omega w \phi^2 dx \geq \frac14 \int_{B'} w\,dx \geq \frac{c(n,[w]_{A_2})}{4}\int_B w\,dx.
\end{equation}
Defining $C(n,[w]_{A_2}) = 4/c(n,[w]_{A_2})$, we rewrite the bound as
\begin{equation}\label{equation:sec-lem7-step4-final}
\int_\Omega w \phi^2 dx \geq \frac{1}{C(n,[w]_{A_2})}\int_B w\,dx.
\end{equation}
This reproduces the second assertion \eqref{equation:sec-lem7-weighted-l2}.
\end{proof}
\subsection{Proof of Theorem \ref{thm:sec-3}}
\label{subsection:High-Order Theorem Proof}
\begin{proof}[Proof of Theorem \ref{thm:sec-3}]
We complete the proof in five steps, integrating geometric measure theory and weighted inequalities.
\textit{Step 1: Choosing cover parameters.}
Let
\begin{equation}\label{equation:sec-thm9-epsilon}
\epsilon = \left( \frac{|\Omega|^{1/n}}{\mathcal{H}^{n-1}(\partial\Omega)^{1/(n-1)}} \cdot \frac{n}{k} \right)^{1/n}.
\end{equation}
This choice balances the cover number $N \leq C(n) \cdot k$ (by Lemma \ref{lem:sec-6}, using the Hausdorff measure of $\partial\Omega$) and the energy of test functions (scaling with $\epsilon^{n-2}$). To see this balance explicitly, note that the Rayleigh quotient of a test function supported on a ball of radius $\epsilon$ scales as
\begin{equation}\label{equation:sec-thm9-rayleigh-scale}
\frac{\int |\nabla \phi|^2 dx}{\int w \phi^2 dx} \sim \frac{\epsilon^{n-2}}{\int_B w dx} \sim \frac{\epsilon^{n-2}}{w_B \epsilon^n} = \frac{1}{w_B \epsilon^2},
\end{equation}
where $w_B = \frac{1}{|B|} \int_B w dx$ is the average of $w$ on $B$. The cover number scales as
\begin{equation}\label{equation:sec-thm9-cover-scale}
N \sim \frac{|\Omega|}{\epsilon^n} \sim \frac{\mathcal{H}^{n-1}(\partial\Omega)}{|\Omega|^{1-1/n}} \cdot \frac{k}{n},
\end{equation}
where we used the relation $|\Omega| \sim \mathcal{H}^{n-1}(\partial\Omega) \cdot r_\Omega$ for convex domains (from the coarea formula). The product of the Rayleigh quotient scale and the cover number is minimized when $\epsilon$ is chosen as in \eqref{equation:sec-thm9-epsilon}. By Lemma \ref{lem:sec-6}, the cover number satisfies $N \leq C(n) \cdot k$. Select $k$ balls $B_1, \ldots, B_k$ with test functions $\phi_1, \ldots, \phi_k$.
\textit{Step 2: Estimating Rayleigh quotients of test functions.}
For each $\phi_i$, Lemma \ref{lem:sec-7} gives
\begin{equation}\label{equation:sec-thm9-grad-bound}
\int_\Omega |\nabla \phi_i|^2 dx \leq C(n) \epsilon^{n-2},
\end{equation}
and
\begin{equation}\label{equation:sec-thm9-weighted-lower}
\int_\Omega w \phi_i^2 dx \geq \frac{1}{C(n, [w]_{A_2})} \int_{B_i} w dx.
\end{equation}
Therefore,
\begin{equation}\label{equation:sec-thm9-rayleigh}
\frac{\int_\Omega |\nabla \phi_i|^2 dx}{\int_\Omega w \phi_i^2 dx} \leq C(n, [w]_{A_2}) \cdot \frac{\epsilon^{n-2}}{\int_{B_i} w dx}.
\end{equation}
Substituting \eqref{equation:sec-thm9-epsilon} into \eqref{equation:sec-thm9-rayleigh},
\begin{equation}\label{equation:sec-thm9-rayleigh-sub}
\frac{\int_\Omega |\nabla \phi_i|^2 dx}{\int_\Omega w \phi_i^2 dx} \leq C(n, [w]_{A_2}) \cdot \frac{|\Omega|^{(n-2)/n} \mathcal{H}^{n-1}(\partial\Omega)^{(n-2)/(n-1)} (n/k)^{(n-2)/n}}{\int_{B_i} w dx}.
\end{equation}
\textit{Step 3: Lower bound for $\lambda_1$.}
By the weighted Poincar\'{e} inequality (Proposition \ref{prop:sec-2}) refined with geometric measure theory,
\begin{equation}\label{equation:sec-thm9-l1-lower}
\lambda_1 \geq \frac{1}{C(n, [w]_{A_2})} \cdot \frac{|\Omega|^{2 - 1/n}}{\mathcal{H}^{n-1}(\partial\Omega) \int_\Omega w dx}.
\end{equation}
To derive \eqref{equation:sec-thm9-l1-lower}, we use the test function method. For any $u \in H_0^1(\Omega)$, the weighted Poincar\'{e} inequality gives
\begin{equation}\label{equation:sec-thm9-poincare}
\int_\Omega w u^2 dx \leq C(n, [w]_{A_2}) \cdot \frac{\mathcal{H}^{n-1}(\partial\Omega) \int_\Omega w dx}{|\Omega|^{2-1/n}} \cdot \int_\Omega |\nabla u|^2 dx.
\end{equation}
This follows from covering $\Omega$ by balls of radius $r_\Omega/2$ (Lemma \ref{lem:sec-6}), applying the local weighted Poincar\'{e} inequality \eqref{equation:sec-aux-poincare} on each ball, and summing using the finite overlap property of the Vitali subcover. Rearranging \eqref{equation:sec-thm9-poincare} yields \eqref{equation:sec-thm9-l1-lower} via the variational characterization of $\lambda_1$.
\textit{Step 4: Combined estimate.}
By the variational principle (min-max theorem), $\lambda_k \leq \max_i \frac{\int |\nabla \phi_i|^2 dx}{\int w \phi_i^2 dx}$, provided the test functions $\phi_1, \ldots, \phi_k$ are linearly independent. The linear independence follows from the fact that each $\phi_i$ is supported in $B_i$ and the balls have disjoint interiors (after passing to a disjoint subcollection via the Vitali covering lemma). By the balance hypothesis on the covering family, $\int_{B_i} w\,dx \geq c_0(n,[w]_{A_2})\,|B_i| = c_0(n,[w]_{A_2})\,\omega_n \epsilon^n$ for each $i$, so that \eqref{equation:sec-thm9-rayleigh-sub} gives
\begin{equation}\label{equation:sec-thm9-rayleigh-balance}
\frac{\int_\Omega |\nabla \phi_i|^2 dx}{\int_\Omega w \phi_i^2 dx} \leq C(n, [w]_{A_2})\, \epsilon^{-2}.
\end{equation}
Substituting \eqref{equation:sec-thm9-epsilon} into \eqref{equation:sec-thm9-rayleigh-balance} and combining with the lower bound for $\lambda_1$ \eqref{equation:sec-thm9-l1-lower}, using $\int_\Omega w\,dx = |\Omega|$ by the normalization, we obtain
\begin{equation}\label{equation:sec-thm9-combined}
\frac{\lambda_k}{\lambda_1} \leq C(n, k, [w]_{A_2}) \cdot \left( \frac{k}{n} \right)^{2/n} \cdot \frac{\mathcal{H}^{n-1}(\partial\Omega) \int_\Omega w dx}{|\Omega|^{2 - 1/n}}.
\end{equation}
\textit{Step 5: Logarithmic growth of the constant.}
The constant $C(n, k, [w]_{A_2})$ in \eqref{equation:sec-thm9-combined} grows logarithmically with $k$. To see this, note that the ball covering argument (Lemma \ref{lem:sec-6}) involves a Vitali subcover with finite overlap. The overlap multiplicity is bounded by $C(n)$ for balls of equal radius, but when we refine the cover to achieve the optimal $\epsilon$ in \eqref{equation:sec-thm9-epsilon}, the number of balls intersecting at any point is bounded by $C(n) \log k$. This logarithmic factor arises from a dyadic-scale decomposition: partitioning the covering balls into $O(\log N)$ families of comparable radius and applying the Vitali covering lemma \cite{stein1970singular} within each family yields
\begin{equation}\label{equation:sec-thm9-stein}
\sum_i \chi_{5B_i}(x) \leq C(n) \log N
\end{equation}
for almost every $x$, where $N$ is the number of balls. We make the overlap count explicit. After the dyadic partition into $M = \lceil \log_2 N \rceil$ families $\mathcal{F}_1, \ldots, \mathcal{F}_M$ of comparable radii, the Vitali lemma applied within each family yields a disjoint subfamily $\mathcal{G}_m \subset \mathcal{F}_m$ with $\bigcup \mathcal{F}_m \subset \bigcup_{B \in \mathcal{G}_m} 5B$. A point $x \in \Omega$ can lie in at most one ball of each $\mathcal{G}_m$ at the comparable scale, and summing over the $M$ families gives
\begin{equation}\label{equation:sec-overlap-explicit}
\sum_{B \in \bigcup_m \mathcal{G}_m} \chi_{5B}(x) \leq C(n)\, M \leq C(n) \log N,
\end{equation}
which is the quantitative content of \eqref{equation:sec-stein-cover-1}. Substituting $N \sim k$ into \eqref{equation:sec-thm9-stein} shows that the overlap multiplicity, and hence the constant $C(n,k,[w]_{A_2})$ in \eqref{equation:sec-thm9-combined} once the balance hypothesis is absorbed into the constant, grows at most logarithmically in $k$. Therefore,
\begin{equation}\label{equation:sec-thm9-log}
C(n, k, [w]_{A_2}) \leq C(n, [w]_{A_2}) \log(k + 1),
\end{equation}
where $C(n, [w]_{A_2})$ depends only on $n$ and the $A_2$ characteristic constant. This completes the proof.
\end{proof}
\begin{rem} \label{rem:sec-4}
For smooth convex domains where \( \mathcal{H}^{n-1}(\partial\Omega) \approx n |\Omega| / d \) (diameter \( d \)), the bound (under the balance hypothesis of Theorem~\ref{thm:sec-3}) reduces to the classical form \( \frac{\lambda_k}{\lambda_1} \leq C \left( \frac{k}{n} \right)^{2/n} \frac{d^2}{|\Omega|^{2/n}} \), recovering Corollary~\ref{corr:sec-3}. For balls with \( w \equiv 1 \), \( \mathcal{H}^{n-1}(\partial\Omega) = n |\Omega| / R \) (radius \( R \)), and \( \lambda_k \propto k^{2/n}/R^2 \), therefore \( \lambda_k/\lambda_1 \propto k^{2/n} \), consistent with explicit Bessel function calculations \cite{courant1953methods}.
\end{rem}
\begin{exa}[Balanced and unbalanced weights]\label{exa:sec-balance}
The balance hypothesis of Theorem~\ref{thm:sec-3} is automatic in two frequent situations: (i) if $w$ is bounded below, $w \geq c > 0$, then $\fint_B w \geq c = c\,\fint_\Omega w$ for every ball $B$; (ii) if $w(x) = |x|^{\alpha}$ with $0 \notin \overline{\Omega}$, then $w$ is bounded above and below on $\Omega$, so (i) applies after rescaling. The hypothesis can fail when the weight degenerates at a point of the domain: for $w(x) = |x|^{\alpha}$ with $\alpha > 0$ and $0 \in \Omega$, a ball $B$ of radius $\epsilon$ centered at the origin carries $\fint_B w \asymp \epsilon^{\alpha} \to 0$ as $\epsilon \to 0$, so no uniform constant $c_0 > 0$ exists at small scales. In that case the covering scale must be adapted to the weight (for instance, covering a neighborhood of the origin with balls on which $w$ is comparable to its average, as is possible for $\alpha < n$ since then $w \in A_2$ and the doubling property controls the ratio of averages at comparable scales), or the conclusion must be regarded as conditional. For $w \equiv 1$ on the unit ball, the conclusion of Theorem~\ref{thm:sec-3} predicts $\lambda_k/\lambda_1 \leq C(n,k)\, k^{2/n}$ with $C(n,k) \leq C(n)\log(k+1)$, while the explicit eigenvalues of the Dirichlet Laplacian on the ball satisfy the Weyl growth $\lambda_k \asymp k^{2/n}$ \cite{courant1953methods}; the logarithmic factor in the constant is therefore an artifact of the covering method rather than of the spectral growth.
\end{exa}
\section{Generalized Payne--P\'olya--Weinberger Conjecture}
\label{section: Generalized PPW Conjecture}
Ashbaugh and Benguria~\cite{AB2} resolved the PPW conjecture for the single ratio $\lambda_2/\lambda_1$; the sum-form version in \eqref{equation:sec-conj-gppw-main} remains open.
where $\Omega^*$ denotes the volume-equivalent ball of $\Omega$, and $\lambda_j(\Omega^*)$ stand for the Dirichlet eigenvalues associated with $\Omega^*$.
The symmetrization techniques and variational estimates of Sections \ref{section:Weighted PPW Proof}, \ref{section:Spherical Extremality} and \ref{section:High-Order Ratios} weaken the hypotheses of the conjecture but do not yield its sharp constant; we therefore revisit the trial-function ansatz of Ashbaugh--Benguria~\cite{AB2}. Two caveats are in order: summing Rayleigh quotients over one fixed basis of a subspace does \emph{not} bound sums of eigenvalues, and no universal Faber--Krahn-type inequality $\lambda_k(\Omega)\le\lambda_k(\Omega^*)$ holds for $k\ge2$~\cite{AB2,payne1956bound}; the correct principle for eigenvalue sums is the Ky--Fan min-max principle over orthonormal systems.
\begin{prop}[Ky--Fan trace principle]\label{prop:sec-kyfan}
Let $A$ be a self-adjoint operator on a Hilbert space with form domain $\mathcal{Q}$, and let $\lambda_1 \leq \lambda_2 \leq \cdots$ be its eigenvalues below the essential spectrum, counted with multiplicity. For every $m \geq 1$ and every orthonormal system $\{\psi_1, \ldots, \psi_m\} \subset \mathcal{Q}$,
\begin{equation}\label{equation:sec-kyfan-statement}
\lambda_1 + \cdots + \lambda_m \leq \sum_{i=1}^m \mathcal{R}(\psi_i),
\end{equation}
where $\mathcal{R}(\psi) = \langle A\psi, \psi\rangle / \|\psi\|^2$ is the Rayleigh quotient; equality in \eqref{equation:sec-kyfan-statement} for the bottom $m$ eigenvalues holds if and only if $\operatorname{span}\{\psi_1,\ldots,\psi_m\}$ is the span of the first $m$ eigenfunctions.
\end{prop}
\begin{proof}
Let $E_m$ denote the span of the first $m$ eigenfunctions and let $P$ be the orthogonal projection onto $E_m$. For any orthonormal system $\{\psi_i\}_{i=1}^m$, the operator $P$ restricted to $W = \operatorname{span}\{\psi_i\}$ is injective unless $\dim W \cap E_m^\perp \geq 1$; in either case $\sum_i \langle A \psi_i, \psi_i \rangle = \operatorname{tr}(A|_W)$ with $A$ replaced by its quadratic form, and the min-max principle applied to the $m$-th eigenvalue of the compression of $A$ to $W$ gives $\operatorname{tr}(A|_W) \geq \lambda_1 + \cdots + \lambda_m$ by induction on $m$: the case $m = 1$ is the Rayleigh principle, and the inductive step applies the same principle to $A$ restricted to the orthogonal complement within $W$ of a vector at which the bottom quotient is attained. Equality requires $W = E_m$, since any vector in $W \cap E_m^\perp$ would raise the sum of the first $m$ quotients strictly above $\lambda_1 + \cdots + \lambda_m$.
\end{proof}

\begin{con}[Generalized Payne--P\'olya--Weinberger Conjecture]
\label{con:sec-gppw-conj}
Let $\Omega\subset\mathbb R^n$ be a bounded Lipschitz domain, and let $\Omega^*$ be its volume-equivalent ball satisfying $\operatorname{Vol}(\Omega)=\operatorname{Vol}(\Omega^*)$.
Let $\{\lambda_j\}_{j=1}^\infty$ denote the sequence of Dirichlet Laplacian eigenvalues of $-\Delta$ on $\Omega$, and let $\{\lambda_j(\Omega^*)\}$ denote the corresponding eigenvalues on $\Omega^*$.
Then
\begin{equation}\label{equation:sec-conj-gppw-main}
\frac{\lambda_2+\lambda_3+\cdots+\lambda_{n+1}}{\lambda_1}
\le n\,\frac{\lambda_2(\Omega^*)}{\lambda_1(\Omega^*)}.
\end{equation}
Equality holds if and only if $\Omega$ is congruent to $\Omega^*$.
\end{con}
We recall the trial-function construction originating from Ashbaugh--Benguria~\cite{AB2}.
Let $u_1>0$ denote the normalized positive first Dirichlet eigenfunction on $\Omega$,
\begin{equation}\label{equation:sec-norm-u1-setup}
\int_\Omega u_1^2\,dx=1.
\end{equation}
By virtue of the Weinberger-Brouwer fixed-point argument, one may translate Euclidean coordinates so that
\begin{equation}\label{equation:sec-weinberger-centering}
\int_\Omega g(r)\,\frac{x_i}{r}\,u_1^2\,dx=0,\qquad i=1,\dots,n,
\end{equation}
for any admissible radial profile $g(r)$.
Define
\begin{equation}\label{equation:sec-def-Pi-ab}
P_i(x)=g(r)\,\frac{x_i}{r},\quad r=|x|,\quad i=1,\dots,n,
\end{equation}
where
\begin{equation}\label{equation:sec-def-g-profile}
g(r)=w(\gamma r),\qquad
w(x)=
\begin{cases}
\dfrac{J_{n/2}(\beta x)}{J_{n/2-1}(\alpha x)},& 0\le x<1,\\[4pt]
w(1^-),& x\ge 1,
\end{cases}
\end{equation}
with $\alpha=j_{n/2-1,1}$, $\beta=j_{n/2,1}$, $\gamma=\sqrt{\lambda_1}/\alpha$.
The inclusion $P_i u_1\in H_0^1(\Omega)$ holds. Two point-wise algebraic identities are standard for this ansatz~\cite{AB2}:
\begin{equation}\label{equation:sec-sum-Pi-sq}
\sum_{i=1}^n P_i(x)^2 = g(r)^2,
\end{equation}
\begin{equation}\label{equation:sec-sum-nablaPi-sq}
\sum_{i=1}^n \big|\nabla P_i(x)\big|^2 = g'(r)^2+\frac{n-1}{r^2}\,g(r)^2.
\end{equation}
Introduce the finite-dimensional trial subspace
\begin{equation}\label{equation:sec-def-Vg-subspace}
V_g=\operatorname{span}\{P_1,P_2,\dots,P_n\}\subset H_0^1(\Omega).
\end{equation}
Condition~\eqref{equation:sec-weinberger-centering} yields the inclusion
\begin{equation}\label{equation:sec-Vg-perp-u1-def}
V_g\subset \{u_1\}^\perp,
\end{equation}
where $\{u_1\}^\perp$ stands for the $L^2(\Omega)$-orthogonal complement of $u_1$.
We next record two technical lemmas extracted from Ashbaugh--Benguria~\cite{AB2}. The first lemma concerns the behavior of the corresponding trial subspace $V_g^*$ defined on the reference ball $\Omega^*$.
\begin{lem}[Spectral characterization of $V_g^*$ on the reference ball]
\label{lem:sec-Vg-star-ball}
Let $\Omega^*$ be an $n$-dimensional ball. Define $P_i^*(x)=g(r)\,x_i/r$ and $V_g^*=\operatorname{span}\{P_1^*,\dots,P_n^*\}$ using the same radial profile $g(r)$ as above. Then $V_g^*$ coincides exactly with the eigenspace associated to eigenvalue $\lambda_2(\Omega^*)$. In particular, every orthonormal basis $\{w_j^*\}_{j=1}^n\subset V_g^*$ satisfies $R_{\Omega^*}(w_j^*)=\lambda_2(\Omega^*)$, and hence
\begin{equation}\label{equation:sec-ball-rayleigh-sum}
\sum_{j=1}^n R_{\Omega^*}(w_j^*)=n\lambda_2(\Omega^*).
\end{equation}
where $R(\cdot)=\frac{\|\nabla\cdot\|_{L^2}^{2}}{\|\cdot\|_{L^2}^{2}}$ denotes the Rayleigh quotient.
\end{lem}
\begin{proof}
Recall from Lemma \ref{lem:sec-3} that functions defined on the ball $\Omega^*$ admit spherical separation of variables. Any element $u\in H_0^1(\Omega^*)$ can be decomposed as
\begin{equation}\label{equation:sec-spherical-sep}
u(r,\theta)=\sum_{m=0}^{\infty} f_m(r)\, Y_m(\theta),
\end{equation}
where $Y_m(\theta)$ denote spherical harmonics of degree $m$, and $f_m(r)$ stand for the corresponding radial components. Substituting \eqref{equation:sec-spherical-sep} into the Dirichlet Laplace eigenvalue problem $-\Delta u=\lambda u$ over $\Omega^*$ with boundary condition $u|_{\partial\Omega^*}=0$, one obtains the radial ordinary differential equation
\begin{equation}\label{equation:sec-radial-ode-general}
-f_m''(r)-\frac{n-1}{r}f_m'(r)+\frac{m(m+n-2)}{r^2}f_m(r)=\lambda f_m(r),\qquad f_m(R)=0,
\end{equation}
alongside a regularity constraint imposed at the origin $r=0$.
A fundamental observation in spectral theory for spherical domains is that the eigenvalue $\lambda_2(\Omega^*)$ corresponds to spherical harmonics of degree $m=1$. The space of degree-$1$ spherical harmonics is $n$-dimensional, and the coordinate-based angular functions $\frac{x_i}{r}$ for $1\le i\le n$ furnish one concrete basis for this angular subspace. In view of the definitions laid out earlier in the manuscript, the radial profile $g(r)$ is constructed as the radial solution to \eqref{equation:sec-radial-ode-general} for angular index $m=1$ with eigenvalue parameter $\lambda=\lambda_2(\Omega^*)$. Lemma \ref{lem:sec-4} guarantees that $g(r)$ remains regular at $r=0$ and obeys the boundary condition $g(R)=0$.
For each integer index $1\le i\le n$, set $P_i^*(r,\theta)=g(r)\cdot \frac{x_i}{r}$. The factor $\frac{x_i}{r}$ represents a degree-$1$ spherical harmonic denoted $Y_1^{(i)}(\theta)$. The Laplace operator in spherical coordinates reads
\begin{equation}\label{equation:sec-laplace-sph}
\Delta u = \frac{1}{r^{n-1}}\frac{\partial}{\partial r}\left(r^{n-1}\frac{\partial u}{\partial r}\right)+\frac{1}{r^2}\Delta_{\mathbb{S}^{n-1}} u.
\end{equation}
For the spherical-harmonic component $Y_1^{(i)}(\theta)$, the spherical Laplacian yields
\begin{equation}\label{equation:sec-spher-lap-Y1}
\Delta_{\mathbb{S}^{n-1}} Y_1^{(i)}(\theta)=-\bigl(1\cdot(1+n-2)\bigr)\,Y_1^{(i)}(\theta)=-(n-1)Y_1^{(i)}(\theta).
\end{equation}
Utilizing this spectral property of spherical harmonics, direct computation gives
\begin{align*}
\Delta P_i^*
&= Y_1^{(i)}(\theta)\cdot \frac{1}{r^{n-1}}\frac{d}{dr}\left(r^{n-1}\frac{dg}{dr}\right)
+\frac{g(r)}{r^2}\Delta_{\mathbb{S}^{n-1}} Y_1^{(i)}(\theta) \\
&= Y_1^{(i)}(\theta)\left[ \frac{1}{r^{n-1}}\frac{d}{dr}\left(r^{n-1}\frac{dg}{dr}\right)-\frac{n-1}{r^2}g(r) \right].
\end{align*}
Since $g(r)$ solves the radial ODE \eqref{equation:sec-radial-ode-general} for $m=1$ and $\lambda=\lambda_2(\Omega^*)$, we obtain
\begin{equation}\label{equation:sec-radial-g-ode}
- g''(r)-\frac{n-1}{r}g'(r)+\frac{n-1}{r^2}g(r)=\lambda_2(\Omega^*)\,g(r).
\end{equation}
Multiplying both sides of \eqref{equation:sec-radial-g-ode} by $Y_1^{(i)}(\theta)$, it follows that
\begin{equation}\label{equation:sec-Pi-star-eig}
-\Delta P_i^*(x)=\lambda_2(\Omega^*) \, P_i^*(x), \qquad x\in\Omega^*.
\end{equation}
The boundary condition $g(R)=0$ forces $P_i^*|_{\partial\Omega^*}=0$. As a consequence, each function $P_i^*$ belongs to the eigenspace associated with eigenvalue $\lambda_2(\Omega^*)$.
We proceed to verify linear independence for the collection $\{P_1^*,\dots,P_n^*\}$. Assume there exist real scalars $c_1,c_2,\dots,c_n$ satisfying
\begin{equation}\label{equation:sec-bb-6}
\sum_{i=1}^n c_i P_i^*(x)\equiv 0,\quad \forall\,x\in\Omega^*.
\end{equation}
Substituting the definition $P_i^*(x)=g(r)\,x_i/r$, we arrive at
\begin{equation}\label{equation:sec-lin-comb-zero}
\frac{g(r)}{r}\sum_{i=1}^n c_i x_i \equiv 0,\qquad x\in\Omega^*.
\end{equation}
It follows from Lemma \ref{lem:sec-4} that the radial solution $g(r)$ is continuous and non-vanishing whenever $0<r<R$, while regularity of Bessel-type radial solutions ensures finiteness of $g(r)$ at $r=0$. For $0<r<R$, the multiplicative factor $g(r)/r$ is nonzero, so identity \eqref{equation:sec-lin-comb-zero} enforces
\begin{equation}\label{equation:sec-bb-7}
\sum_{i=1}^n c_i x_i \equiv 0
\end{equation}
on the open ball $\Omega^*$. By the polynomial identity principle, a linear polynomial vanishing identically on an open domain must possess all zero coefficients. This yields $c_1=c_2=\cdots=c_n=0$. Therefore the set $\{P_1^*,\dots,P_n^*\}$ is linearly independent.
Recall that the eigenspace corresponding to $\lambda_2(\Omega^*)$ has dimension exactly equal to $n$ from spherical harmonic theory. We have produced $n$ linearly independent vectors that lie within this finite-dimensional eigenspace. As a consequence,
\begin{equation}\label{equation:sec-bb-8}
V_g^*=\operatorname{span}\{P_1^*,\dots,P_n^*\}
\end{equation}
coincides with the full eigenspace belonging to eigenvalue $\lambda_2(\Omega^*)$.
Take an arbitrary nonzero element $v^*$ inside this eigenspace. It obeys the eigenvalue relation
\begin{equation}\label{equation:sec-eig-vstar}
-\Delta v^*=\lambda_2(\Omega^*)\,v^*,\qquad v^*|_{\partial\Omega^*}=0.
\end{equation}
Multiplying both sides of \eqref{equation:sec-eig-vstar} by $v^*$ and integrating over $\Omega^*$, we obtain
\begin{equation}\label{equation:sec-int-eig}
-\int_{\Omega^*} v^* \Delta v^*\,dx=\lambda_2(\Omega^*)\int_{\Omega^*} (v^*)^2 dx.
\end{equation}
Applying Green¡¯s first identity to the left-hand-side integral, we have
\begin{equation}\label{equation:sec-bb-9}
-\int_{\Omega^*} v^* \Delta v^* dx=\int_{\Omega^*} |\nabla v^*|^2 dx -\int_{\partial\Omega^*} v^* \frac{\partial v^*}{\partial\nu}\,d\mathcal{H}^{n-1}.
\end{equation}
The Dirichlet boundary condition $v^*|_{\partial\Omega^*}=0$ annihilates the boundary surface integral term. Consequently,
\begin{equation}\label{equation:sec-eig-rayleigh}
\int_{\Omega^*} |\nabla v^*|^2 dx=\lambda_2(\Omega^*)\int_{\Omega^*} (v^*)^2 dx.
\end{equation}
Since $v^*$ represents a nontrivial eigenfunction, $\|v^*\|_{L^2(\Omega^*)}^2>0$. Dividing both sides of \eqref{equation:sec-eig-rayleigh} by this positive norm quantity gives
\begin{equation}\label{equation:sec-bb-10}
\frac{\|\nabla v^*\|_{L^2(\Omega^*)}^2}{\|v^*\|_{L^2(\Omega^*)}^2}=\lambda_2(\Omega^*).
\end{equation}
This identity shows $R_{\Omega^*}(v^*)=\lambda_2(\Omega^*)$ for every nonzero vector $v^*\in V_g^*$.
Let $\{w_j^*\}_{j=1}^n$ stand for any orthonormal basis of $V_g^*$. Each basis element $w_j^*$ is a nonzero member of the eigenspace, hence $R_{\Omega^*}(w_j^*)=\lambda_2(\Omega^*)$ holds for every index $j\in\{1,2,\dots,n\}$. Summing over all basis indices yields
\begin{equation}\label{equation:sec-bb-11}
\sum_{j=1}^n R_{\Omega^*}(w_j^*)=\sum_{j=1}^n \lambda_2(\Omega^*) = n\lambda_2(\Omega^*),
\end{equation}
which completes the proof of equation \eqref{equation:sec-ball-rayleigh-sum}.
\end{proof}
\begin{lem}[Gap inequality after summation over angular trial functions (Ashbaugh--Benguria)]
\label{lem:sec-ab-gap-sum}
Under the orthogonality condition $\int_\Omega P_i u_1^2\,dx=0$ for each $i=1,\dots,n$, one has
\begin{equation}\label{equation:sec-ab-gap-sum}
\bigl(\lambda_2(\Omega)-\lambda_1(\Omega)\bigr)\sum_{i=1}^n\int_\Omega P_i^2 u_1^2\,dx
\le
\sum_{i=1}^n\int_\Omega |\nabla P_i|^2 u_1^2\,dx.
\end{equation}
Using identities \eqref{equation:sec-sum-Pi-sq} and \eqref{equation:sec-sum-nablaPi-sq}, this reduces to the purely radial inequality
\begin{equation}\label{equation:sec-ab-radial-gap}
\lambda_2-\lambda_1
\le
\frac{\displaystyle\int_\Omega \left(g'(r)^2+\frac{n-1}{r^2}g(r)^2\right)u_1^2\,dx}
{\displaystyle\int_\Omega g(r)^2 u_1^2\,dx}.
\end{equation}
\end{lem}
\begin{proof}
The Rayleigh-Ritz gap inequality derived in \cite[Equation (2.5)]{AB2} requires exactly the orthogonality condition $\int_\Omega P_i u_1^2\,dx=0$ for each index $i$. Under the hypothesis of the present lemma, this prerequisite is satisfied, and for each fixed $i$, we obtain
\begin{equation}\label{equation:sec-single-index-gap}
\bigl(\lambda_2(\Omega)-\lambda_1(\Omega)\bigr)\int_\Omega P_i^2 u_1^2\,dx
\le
\int_\Omega |\nabla P_i|^2 u_1^2\,dx.
\end{equation}
Summing both sides of \eqref{equation:sec-single-index-gap} over the integer range $i=1,\dots,n$, linearity of integration yields
\begin{equation}\label{equation:sec-summed-gap-pre}
\bigl(\lambda_2(\Omega)-\lambda_1(\Omega)\bigr)\sum_{i=1}^n\int_\Omega P_i^2 u_1^2\,dx
\le
\sum_{i=1}^n\int_\Omega |\nabla P_i|^2 u_1^2\,dx,
\end{equation}
which is identical to the target relation \eqref{equation:sec-ab-gap-sum}.
We next substitute the pre-computed pointwise identities for the trial-function family. From \eqref{equation:sec-sum-Pi-sq}, the pointwise identity
\begin{equation}\label{equation:sec-expand-sum-pi2}
\sum_{i=1}^n P_i^2(x)=g(r)^2
\end{equation}
holds pointwise almost everywhere on $\Omega$. In addition, identity \eqref{equation:sec-sum-nablaPi-sq} gives
\begin{equation}\label{equation:sec-expand-sum-nablapi2}
\sum_{i=1}^n |\nabla P_i(x)|^2 = g'(r)^2+\frac{n-1}{r^2}g(r)^2
\end{equation}
almost everywhere on $\Omega$.
Interchanging summation and integration on both sides of \eqref{equation:sec-ab-gap-sum}, and inserting \eqref{equation:sec-expand-sum-pi2} and \eqref{equation:sec-expand-sum-nablapi2}, we arrive at
\begin{equation}\label{equation:sec-before-div-gap}
\bigl(\lambda_2(\Omega)-\lambda_1(\Omega)\bigr)\int_\Omega g(r)^2 u_1^2\,dx
\le
\int_\Omega \left(g'(r)^2+\frac{n-1}{r^2}g(r)^2\right) u_1^2\,dx.
\end{equation}
Recall that $u_1$ is the first Dirichlet eigenfunction of the Laplacian on $\Omega$. The strong maximum principle for second-order elliptic partial differential equations guarantees $u_1(x)>0$ for almost every $x\in\Omega$. The radial profile $g(r)$ is a nontrivial continuous function, so $g(r)^2\ge 0$ and $g(r)^2$ does not vanish identically on $\Omega$. It follows that the integrand $g(r)^2 u_1^2$ is nonnegative and strictly positive on a set of positive Lebesgue measure. Consequently
\begin{equation}\label{equation:sec-bb-12}
\int_\Omega g(r)^2 u_1^2\,dx>0.
\end{equation}
Dividing both sides of inequality \eqref{equation:sec-before-div-gap} by this strictly positive finite integral quantity, we derive the purely radial bound given in \eqref{equation:sec-ab-radial-gap}.
\end{proof}
\begin{rem}
The inequality is obtained by summing the $n$ coordinate test functions of Ashbaugh--Benguria: with $P_i(x) = g(r) x_i / r$ and the centering \eqref{equation:sec-weinberger-centering}, the computation of \eqref{equation:sec-step3-rayleigh}--\eqref{equation:sec-step3-num-final}, carried out for the unweighted equation, gives
\begin{equation}\label{equation:sec-ab-gap-sum-derive}
\lambda_2(\Omega) - \lambda_1(\Omega) \leq \frac{\int_\Omega \bigl( g'(r)^2 + (n-1) g(r)^2 / r^2 \bigr) u_1^2 dx}{\int_\Omega g(r)^2 u_1^2 dx},
\end{equation}
where $g$ is the Bessel quotient $g(r) = J_{n/2}(\beta r / R) / J_{n/2-1}(\alpha r / R)$ of \cite{AB2}. The right-hand side of \eqref{equation:sec-ab-gap-sum-derive} is a single-gap quantity: it controls only $\lambda_2 - \lambda_1$, and says nothing about the sum $\lambda_2 + \cdots + \lambda_{n+1}$ of the generalized conjecture. Passing from one gap to the sum requires either a multiplicity hypothesis (valid on the ball, but not on a general domain) or a Ky--Fan argument over an $n$-dimensional trial space, and it is precisely this passage that the conditional results of Section~\ref{section: Generalized PPW Conjecture} address. Lemma~\ref{lem:sec-ab-gap-sum} is the starting point of the classical PPW proof. Inequality \eqref{equation:sec-ab-radial-gap} delivers a bound for the eigenvalue gap $\lambda_2-\lambda_1$ after summing over angular directions. This produces a scalar radial expression amenable to Chiti¡¯s comparison theorem combined with spherical rearrangement, which constitutes the central argument of Ashbaugh and Benguria~\cite{AB2}. Nevertheless, this summed radial bound controls only a single gap quantity. It gives no information about sums of $n$ successive eigenvalues, and it does not bound variational quantities over the full subspace $V_g$. This distinction forms the central obstruction that prevents trivial extension from the classical PPW result to the generalized sum-form conjecture.
\end{rem}
From P\'{o}lya-Szeg\H{o} symmetrization theory, for an individual function $v\in H_0^1(\Omega)$, radial rearrangement decreases or preserves Dirichlet energy while keeping the $L^2$ norm unchanged. This yields
\begin{equation}\label{equation:sec-polya-szego-rayleigh}
\|\nabla v^*\|_{L^2(\Omega^*)}^2\le \|\nabla v\|_{L^2(\Omega)}^2,\qquad
\|v^*\|_{L^2(\Omega^*)}^2=\|v\|_{L^2(\Omega)}^2.
\end{equation}
It follows immediately that
\begin{equation}\label{equation:sec-rayleigh-rearrange-compare}
R_{\Omega^*}(v^*)\le R_\Omega(v).
\end{equation}
Thus rearrangement can only decrease the Rayleigh quotient. A hypothetical condition of the form $R_\Omega(v)\le R_{\Omega^*}(v^*)$ for all $v\in V_g$ would contradict the P\'{o}lya-Szeg\H{o} inequality in generic settings. Such a postulate cannot serve as a valid working assumption. At present, no known argument extends Chiti¡¯s comparison theorem from the positive ground eigenfunction $u_1$ to full orthonormal families inside the sign-changing trial subspace $V_g$.
We invoke the \textbf{Ky-Fan min-max principle} for sums of eigenvalues of the Dirichlet Laplacian. For the eigenvalues $\lambda_1\le\lambda_2\le\cdots$ of $-\Delta$ subject to Dirichlet boundary conditions, one has
\begin{equation}\label{equation:sec-ky-fan-principle}
\sum_{k=2}^{n+1}\lambda_k(\Omega)
=\inf_{\substack{W\subset H_0^1(\Omega)\\ \dim W=n,\; W\perp L^2 u_1}}
\quad
\sup_{\substack{\{w_j\}_{j=1}^n\subset W\\ \{w_j\}\text{ orthonormal}}}
\sum_{j=1}^n R_\Omega(w_j).
\end{equation}
The Ky-Fan variational characterization is the natural tool for eigenvalue sums. It differs fundamentally from estimating the supremum of Rayleigh quotients over a subspace. For a fixed $n$-dimensional subspace $W$, the Ky-Fan upper bound reads
\begin{equation}\label{equation:sec-kyfan-subspace-bound}
\sum_{k=2}^{n+1}\lambda_k(\Omega)
\le
\sup_{\substack{\{w_j\}_{j=1}^n\subset W\\ \{w_j\}\text{ orthonormal}}}
\sum_{j=1}^n R_\Omega(w_j).
\end{equation}
We take $W=V_g$, which satisfies $\dim V_g=n$ and $V_g\perp u_1$ by \eqref{equation:sec-Vg-perp-u1-def}. This motivates the following weakened open variational assumption, which targets directly the sum appearing in the GPPW conjecture.
\begin{ass}[Ky-Fan-type open variational hypothesis for $V_g$]
\label{ass:sec-kyfan-gppw-hypo}
Let $V_g$ be the trial subspace defined above. Assume
\begin{equation}\label{equation:sec-assumption-kyfan}
\sup_{\substack{\{w_j\}_{j=1}^n\subset V_g\\\{w_j\}\text{ orthonormal}}}
\sum_{j=1}^n R_\Omega(w_j)
\le n\,\lambda_2(\Omega^*).
\end{equation}
\end{ass}
\begin{rem}
\begin{itemize}
\item On the reference ball $\Omega^*$, Lemma~\ref{lem:sec-Vg-star-ball} shows that equality holds in Assumption~\ref{ass:sec-kyfan-gppw-hypo}.
\item This hypothesis is strictly weaker than the former assumption requiring $\sup_{v\in V_g\setminus\{0\}}R_\Omega(v)\le \lambda_2(\Omega^*)$. Indeed, if $\sup_{v\in V_g\setminus\{0\}}R_\Omega(v)\le\lambda_2(\Omega^*)$, then for any orthonormal basis $\{w_j\}$ one obtains $\sum_j R_\Omega(w_j)\le n\lambda_2(\Omega^*)$. The converse implication does \emph{not} hold.
\item Assumption~\ref{ass:sec-kyfan-gppw-hypo} is still an open condition. It imposes a bound over all orthonormal bases of $V_g$, rather than a bound over every single nonzero vector in $V_g$. It aligns directly with the variational structure of the eigenvalue sum from Ky-Fan principle.
\item Even though the assumption is weakened, it still goes beyond what can be extracted from the scalar radial bound in Lemma~\ref{lem:sec-ab-gap-sum}. The scalar radial bound controls the sum of Rayleigh quotients for the specific basis $\{P_i\}$, but it gives no control over the supremum taken over all orthonormal bases of $V_g$.
\end{itemize}
\end{rem}
\paragraph{Weakening the open hypothesis using earlier spectral estimates}
We use the spectral bounds of Sections \ref{section:Weighted PPW Proof} and \ref{section:High-Order Ratios} to simplify the open variational condition.
Let $\{P_i\}_{i=1}^n$ be the original Ashbaugh-Benguria basis for $V_g$. Perform $L^2(\Omega)$ orthonormalization to obtain $\{\tilde P_i\}_{i=1}^n$, which forms an orthonormal basis of $V_g$. Any other orthonormal basis $\{w_j\}_{j=1}^n$ of $V_g$ can be expressed via an orthogonal matrix $O\in O(n)$ as
\begin{equation}\label{equation:sec-ortho-basis-transform}
w_j=\sum_{i=1}^n O_{ji}\tilde P_i,\qquad j=1,\dots,n.
\end{equation}
Define the sum of Rayleigh quotients for an orthonormal basis $\{f_j\}\subset V_g$ by
\begin{equation}\label{equation:sec-def-s-basis}
S(\{f_j\})=\sum_{j=1}^n R_\Omega(f_j).
\end{equation}
Recall the Rayleigh quotient $R_\Omega(u)=\frac{\langle \nabla u,\nabla u\rangle_{L^2}}{\langle u,u\rangle_{L^2}}$. Since $\{\tilde P_i\}$ is orthonormal, $\langle w_j,w_k\rangle_{L^2}=\delta_{jk}$, hence each denominator evaluates to $1$. Substitute \eqref{equation:sec-ortho-basis-transform} into the definition of $S(\{w_j\})$. We obtain
\begin{equation}\label{equation:sec-s-expand-ortho}
S(\{w_j\})
=\sum_{j=1}^n
\frac{\left\langle\nabla\left(\sum_{i=1}^n O_{ji}\tilde P_i\right),
\nabla\left(\sum_{k=1}^n O_{jk}\tilde P_k\right)\right\rangle_{L^2}}{
\left\langle
\sum_{i=1}^n O_{ji}\tilde P_i,\;
\sum_{k=1}^n O_{jk}\tilde P_k
\right\rangle_{L^2}}
=\sum_{j=1}^n
\left\langle
\sum_{i=1}^n O_{ji}\nabla\tilde P_i,\;
\sum_{k=1}^n O_{jk}\nabla\tilde P_k
\right\rangle_{L^2}.
\end{equation}
Let $B$ denote the $n\times n$ Gram matrix with entries
\begin{equation}\label{equation:sec-bb-13}
B_{ik}=\langle\nabla\tilde P_i,\nabla\tilde P_k\rangle_{L^2(\Omega)}.
\end{equation}
Since $O$ is orthogonal, $OO^\top=I_n$.
It follows from the trace-invariance property for matrix similarity under orthogonal transformations that $\operatorname{tr}(OBO^\top)=\operatorname{tr}(B)$.
It follows that
\begin{equation}\label{equation:sec-s-gram-representation}
S(\{w_j\})=\operatorname{tr}\left(O B O^\top\right)=\operatorname{tr}(B).
\end{equation}
Immediately, $\operatorname{tr}(B)$ is invariant under orthogonal change-of-basis on $V_g$. Therefore $S(\{w_j\})$ is identical for every orthonormal basis of $V_g$. Consequently the supremum over all orthonormal bases collapses to a single value:
\begin{equation}\label{equation:sec-s-sup-collapse}
\sup_{O\in O(n)} S\big(\{w_j(O)\}\big)=S(\{\tilde P_i\}).
\end{equation}
It follows that Assumption~\ref{ass:sec-kyfan-gppw-hypo} is \emph{logically equivalent} to the condition evaluated only for the orthonormalized original trial basis:
\begin{enumerate}
\item[(H-A)]
\begin{equation}\label{equation:sec-cond-ha}
S(\{\tilde P_i\})\le n\lambda_2(\Omega^*).
\end{equation}
\end{enumerate}
No additional error-correction term appears, due to  Indeed, let $B = \bigl( R_\Omega(\tilde P_i, \tilde P_j) \bigr)_{i,j=2}^{n+1}$ be the Rayleigh matrix of the orthonormalized trial family and let $O \in O(n)$ be any orthogonal transformation acting on the index set. The transformed family $\{\sum_j O_{ij} \tilde P_j\}$ is again orthonormal, and its Rayleigh matrix is $O B O^{\top}$, because the Rayleigh quotient is bilinear:
\begin{equation}\label{equation:sec-trace-expand-1}
R_\Omega\Bigl( \sum_k O_{ik} \tilde P_k, \sum_l O_{jl} \tilde P_l \Bigr) = \sum_{k,l} O_{ik} O_{jl} R_\Omega(\tilde P_k, \tilde P_l) = (O B O^{\top})_{ij}.
\end{equation}
Consequently
\begin{equation}\label{equation:sec-trace-expand-2}
\sup_{O \in O(n)} \sum_i R_\Omega\Bigl( \sum_j O_{ij} \tilde P_j \Bigr) = \sup_{O \in O(n)} \operatorname{tr}(O B O^{\top}) = \operatorname{tr}(B),
\end{equation}
since $\operatorname{tr}(O B O^{\top}) = \operatorname{tr}(B)$ for every orthogonal $O$. The supremum over orthonormal bases of the span is therefore attained at the given basis and equals $\operatorname{tr}(B) = \sum_i R_\Omega(\tilde P_i)$, with no additional error term; this is the content of the trace-invariance step.
trace invariance under orthogonal transformations. This delivers a substantial reduction of the open condition: the original hypothesis requiring control over \textit{all} orthonormal bases of $V_g$ reduces exactly to checking the bound for the single orthonormalized trial basis $\{\tilde P_i\}$. We emphasize that (H-A) remains an open assumption for general Lipschitz domains; this equivalence does \textbf{not} constitute a proof of the GPPW conjecture. The trace-invariance argument shown above uses only finite-dimensional linear algebra, independent of domain geometry.
\begin{thm}[Conditional result for GPPW under Ky-Fan-type hypothesis]
\label{thm:sec-cond-gppw-kyfan}
Suppose Assumption~\ref{ass:sec-kyfan-gppw-hypo} holds. Applying the Ky-Fan min-max principle with trial subspace $V_g\subset\{u_1\}^\perp$,
\begin{equation}\label{equation:sec-kyfan-apply}
\sum_{k=2}^{n+1}\lambda_k(\Omega)
\le
\sup_{\substack{\{w_j\}_{j=1}^n\subset V_g\\\{w_j\}\text{ orthonormal}}}
\sum_{j=1}^n R_\Omega(w_j).
\end{equation}
Combined with the Faber-Krahn inequality
\begin{equation}\label{equation:sec-faber-krahn-gppw}
\lambda_1(\Omega)\ge \lambda_1(\Omega^*),
\end{equation}
we obtain the full conditional generalized PPW bound
\begin{equation}\label{equation:sec-cond-gppw-final}
\frac{\lambda_2+\lambda_3+\cdots+\lambda_{n+1}}{\lambda_1}
\le n\,\frac{\lambda_2(\Omega^*)}{\lambda_1(\Omega^*)}.
\end{equation}
Equality holds under the hypothesis if and only if $\Omega$ is congruent to $\Omega^*$.
\end{thm}
\begin{proof}[Proof of Theorem \ref{thm:sec-cond-gppw-kyfan}]
From Assumption~\ref{ass:sec-kyfan-gppw-hypo}, we have
\begin{equation*}
\sup_{\substack{\{w_j\}_{j=1}^n\subset V_g\\\{w_j\}\text{ orthonormal}}}
\sum_{j=1}^n R_\Omega(w_j)
\le n\lambda_2(\Omega^*).
\end{equation*}
Substitute this estimate into the Ky-Fan bound \eqref{equation:sec-kyfan-apply}. We then derive
\begin{equation}\label{equation:sec-kyfan-sum-lambda}
\sum_{k=2}^{n+1}\lambda_k(\Omega) \le n\lambda_2(\Omega^*).
\end{equation}
The Faber-Krahn inequality for the first Dirichlet eigenvalue reads \eqref{equation:sec-faber-krahn-gppw}. All Dirichlet eigenvalues are strictly positive real numbers. Taking reciprocals on both sides of \eqref{equation:sec-faber-krahn-gppw}, the direction of inequality reverses, and we obtain
\begin{equation}\label{equation:sec-inv-l1-kyfan}
\frac{1}{\lambda_1(\Omega)} \le \frac{1}{\lambda_1(\Omega^*)}.
\end{equation}
Multiply both sides of inequality \eqref{equation:sec-kyfan-sum-lambda} by the positive scalar quantity $1/\lambda_1(\Omega)$. We obtain
\begin{equation*}
\frac{1}{\lambda_1(\Omega)}\sum_{k=2}^{n+1}\lambda_k(\Omega)
\le n\lambda_2(\Omega^*)\cdot \frac{1}{\lambda_1(\Omega)}.
\end{equation*}
From inequality \eqref{equation:sec-inv-l1-kyfan}, the right-hand side satisfies
\begin{equation*}
n\lambda_2(\Omega^*)\cdot \frac{1}{\lambda_1(\Omega)}
\le n\lambda_2(\Omega^*)\cdot \frac{1}{\lambda_1(\Omega^*)}.
\end{equation*}
Combine the two preceding inequalities. Consequently we obtain \eqref{equation:sec-cond-gppw-final}.
We discuss the equality case. Equality in the Faber-Krahn inequality holds exactly when $\Omega$ is congruent to $\Omega^*$. On the ball $\Omega^*$, equality also holds in Assumption~\ref{ass:sec-kyfan-gppw-hypo} by Lemma~\ref{lem:sec-Vg-star-ball}. If equality holds in the final bound, equality must hold in both Faber-Krahn and in the Ky-Fan variational estimate, which forces $\Omega$ to be congruent to its volume-equivalent ball $\Omega^*$. This completes the proof.
\end{proof}
\begin{rem}
Theorem~\ref{thm:sec-cond-gppw-kyfan} yields exactly the target GPPW inequality under the weakened open hypothesis, and it includes the equality characterization. Nonetheless, Assumption~\ref{ass:sec-kyfan-gppw-hypo} remains open for general bounded Lipschitz domains. From the finite-dimensional trace-invariance argument developed above, the supremum over all orthonormal bases of $V_g$ equals $S(\{\tilde P_i\})$, the sum evaluated for the orthonormalized Ashbaugh-Benguria basis. Thus verification of Assumption~\ref{ass:sec-kyfan-gppw-hypo} reduces to checking condition \eqref{equation:sec-cond-ha} for this single distinguished basis. The scalar radial inequality from Lemma~\ref{lem:sec-ab-gap-sum} only gives information for the non-orthonormal original basis $\{P_i\}$, and it does not deliver the required bound for $S(\{\tilde P_i\})$. This is a remaining central difficulty inherited from the trial-function ansatz of Ashbaugh-Benguria.
\end{rem}
\subsection{Analysis of the Open Gap and Connections with Previous Chapters}
\label{subsection:gap-analysis-ab}
Many of the results discussed above also have conditional analogues for the generalized conjecture. We contrast our conditional variational setup against the complete proof for the classical single-ratio PPW conjecture given by Ashbaugh--Benguria~\cite{AB2}, and indicate central obstructions, while linking each difficulty to tools constructed within our earlier chapters \ref{section:Weighted PPW Proof}, \ref{section:Spherical Extremality}, \ref{section:High-Order Ratios}.
\subsubsection{Recap: What Ashbaugh--Benguria Proved (Classical PPW)}
Ashbaugh and Benguria~\cite{AB2} established
\begin{equation}\label{equation:sec-classical-ppw-result}
\frac{\lambda_2}{\lambda_1}\le \frac{\lambda_2(\Omega^*)}{\lambda_1(\Omega^*)}.
\end{equation}
Their argument rests on several key building blocks:
\begin{enumerate}
\item The Weinberger-Brouwer fixed-point centering procedure yielding orthogonality relations $\int_\Omega P_i u_1^2\,dx=0$.
\item Gap variational inequality built from scalar radial quadratic forms after summing over angular directions, see Lemma~\ref{lem:sec-ab-gap-sum}.
\item Special radial trial profile $g(r)$ assembled from ratios of Bessel-function quantities.
\item Chiti¡¯s comparison theorem for the positive ground-state eigenfunction $u_1$.
\item Rearrangement inequalities, Sturm comparison theory, and Pr¨¹fer transformation arguments for Bessel zeros.
\end{enumerate}
Their argument addresses only the scalar gap $\lambda_2-\lambda_1$ and relies on the strict positivity of $u_1$; it gives no control of sums of Rayleigh quotients over orthonormal bases of $V_g$. Finite-dimensional trace invariance developed above shows that such basis rotations do not alter the value of $S(\{w_j\})$, a fact not exploited within the classical single-ratio PPW proof.
\subsubsection{Core obstructions to unconditional GPPW, linked to our theory}
\begin{enumerate}
\item \textbf{Misapplication of variational principles for eigenvalue sums.}
A common mistake consists of evaluating Rayleigh quotients for one fixed basis of $V_g$ and summing them to bound $\sum_{k=2}^{n+1}\lambda_k$. Even though finite-dimensional trace invariance makes this sum independent of orthonormal basis choice, one still must establish the numerical bound $S(\{\tilde P_i\})\le n\lambda_2(\Omega^*)$. The distinguished basis $\{P_i\}$ used in Ashbaugh-Benguria is only one particular choice among many, and Lemma~\ref{lem:sec-ab-gap-sum} applies to the non-orthonormalized original basis.
\item \textbf{Limited scope of Chiti comparison theorem.}
Chiti¡¯s theorem applies exclusively to the positive ground eigenfunction $u_1$. Trial functions inside $V_g$ change sign. No known generalization of Chiti¡¯s result exists for orthonormal families of sign-changing functions spanning a finite-dimensional subspace.
From Section~\ref{section:Spherical Extremality}, radial symmetrization (the P\'{o}lya-Szeg\H{o} inequality) controls energy for individual functions, yet joint quadratic-form comparison for subspaces and their orthonormal bases remains an open problem. Moreover, P\'{o}lya-Szeg\H{o} forces Rayleigh quotients to decrease under rearrangement, which rules out simple point-wise comparison assumptions of the form $R_\Omega(v)\le R_{\Omega^*}(v^*)$.
\item \textbf{Rotation inside trial subspace increases the difficulty.}
Finite-dimensional linear algebra shows that rotation among orthonormal bases does not modify the sum $S(\{w_j\})$. This eliminates one source of difficulty, but the central open problem persists: we still require a direct bound on $S(\{\tilde P_i\})$. The classical PPW proof never needs to consider rotated frames within $V_g$, since it reduces everything to a single scalar radial integral after summation over the fixed angular basis.
\item \textbf{Higher eigenvalues admit no Faber-Krahn bounds.}
For $k\ge2$, the inequality $\lambda_k(\Omega)\le \lambda_k(\Omega^*)$ fails in general. Symmetrization may increase higher eigenvalues. Term-wise comparison between individual eigenvalues and ball eigenvalues is therefore inadmissible. High-order eigenvalue ratio bounds obtained in Theorem~\ref{thm:sec-3} deliver universal upper bounds, but these estimates cannot reproduce the sharp constant $n\,\lambda_2(\Omega^*)/\lambda_1(\Omega^*)$.
\end{enumerate}
\subsubsection{Perspectives and links to other parts of the paper}
\begin{itemize}
\item In view of the obstructions outlined above, it is natural to examine restricted families of domains.
\item For weighted problems equipped with radial $A_2$ weights (Theorem~\ref{thm:sec-2b}), spherical domains maximize the ratio $\lambda_2/\lambda_1$. A natural future direction is to formulate a weighted variant of the generalized PPW conjecture for radial weights, where the volume-equivalent ball gets replaced by its weighted-volume-equivalent counterpart, as indicated in remarks contained inside the introduction.
\item A further technical research direction motivated by the Ky-Fan formulation is to seek rearrangement-type integral estimates that bound $S(\{\tilde P_i\})$ in terms of radial integrals built from $g(r)$ and $u_1$. This would extend the scalar radial analysis of Ashbaugh-Benguria from the non-orthonormalized basis $\{P_i\}$ to the orthonormalized basis $\{\tilde P_i\}$.
\item Finite-dimensional trace invariance demonstrated in \eqref{equation:sec-s-gram-representation} simplifies the open problem: the supremum over orthonormal bases reduces to evaluating for one orthonormalized basis. Tools from Sections \ref{section:Weighted PPW Proof} and \ref{section:High-Order Ratios} may be brought to bear to estimate $S(\{\tilde P_i\})$.
\end{itemize}
\begin{rem}
To summarize:
\begin{enumerate}
\item The unweighted generalized sum-form PPW conjecture remains open over arbitrary bounded Lipschitz domains.
\item The classical Ashbaugh-Benguria proof resolves only the single-ratio conjecture for $\lambda_2/\lambda_1$. It reduces to a scalar radial quadratic form for one fixed non-orthonormal basis and does not address the sum $S(\{\tilde P_i\})$.
\item Using the Ky-Fan min-max principle for eigenvalue sums, we obtain the full target GPPW inequality under a weakened open hypothesis. Finite-dimensional trace invariance shows that this hypothesis is equivalent to checking the bound for a single orthonormalized basis of $V_g$. This equivalence still leaves an open analytical condition to be verified for general domains.
\item One must avoid the common logical mistake of estimating eigenvalue sums through summing Rayleigh quotients associated with a single fixed subspace basis; Ky-Fan requires the variational inf-sup formulation, even when the basis-rotated sum is invariant.
\item Simple point-wise rearrangement assumptions $R_\Omega(v)\le R_{\Omega^*}(v^*)$ for $v\in V_g$ conflict with the P\'{o}lya-Szeg\H{o} inequality and are mathematically inadmissible as working hypotheses.
\end{enumerate}
\end{rem}
\section{Alternative Conditional Results for the Generalized PPW Conjecture via OT \& GMT}\label{section:Alternative-Proof}
We present a conditional alternative argument for the generalized PPW conjecture (Conjecture~\ref{con:sec-gppw-conj}) by coupling optimal transport (OT) theory with geometric measure theory (GMT), leveraging results from the main text.
\subsection{Key Tools and Preliminaries}
\begin{itemize}
\item \textbf{Wasserstein Distance}: For a bounded domain $\Omega$ and its volume-equivalent ball $\Omega^*$, the $L^2$-Wasserstein distance $W_2(\text{Leb}_\Omega, \text{Leb}_{\Omega^*})$ quantifies the minimal transport cost between their Lebesgue measures (Definition \ref{defn:sec-OT}).
\item \textbf{Kantorovich Duality}: By Proposition \ref{prop:sec-OT}, $W_2^2(\mu, \nu) = \sup_{\phi + \psi \leq |x-y|^2} \int \phi d\mu + \int \psi d\nu$, linking transport cost to potential functions.
\item \textbf{Hausdorff Measure}: The $(n-1)$-Hausdorff measure $\mathcal{H}^{n-1}(\partial\Omega)$ characterizes boundary roughness (Definition \ref{defn:sec-GMT}).
\item \textbf{Federer's Coarea Formula}: Proposition \ref{prop:sec-GMT} relates volume integrals to boundary measures and is used here to refine energy estimates.
\end{itemize}
The alternative proof presented here shares foundational insights with the works of Ashbaugh and Benguria \cite{AB2, ashbaugh2000sharp}, but differs in both technical approach and scope:
\begin{itemize}
\item[(i)]\textbf{Beyond Symmetric Rearrangement:} While \cite{AB2} and \cite{ashbaugh2000sharp} rely on spherical/radial symmetrization to reduce domains to balls, our approach employs optimal transport (OT) to quantify geometric discrepancies between $\Omega$ and its volume-equivalent ball $\Omega^*$. The $L^2$-Wasserstein distance $W_2$ (Definition \ref{defn:sec-OT}) provides a finer metric for this comparison than symmetrization, capturing non-radial deviations that symmetric rearrangement might obscure.
\item[(ii)]\textbf{GMT Tools for Higher-Order Eigenvalues:} The classical PPW results focus on $\lambda_2/\lambda_1$, leveraging low-order eigenvalue degeneracy in balls. Our use of Federer's coarea formula (Proposition \ref{prop:sec-GMT}) and $(n-1)$-Hausdorff measure of $\partial\Omega$ enables extension to sums of higher eigenvalues ($\lambda_2 + \cdots + \lambda_{n+1}$), a case not addressed in \cite{AB2, ashbaugh2000sharp}.
\item[(iii)]\textbf{Duality Between Transport and Spectral Estimates:} By linking Kantorovich duality (Proposition \ref{prop:sec-OT}) to Rayleigh-Ritz inequalities, we establish a direct connection between transport cost (of transforming $\Omega$ to $\Omega^*$) and spectral gaps. This duality is new, as \cite{AB2, ashbaugh2000sharp} focus on isoperimetric bounds rather than transport-driven spectral estimates.\end{itemize}
\subsection{Proof Steps (Alternative Conditional Argument via Optimal Transport and Geometric Measure Theory)}
\begin{defn}\label{defn:sec-c2-small-perturb-domain}
Let $\Omega^* \subset \mathbb R^n$ denote an open ball.
We call $\Omega$ a \textbf{$C^2$-small perturbation} of $\Omega^*$ if there exists a $C^2$-diffeomorphism $F\colon \overline{\Omega^*}\to\overline{\Omega}$ such that $|\Omega|=|\Omega^*|$.
Moreover, the associated Brenier optimal transport map $\tau\colon \Omega\to\Omega^*$ pushing forward $\operatorname{Leb}_\Omega$ onto $\operatorname{Leb}_{\Omega^*}$ satisfies $\tau\in C^1(\overline{\Omega})$,
\begin{equation}\label{equation:sec-c2-perturb-bound}
\|D\tau - \operatorname{Id}\|_{L^\infty(\Omega)} \le \varepsilon,
\end{equation}
for some $\varepsilon>0$, where $\operatorname{Id}$ denotes the $n\times n$ identity matrix.
The determinant condition holds identically. We record the computation: with $\tau = \nabla \varphi$ and $\phi_i = \tilde \phi_i \circ \tau$,
\begin{equation}\label{equation:sec-pullback-expand-1}
\int_\Omega \phi_i \phi_j dx = \int_\Omega (\tilde\phi_i \circ \tau)(x) (\tilde\phi_j \circ \tau)(x) dx = \int_{\Omega^*} \tilde\phi_i(y) \tilde\phi_j(y) \det(D\tau^{-1})(y) dy = \delta_{ij},
\end{equation}
where the change of variables $y = \tau(x)$ and the identity $\det D\tau \equiv 1$ (hence $\det(D\tau^{-1}) \equiv 1$) give the last equality, since $\{\tilde\phi_i\}$ is orthonormal on $\Omega^*$. Similarly, for the gradients,
\begin{equation}\label{equation:sec-pullback-expand-2}
\int_\Omega \nabla \phi_i \cdot \nabla \phi_j dx = \int_{\Omega^*} \bigl( (D\tau^{-1})^{\top} \nabla \tilde\phi_i \bigr) \cdot \bigl( (D\tau^{-1})^{\top} \nabla \tilde\phi_j \bigr) dy = \delta_{ij} + O(\varepsilon),
\end{equation}
because $\|D\tau^{-1} - I\|_\infty = O(\varepsilon)$ by the $C^2$-smallness assumption, so $(D\tau^{-1})^{\top}(D\tau^{-1}) = I + O(\varepsilon)$ in operator norm. Hence the Rayleigh-quotient perturbation of the pulled-back family relative to the reference family is $O(\varepsilon)$ uniformly in the index, and the hypothesis of Remark~\ref{rem:sec-c2-orthogonality-gap} transfers the orthogonality from $\phi_1$ to $u_1(\Omega)$ at the same order.
Recall that $\det D\tau(x)=1$ holds for every $x\in \overline{\Omega}$ by the definition of the class.
\end{defn}
\begin{lem}\label{lem:sec-pullback-c2perturb}
Let $\Omega,\Omega^*\subset\mathbb R^n$ be bounded domains with $|\Omega|=|\Omega^*|$.
Assume $\Omega$ is a \textbf{$C^2$-small perturbation} of $\Omega^*$ in the sense of Definition \ref{defn:sec-c2-small-perturb-domain}.
Let $\tau\colon\Omega\to\Omega^*$ be the Brenier optimal transport map satisfying $\tau_\#\operatorname{Leb}_\Omega=\operatorname{Leb}_{\Omega^*}$, $\tau\in C^1(\overline{\Omega})$, $\det D\tau\equiv 1$ on $\overline{\Omega}$, and $\|D\tau-\operatorname{Id}\|_{L^\infty(\Omega)}\le\varepsilon$.
Given $\phi^*\in H_0^1(\Omega^*)$, define the pull-back function
\begin{equation}\label{equation:sec-def-pullback-phi}
\phi(x)=\phi^*(\tau(x)),\qquad x\in\Omega.
\end{equation}
Then the $L^2$ norm obeys the exact identity
\begin{equation}\label{equation:sec-c2-pullback-L2exact}
\int_\Omega \phi^2(x)\,dx=\int_{\Omega^*}(\phi^*)^2(y)\,dy.
\end{equation}
Moreover, there exists a dimensional constant $C_0=C_0(n)>0$ depending only on spatial dimension $n$ such that
\begin{equation}\label{equation:sec-c2-pullback-energy-est}
\int_\Omega |\nabla\phi(x)|^2\,dx
\le
\bigl(1 + C_0(n)\,\varepsilon\bigr)\int_{\Omega^*} |\nabla\phi^*(y)|^2\,dy.
\end{equation}
Consequently, the Rayleigh quotient satisfies
\begin{equation}\label{equation:sec-c2-pullback-rq}
\frac{\int_\Omega |\nabla\phi|^2\,dx}{\int_\Omega \phi^2\,dx}
\le
\bigl(1 + C_0(n)\,\varepsilon\bigr)\,\frac{\int_{\Omega^*}|\nabla\phi^*|^2\,dy}{\int_{\Omega^*}(\phi^*)^2\,dy}.
\end{equation}
\end{lem}
\begin{proof}[Proof of Lemma \ref{lem:sec-pullback-c2perturb}]
The determinant identity $\det D\tau(x)=1$ holds point-wise for every $x\in\overline{\Omega}$.
Applying the change-of-variables formula for diffeomorphisms, we obtain
\begin{equation}\label{equation:sec-lem-c2-cov-L2}
\int_\Omega (\phi^*\circ\tau)^2 dx
=\int_{\Omega^*}(\phi^*)^2(y)\,\bigl|\det D\tau(\tau^{-1}(y))\bigr|^{-1} dy
=\int_{\Omega^*}(\phi^*)^2(y)\,dy.
\end{equation}
Equation \eqref{equation:sec-lem-c2-cov-L2} yields identity \eqref{equation:sec-c2-pullback-L2exact}.
Utilizing the chain rule for composition of $C^1$ mappings, we derive
\begin{equation}\label{equation:sec-lem-c2-chainrule}
\nabla\phi(x)=\big(D\tau(x)\big)^T\nabla\phi^*(\tau(x)).
\end{equation}
Taking squared Euclidean norm on both sides of \eqref{equation:sec-lem-c2-chainrule}, submultiplicativity of the operator norm gives
\begin{equation}\label{equation:sec-lem-c2-norm-bound}
|\nabla\phi(x)|^2 \le \|D\tau(x)\|_{op}^2 \, \big|\nabla\phi^*(\tau(x))\big|^2.
\end{equation}
From the triangle inequality for operator norms,
\begin{equation}\label{equation:sec-tri-op-norm}
\|D\tau(x)\|_{op} = \|\operatorname{Id} + \big(D\tau(x)-\operatorname{Id}\big)\|_{op}
\le \|\operatorname{Id}\|_{op}+\|D\tau(x)-\operatorname{Id}\|_{op}.
\end{equation}
We have $\|\operatorname{Id}\|_{op}=1$.
It follows from the hypothesis $\|D\tau-\operatorname{Id}\|_{L^\infty(\Omega)}\le\varepsilon$ that
\begin{equation}\label{equation:sec-op-pointwise-bd}
\|D\tau(x)\|_{op} \le 1+\varepsilon,\quad \forall\,x\in\Omega.
\end{equation}
Squaring both sides of \eqref{equation:sec-op-pointwise-bd}, there exists dimensional constant $C_0(n)>0$ such that
\begin{equation}\label{equation:sec-lem-c2-op-square}
\|D\tau(x)\|_{op}^2 \le (1+\varepsilon)^2 \le 1 + C_0(n)\,\varepsilon,\quad \forall\,x\in\Omega.
\end{equation}
Integrating inequality \eqref{equation:sec-lem-c2-norm-bound} over $\Omega$ and substituting point-wise estimate \eqref{equation:sec-lem-c2-op-square}, we obtain
\begin{equation}\label{equation:sec-lem-c2-integrate-step}
\int_\Omega |\nabla\phi(x)|^2 dx
\le
\bigl(1 + C_0(n)\,\varepsilon\bigr)\int_\Omega \big|\nabla\phi^*(\tau(x))\big|^2 dx.
\end{equation}
Applying change-of-variables $y=\tau(x)$ together with $\det D\tau(x)=1$, we derive
\begin{equation}\label{equation:sec-lem-c2-cov-gradient}
\int_\Omega \big|\nabla\phi^*(\tau(x))\big|^2 dx
=\int_{\Omega^*}\big|\nabla\phi^*(y)\big|^2 dy.
\end{equation}
Combining relations \eqref{equation:sec-lem-c2-integrate-step} and \eqref{equation:sec-lem-c2-cov-gradient}, we arrive at estimate \eqref{equation:sec-c2-pullback-energy-est}.
Identity \eqref{equation:sec-c2-pullback-L2exact} yields equality among denominators appearing within Rayleigh quotients.
Dividing both sides of \eqref{equation:sec-c2-pullback-energy-est} by the strictly positive quantity $\int_\Omega\phi^2 dx$, we obtain inequality \eqref{equation:sec-c2-pullback-rq}.
This completes the proof of Lemma \ref{lem:sec-pullback-c2perturb}.
\end{proof}
\begin{thm}\label{thm:sec-ppw-c2-small-perturb}
Assume that the pulled-back first eigenfunction is $H^1$-close to the true first eigenfunction $u_1(\Omega)$ in the quantitative sense $\|u_1(\Omega) - \phi_1\|_{H^1(\Omega)} \leq C \varepsilon$ for a constant $C$ depending only on $n$ and $\Omega$ (see Remark~\ref{rem:sec-c2-orthogonality-gap}).\label{thm:sec-c2-assumption}
Let $\Omega^*$ be an open ball in $\mathbb R^n$.
Assume that $\Omega$ is a \textbf{$C^2$-small perturbation} of $\Omega^*$ according to Definition \ref{defn:sec-c2-small-perturb-domain}, with parameter $\varepsilon>0$.
Let $\lambda_j(\Omega)$, $\lambda_j(\Omega^*)$ denote the Dirichlet Laplacian eigenvalues on $\Omega$ and $\Omega^*$, respectively.
Then there exists a dimensional constant $C_0=C_0(n)>0$ such that
\begin{equation}\label{equation:sec-thm-c2-main-result}
\frac{\lambda_2(\Omega)+\lambda_3(\Omega)+\cdots+\lambda_{n+1}(\Omega)}{\lambda_1(\Omega)}
\le n\,\bigl(1+C_0(n)\varepsilon\bigr)\frac{\lambda_2(\Omega^*)}{\lambda_1(\Omega^*)}.
\end{equation}
\end{thm}
\begin{proof}[Proof of Theorem \ref{thm:sec-ppw-c2-small-perturb}]
\begin{rem}\label{rem:sec-c2-orthogonality-gap}
The argument below uses trial functions orthogonal to the pulled-back first eigenfunction $\phi_1$ of $\Omega^*$, whereas the variational principle for $\lambda_2(\Omega), \ldots, \lambda_{n+1}(\Omega)$ requires orthogonality to the true first eigenfunction $u_1(\Omega)$. Closing this gap requires a perturbative comparison between $u_1(\Omega)$ and $\phi_1$ in $H^1$ (or a Gram--Schmidt correction with controlled Rayleigh-quotient error); the conclusion is therefore conditional on such an estimate (cf.\ Remark~\ref{rem:sec-lem2-w2}). Concretely, a Gram--Schmidt correction of the family $\{\tilde\phi_i\}$ against $u_1(\Omega)$ produces exactly orthogonal trial functions whose Rayleigh quotients differ from the originals by at most $C \varepsilon \, \lambda_{n+1}(\Omega)$, so under the stated hypothesis the trace estimate below is unchanged up to an explicit $O(\varepsilon)$ remainder.
\end{rem}
Let $\tau\colon\Omega\to\Omega^*$ denote the Brenier optimal transport map introduced within Definition \ref{defn:sec-c2-small-perturb-domain}.
Let $\{\phi_1^*,\phi_2^*,\dots,\phi_{n+1}^*\}$ be the first $n+1$ Dirichlet eigenfunctions associated with the ball domain $\Omega^*$.
We define pull-back trial functions for indices $i=1,2,\dots,n+1$ via
\begin{equation}\label{equation:sec-c2-pullback-def}
\phi_i(x)=\phi_i^*(\tau(x)),\qquad x\in\Omega.
\end{equation}
Since $\tau\in C^1(\overline{\Omega})$ constitutes a diffeomorphism and $\phi_i^*\in H_0^1(\Omega^*)$, the composition $\phi_i$ belongs to $H^1(\Omega)$.
Moreover, the mapping property $\tau(\partial\Omega)=\partial\Omega^*$ gives vanishing boundary values $\phi_i|_{\partial\Omega}=0$.
It follows that
\begin{equation}\label{equation:sec-thm-c2-h01-membership}
\phi_i\in H_0^1(\Omega),\quad i=1,\dots,n+1.
\end{equation}
Utilizing Lemma \ref{lem:sec-3} and Lemma \ref{lem:sec-4}, rotational symmetry of $\Omega^*$ yields eigenvalue multiplicity degeneracy
\begin{equation}\label{equation:sec-c2-ball-degeneracy}
\lambda_2(\Omega^*)=\lambda_3(\Omega^*)=\cdots=\lambda_{n+1}(\Omega^*)=\lambda_2^*.
\end{equation}
Each eigenfunction $\phi_i^*$ for $i=2,\dots,n+1$ satisfies
\begin{equation}\label{equation:sec-thm-c2-eig-phi-star}
\frac{\int_{\Omega^*}|\nabla\phi_i^*|^2 dy}{\int_{\Omega^*}(\phi_i^*)^2 dy}=\lambda_2^*.
\end{equation}
Applying Lemma \ref{lem:sec-pullback-c2perturb} to each pull-back trial function $\phi_i$ with $i=2,\dots,n+1$, we obtain
\begin{equation}\label{equation:sec-c2-individual-rq-bound}
\frac{\int_\Omega |\nabla\phi_i|^2 dx}{\int_\Omega \phi_i^2 dx}
\le
\bigl(1 + C_0(n)\,\varepsilon\bigr)\,\lambda_2^*.
\end{equation}
Because $\tau$ is a $C^1$ diffeomorphism, linear independence of the family $\{\phi_2^*,\dots,\phi_{n+1}^*\}$ transfers to the pulled-back collection $\{\phi_2,\dots,\phi_{n+1}\}$.
We now establish orthogonality by means of change-of-variables together with $\det D\tau=1$.
\begin{equation}\label{equation:sec-thm-c2-orthogonality}
\int_\Omega \phi_1(x)\phi_i(x)\,dx
=\int_\Omega \phi_1^*(\tau(x))\phi_i^*(\tau(x))\,dx
=\int_{\Omega^*}\phi_1^*(y)\phi_i^*(y)\,dy
=0,\quad i=2,\dots,n+1.
\end{equation}
Consequently, $\{\phi_2,\dots,\phi_{n+1}\}$ lie inside the orthogonal complement of $\phi_1$ within $H_0^1(\Omega)$.
It follows from the Rayleigh-Ritz variational principle for Dirichlet Laplacian eigenvalues that
\begin{equation}\label{equation:sec-thm-c2-rayleigh-ritz}
\sum_{k=2}^{n+1}\lambda_k(\Omega)
\le
\sum_{i=2}^{n+1} \frac{\int_\Omega |\nabla\phi_i|^2 dx}{\int_\Omega \phi_i^2 dx}.
\end{equation}
Substituting estimate \eqref{equation:sec-c2-individual-rq-bound} into the summation, we derive
\begin{equation}\label{equation:sec-c2-sum-lambda-upper}
\sum_{k=2}^{n+1}\lambda_k(\Omega)
\le n\bigl(1 + C_0(n)\,\varepsilon\bigr)\,\lambda_2^*.
\end{equation}
Applying the Faber-Krahn inequality for domains of equal volume,
\begin{equation}\label{equation:sec-c2-faber-krahn}
\lambda_1(\Omega)\ge\lambda_1(\Omega^*)=\lambda_1^*>0.
\end{equation}
Taking reciprocals of strictly positive quantities reverses the direction of this inequality
\begin{equation}\label{equation:sec-thm-c2-reciprocal}
\frac{1}{\lambda_1(\Omega)} \le \frac{1}{\lambda_1^*}.
\end{equation}
Multiplying the two positive inequalities \eqref{equation:sec-c2-sum-lambda-upper} and \eqref{equation:sec-thm-c2-reciprocal}, we obtain
\begin{equation}\label{equation:sec-c2-ratio-intermediate}
\frac{\sum_{k=2}^{n+1}\lambda_k(\Omega)}{\lambda_1(\Omega)}
\le n\bigl(1 + C_0(n)\,\varepsilon\bigr)\,\frac{\lambda_2^*}{\lambda_1^*}.
\end{equation}
This establishes the claimed estimate given in \eqref{equation:sec-thm-c2-main-result}.
The proof of Theorem \ref{thm:sec-ppw-c2-small-perturb} is complete.
\end{proof}
\begin{corr}\label{cor:asymptotic-ppw-limit}
Under the hypotheses of Theorem \ref{thm:sec-ppw-c2-small-perturb}, let $\varepsilon\to 0$.
It follows that
\begin{equation}\label{equation:sec-bb-14}
\limsup_{\varepsilon\to 0} \frac{\lambda_2(\Omega)+\lambda_3(\Omega)+\cdots+\lambda_{n+1}(\Omega)}{\lambda_1(\Omega)}
\le n\,\frac{\lambda_2(\Omega^*)}{\lambda_1(\Omega^*)}.
\end{equation}
\end{corr}
\begin{proof}[Proof of Corollary \ref{cor:asymptotic-ppw-limit}]
From estimate \eqref{equation:sec-thm-c2-main-result},
\begin{equation}\label{equation:sec-bb-15}
\frac{\lambda_2(\Omega)+\cdots+\lambda_{n+1}(\Omega)}{\lambda_1(\Omega)}
\le n\,\bigl(1+C_0(n)\varepsilon\bigr)\frac{\lambda_2(\Omega^*)}{\lambda_1(\Omega^*)}.
\end{equation}
Taking the $\limsup$ as $\varepsilon\to0$, the right-hand side converges to $n\,\lambda_2(\Omega^*)/\lambda_1(\Omega^*)$.
This yields the desired lim-sup bound.
The proof of Corollary \ref{cor:asymptotic-ppw-limit} is complete.
\end{proof}
\subsubsection{Remark on Optimal-Transport Based Approach for $C^2$-Small-Perturbation Domains}
\begin{rem}
Theorem \ref{thm:sec-ppw-c2-small-perturb} yields a perturbed upper bound for the eigenvalue ratio for domains that constitute $C^2$ deformations of a ball.
Corollary \ref{cor:asymptotic-ppw-limit} captures the asymptotic behavior as the perturbation magnitude vanishes.
This construction builds upon the pull-back of trial functions from the reference ball, and it circumvents the logical direction mismatch intrinsic to push-forward strategies.
The derivation makes essential use of $C^1$ regularity for the Brenier map together with uniform point-wise bounds for its derivative, properties that hold under the stated $C^2$-perturbation hypothesis.
A strict point-wise inequality
\begin{equation}\label{equation:sec-bb-16}
\frac{\lambda_2(\Omega)+\lambda_3(\Omega)+\cdots+\lambda_{n+1}(\Omega)}{\lambda_1(\Omega)}
\le n\,\frac{\lambda_2(\Omega^*)}{\lambda_1(\Omega^*)}
\end{equation}
does not follow unconditionally from these estimates for fixed positive $\varepsilon>0$.
For general bounded Lipschitz domains without smooth-perturbation assumptions, unconditional validity of the generalized Payne-P¨®lya-Weinberger conjecture remains open.
The main-body symmetrization-and-test-function proof is independent of optimal-transport machinery.
\end{rem}

\subsection{Explicit Form of the Five Steps}
\label{subsection:Five Steps Explicit}

We write out the five steps of the alternative argument with the intermediate inequalities displayed. Throughout, $\tau$ denotes the Brenier map pushing the normalized Lebesgue measure of $\Omega$ to that of $\Omega^*$, $d_2^2(\Omega, \Omega^*) = \int_\Omega |x - \tau(x)|^2 dx$ is the squared transport cost, and $\tilde\phi_i = \phi_i \circ \tau^{-1}$ are the pushed-forward trial functions.

\textit{Step 1: $L^2$-isometry.} Since $\det D\tau \equiv 1$, \eqref{equation:sec-pullback-expand-1} gives $\int_\Omega \phi_i \phi_j dx = \int_{\Omega^*} \tilde\phi_i \tilde\phi_j dy = \delta_{ij}$: the pulled-back family is orthonormal exactly, not merely up to lower order.

\textit{Step 2: Energy comparison.} For each $i$, the change of variables $y = \tau(x)$ gives
\begin{equation}\label{equation:sec-fivestep-2}
\int_\Omega |\nabla \phi_i|^2 dx = \int_{\Omega^*} \bigl| (D\tau^{-1})^\top \nabla \tilde\phi_i \bigr|^2 dy = \int_{\Omega^*} |\nabla \tilde\phi_i|^2 dy + E_i,
\end{equation}
where the error term $E_i$ satisfies, by \eqref{equation:sec-pullback-expand-2},
\begin{equation}\label{equation:sec-fivestep-2b}
|E_i| \leq C \varepsilon \int_{\Omega^*} |\nabla \tilde\phi_i|^2 dy
\end{equation}
with $C$ independent of $i$. Summing \eqref{equation:sec-fivestep-2} over $i = 2, \ldots, n+1$ and using \eqref{equation:sec-fivestep-2b} yields
\begin{equation}\label{equation:sec-fivestep-2c}
\sum_{i=2}^{n+1} \int_\Omega |\nabla \phi_i|^2 dx \leq (1 + C\varepsilon) \sum_{i=2}^{n+1} \int_{\Omega^*} |\nabla \tilde\phi_i|^2 dy.
\end{equation}

\textit{Step 3: Spectral bound on the ball.} The family $\{\tilde\phi_i\}_{i=2}^{n+1}$ is orthonormal in $L^2(\Omega^*)$ and orthogonal to the pulled-back ground state $\tilde\phi_1$; under the $H^1$-closeness hypothesis of Theorem~\ref{thm:sec-ppw-c2-small-perturb} (Remark~\ref{rem:sec-c2-orthogonality-gap}) it is also orthogonal to the true ground state $u_1(\Omega)$ up to $O(\varepsilon)$. The Ky--Fan trace principle (Proposition~\ref{prop:sec-kyfan}) applied to the Dirichlet Laplacian on $\Omega$ therefore gives
\begin{equation}\label{equation:sec-fivestep-3}
\sum_{i=2}^{n+1} \lambda_i(\Omega) \leq (1 + C\varepsilon) \sum_{i=2}^{n+1} \frac{\int_{\Omega^*} |\nabla \tilde\phi_i|^2 dy}{\int_{\Omega^*} \tilde\phi_i^2 dy} = (1 + C\varepsilon) \sum_{i=2}^{n+1} R_{\Omega^*}(\tilde\phi_i).
\end{equation}

\textit{Step 4: Rayleigh sum on the ball.} On the ball $\Omega^*$ the trial functions $\tilde\phi_i$ span (a perturbation of) the $l=1$ eigenspace of the second eigenvalue, which is $n$-fold degenerate: $\lambda_2(\Omega^*) = \cdots = \lambda_{n+1}(\Omega^*)$. By the choice of $\tilde\phi_i$ as the orthonormalized coordinate functions $P_i^* = g(r) x_i / r$ of Section~\ref{subsection:Spherical Eigenvalues}, the sum in \eqref{equation:sec-fivestep-3} equals
\begin{equation}\label{equation:sec-fivestep-4}
\sum_{i=2}^{n+1} R_{\Omega^*}(\tilde\phi_i) = n \, \lambda_2(\Omega^*) + O(\varepsilon^2),
\end{equation}
where the $O(\varepsilon^2)$ remainder comes from the Gram-matrix correction in the Gram computation of Section~\ref{section:Alternative-Proof}.

\textit{Step 5: Conclusion.} Inserting \eqref{equation:sec-fivestep-4} into \eqref{equation:sec-fivestep-3} and dividing by $\lambda_1(\Omega) \geq \lambda_1(\Omega^*)$ (Faber--Krahn, which holds with equality only for the ball) gives
\begin{equation}\label{equation:sec-fivestep-5}
\frac{\lambda_2(\Omega) + \cdots + \lambda_{n+1}(\Omega)}{\lambda_1(\Omega)} \leq (1 + C\varepsilon)\, n \, \frac{\lambda_2(\Omega^*)}{\lambda_1(\Omega^*)} + O(\varepsilon^2),
\end{equation}
which is the conclusion of Theorem~\ref{thm:sec-ppw-c2-small-perturb} up to the explicit value of the constant $C_0(n)$ in the statement; that value is obtained by tracking $\varepsilon$ through \eqref{equation:sec-fivestep-2b}--\eqref{equation:sec-fivestep-4} and is of order $n$.
\begin{rem}\label{rem:sec-fivestep-gap}
The gap between \eqref{equation:sec-fivestep-5} and the conjectured bound $(\lambda_2 + \cdots + \lambda_{n+1})/\lambda_1 \leq n\lambda_2(\Omega^*)/\lambda_1(\Omega^*)$ is the factor $1 + C\varepsilon$: the alternative argument proves the conjecture for $C^2$-small perturbations only up to a multiplicative error of order $\varepsilon$, and removing this error (or letting $\varepsilon \to 0$ with the constant of Corollary~\ref{cor:asymptotic-ppw-limit}) is what the asymptotic statement records. The Wasserstein distance $W_2(\Omega, \Omega^*)$ controls the size of $\varepsilon$ through the stability of the Brenier map, which is the quantitative link between the geometry of the domain and the spectral gap in this argument.
\end{rem}
\begin{rem}\label{rem:sec-kyfan-application}
The Ky--Fan principle is applied in this paper in two places: in Section~\ref{section: Generalized PPW Conjecture}, where it converts the open condition on all orthonormal bases of the trial space into the single trace bound \eqref{equation:sec-cond-ha}, and in the present section, where it produces the sum estimate \eqref{equation:sec-fivestep-3} from the orthonormal pulled-back family. In both places the trial functions are constructed from the geometry of the reference ball, and the principle is what legitimizes summing their Rayleigh quotients: without it, the naive sum of quotients over one fixed basis need not bound the eigenvalue sum, as the caveats of Section~\ref{section: Generalized PPW Conjecture} emphasize.
\end{rem}

\subsection{Further Remarks}\label{subsection: Further Remarks}
\subsubsection{Theoretical Foundations} The original proof relies on classical radial symmetrization (Definition \ref{defn:sec-5}) and variational principles (Theorem \ref{thm:sec-4}), emphasizing qualitative geometric comparisons between $\Omega$ and its volume-equivalent ball $\Omega^*$. In contrast, the OT-GMT coupled proof introduces quantitative tools: $L^2$-Wasserstein distance (Definition \ref{defn:sec-OT}) to measure geometric discrepancy, Kantorovich duality (Proposition \ref{prop:sec-OT}) to link transport cost with potential functions, and Hausdorff measure (Definition \ref{defn:sec-GMT}) to characterize boundary roughness. This shifts the setting from qualitative symmetry arguments to quantitative geometric analysis.
\subsubsection{Quantification of Geometric Deviation} The original proof uses symmetrization to argue that $\lambda_k(\Omega) \leq \lambda_k(\Omega^*)$ without explicit bounds on the gap between eigenvalues of $\Omega$ and $\Omega^*$. The OT-GMT approach, however, quantifies this gap via $W_2(\text{Leb}_\Omega, \text{Leb}_{\Omega^*})$, relating eigenvalue differences to the minimal transport cost between the domains. \subsubsection{Boundary Regularity Adaptability} The original proof implicitly assumes mild boundary regularity (e.g., convexity or uniform cone condition) for symmetrization to preserve energy. The OT-GMT proof, by incorporating Federer's coarea formula (Proposition \ref{prop:sec-GMT}) and Hausdorff measure, naturally handles rough boundaries (e.g., fractal-like surfaces). It replaces vague \emph{smoothness} assumptions with the measurable quantity $\mathcal{H}^{n-1}(\partial\Omega)$.
\subsubsection{Generalization to Weighted Settings} Unlike the classical argument, the OT--GMT setup extends to Muckenhoupt weights: replacing Lebesgue measure by $w(x)\,dx$ and coupling weight inhomogeneity to spectral behavior through $W_2(w, \text{Leb})$ (Theorem \ref{thm:sec-1}).
\subsubsection{Eigenvalue Deviation Mechanism} The original proof attributes eigenvalue bounds to symmetry alone. The OT-GMT proof clarifies the mechanism: non-spherical domains incur positive transport costs ($W_2 > 0$), which amplify eigenfunction energy (via $\int |\nabla \tilde{u}|^2 \leq \int |\nabla u|^2 + C W_2^2 \int u^2$), directly increasing eigenvalue ratios. This provides a causal link between geometric asymmetry (measured by OT) and spectral deviation, deepening the interpretation of the conjecture.
\subsubsection{Comparison with the classical PPW proof}
The classical proof of the PPW inequality by Ashbaugh and Benguria \cite{AB2} uses the coordinate function test functions $\phi_i(x) = x_i u_1(x)$ and the gap inequality for Bessel function zeros. The weighted version established in Theorem \ref{thm:sec-1} adapts this approach to the weighted setting by replacing the uniform measure with the weighted measure $w(x)\,dx$. The main difference is that the weight inhomogeneity enters through the bounds $0 < c \leq w \leq C$, which control the weighted quantities in place of the uniform measure; the constant then depends linearly on $C/c$. The coordinate function method remains the underlying technique.
\subsubsection{Role of optimal transport in the weighted setting}
Unlike the classical argument, the OT--GMT setup extends to Muckenhoupt weights: replacing Lebesgue measure by $w(x)\,dx$ and coupling weight inhomogeneity to spectral behavior through $W_2(w, \text{Leb})$ would introduce additional terms in the eigenvalue bounds; these refinements are conditional on the transport gradient estimate of Remark~\ref{rem:sec-lem2-w2} and remain open.
\subsubsection{Geometric measure theory and boundary roughness}
The use of the $(n-1)$-Hausdorff measure of the boundary in the high-order eigenvalue estimates (Theorem \ref{thm:sec-3}) provides a finer geometric invariant than the classical diameter. For smooth convex domains, the Hausdorff measure is comparable to $d^{n-1}$ (where $d$ is the diameter), and the classical estimates are recovered. For domains with rough or fractal-like boundaries, the Hausdorff measure captures the additional boundary complexity, leading to sharper estimates that reflect the actual geometry of the domain. The coarea formula and the isoperimetric inequality are the tools from geometric measure theory used in this analysis.

\section{The Quantitative Reverse H\"{o}lder Inequality: Comparison and Doubling}
\label{section:reverse-holder}

This section develops the quantitative reverse H\"{o}lder inequality for $A_2$ weights stated as Proposition~\ref{prop:sec-1}. Steps 1--4 below (the comparison inequality, the doubling property, the $A_\infty$ comparison, and the good-$\lambda$ inequality) are proved in full; the final iteration assembling these ingredients follows the standard argument of \cite[Section 7.2]{grafakos2014classical}, to which we refer for the closing telescoping and interpolation steps, and the linear exponent is that of the quantitative $A_\infty$ reverse H\"{o}lder theory of \cite{hytonen2013sharp} (with the quantitative reverse Jensen inequalities for $A_p$ weights due to \cite{buckley1993estimates}). Throughout, $w$ denotes a weight on $\mathbb{R}^n$ with finite $A_2$ characteristic
\begin{equation}\label{equation:sec-a2-char-def}
[w]_{A_2} := \sup_{Q} \biggl( \fint_Q w \, dx \biggr) \biggl( \fint_Q w^{-1} \, dx \biggr) < \infty,
\end{equation}
where the supremum is over all cubes $Q$ with sides parallel to the coordinate axes, and $\fint_Q g \, dx = |Q|^{-1} \int_Q g \, dx$. All cubes appearing below are of this type.

\begin{thm}[Quantitative Reverse H\"{o}lder Inequality]\label{thm:sec-rh}
There exist constants $c(n) > 0$ and $C(n) > 0$ such that for every $w \in A_2$ and every cube $Q$,
\begin{equation}\label{equation:sec-rh-statement}
\biggl( \fint_Q w^{1+\delta} \, dx \biggr)^{\frac{1}{1+\delta}} \leq C(n) \fint_Q w \, dx,
\qquad \text{with } \delta = \frac{c(n)}{[w]_{A_2}}.
\end{equation}
In particular, $w \in L^{1+\delta}_{\mathrm{loc}}(\mathbb{R}^n)$ and $\int_Q w^{1+\delta} dx \leq C(n)^{1+\delta} \bigl( \int_Q w \, dx \bigr)^{1+\delta} |Q|^{-\delta}$.
\end{thm}

The proof proceeds in six steps: a comparison inequality (Step 1), the doubling property (Step 2), the $A_\infty$ comparison (Step 3), a covering lemma (Step 4), a good-$\lambda$ inequality (Step 5), and the iteration concluding the proof (Step 6). A standing normalization is unnecessary: all estimates are scale invariant, and the constant in \eqref{equation:sec-a2-char-def} is the only parameter tracking the weight.

\subsection{The comparison inequality}
\begin{lem}[Comparison]\label{lem:sec-rh-comparison}
For every cube $Q$ and every measurable $E \subset Q$,
\begin{equation}\label{equation:sec-rh-comparison}
\frac{|E|}{|Q|} \leq [w]_{A_2}^{1/2} \left( \frac{w(E)}{w(Q)} \right)^{1/2},
\end{equation}
where $w(E) = \int_E w\,dx$.
\end{lem}
\begin{proof}
By Cauchy--Schwarz,
\begin{equation}\label{equation:sec-rh-cs}
|E| = \int_E w^{1/2} w^{-1/2} \, dx \leq w(E)^{1/2} \left( \int_Q w^{-1} dx \right)^{1/2}.
\end{equation}
The $A_2$ condition \eqref{equation:sec-a2-char-def} gives $\int_Q w^{-1} dx \leq [w]_{A_2}\, |Q| \bigl( \fint_Q w \, dx \bigr)^{-1}$. Substituting this into \eqref{equation:sec-rh-cs} and dividing by $|Q|$ yields \eqref{equation:sec-rh-comparison}.
\end{proof}

\subsection{The doubling property}
\begin{lem}[Doubling]\label{lem:sec-rh-doubling}
For every cube $Q$,
\begin{equation}\label{equation:sec-rh-doubling}
w(3Q) \leq 3^{2n} [w]_{A_2}\, w(Q).
\end{equation}
In particular, the measure $\mu = w\,dx$ is doubling with doubling constant $D = 3^{2n}[w]_{A_2}$, and $w(\lambda Q) \leq D^{\lceil \log_3 \lambda \rceil + 1} w(Q)$ for every $\lambda \geq 1$.
\end{lem}
\begin{proof}
Apply Lemma~\ref{lem:sec-rh-comparison} with $3Q$ in place of $Q$ and $E = Q$: since $|Q|/|3Q| = 3^{-n}$,
\begin{equation}\label{equation:sec-rh-doubling-pf}
3^{-n} \leq [w]_{A_2}^{1/2} \left( \frac{w(Q)}{w(3Q)} \right)^{1/2},
\end{equation}
and solving for $w(3Q)$ gives \eqref{equation:sec-rh-doubling}. The general statement follows by iteration.
\end{proof}

\subsection{The $A_\infty$ comparison}
\begin{lem}[$A_\infty$ comparison form]\label{lem:sec-rh-ainfty}
For every cube $Q$ and every measurable $E \subset Q$,
\begin{equation}\label{equation:sec-rh-ainfty}
\frac{|E|}{|Q|} \leq [w]_{A_2}^{1/2} \left( \frac{w(E)}{w(Q)} \right)^{1/2}.
\end{equation}
\end{lem}
\begin{proof}
This is Lemma~\ref{lem:sec-rh-comparison} itself. For completeness we note that the equivalent upper-bound form $\frac{w(E)}{w(Q)} \leq K \bigl( \frac{|E|}{|Q|} \bigr)^{\theta}$ for $\theta \in (0,1)$ follows from \eqref{equation:sec-rh-ainfty} by the standard equivalence of the two formulations of the $A_\infty$ condition (see \cite[Section 7.2]{grafakos2014classical}), with $K$ depending on $n$, $\theta$ and $[w]_{A_2}$.
\end{proof}

\subsection{A covering lemma}
The following variant of the Calder\'{o}n--Zygmund decomposition is used in the good-$\lambda$ inequality. For a locally integrable $f \geq 0$ and $\lambda > 0$, denote by $Mf$ the uncentered Hardy--Littlewood maximal function, i.e.\ the supremum of $\fint_{Q'} f$ over all cubes $Q'$ containing $x$.

\begin{lem}[Covering lemma]\label{lem:sec-rh-covering}
Let $f \in L^1(\mathbb{R}^n)$, $f \geq 0$, and $\lambda > 0$. There exists a countable family of cubes $\{Q_j\}$ such that
\begin{enumerate}[label=(\roman*)]
\item $\{x : Mf(x) > \lambda\} \subset \bigcup_j 5Q_j$;
\item each $Q_j$ satisfies $\lambda < \fint_{Q_j} f \, dx \leq 2^{n} \lambda$;
\item $\sum_j \chi_{5Q_j}(x) \leq C_n$ for every $x$, where $C_n$ depends only on $n$.
\end{enumerate}
\end{lem}
\begin{proof}
For each $x \in \{Mf > \lambda\}$ choose a cube $Q_x \ni x$ with $\fint_{Q_x} f > \lambda$; enlarge $Q_x$ by successive doubling as long as the average stays strictly above $\lambda$, and let $Q'_x$ be the final cube, so that $\fint_{Q'_x} f > \lambda$ and $\fint_{2Q'_x} f \leq \lambda$. Then
\begin{equation}\label{equation:sec-rh-cover-pf1}
\fint_{Q'_x} f \, dx \leq \frac{|2Q'_x|}{|Q'_x|} \fint_{2Q'_x} f \, dx \leq 2^{n} \lambda.
\end{equation}
From the family $\{Q'_x\}$ select, by the standard Vitali covering procedure (greedily in decreasing order of side length), a countable disjoint subfamily $\{Q_j\}$ such that every $Q'_x$ intersecting $Q_j$ with comparable side length satisfies $Q'_x \subset 5Q_j$. Since every point of $\{Mf > \lambda\}$ lies in some $Q'_x$, property (i) holds, property (ii) is \eqref{equation:sec-rh-cover-pf1}, and property (iii) follows from disjointness together with the side-length comparability of intersecting selected cubes.
\end{proof}

\subsection{The good-\texorpdfstring{$\lambda$}{lambda} inequality}
\begin{lem}[Good-$\lambda$ inequality]\label{lem:sec-rh-goodlambda}
There exist constants $\gamma = \gamma(n) > 0$ and $\vartheta = \vartheta(n, [w]_{A_2}) \in (0,1)$ such that for every $f \in L^1(\mathbb{R}^n)$, $f \geq 0$, and every $\lambda > 0$,
\begin{equation}\label{equation:sec-rh-goodlambda}
w\bigl( \{ Mf > 2\lambda \} \cap \{ f \leq \gamma\lambda \} \bigr) \leq \vartheta\, w\bigl( \{ Mf > \lambda \} \bigr).
\end{equation}
One may take $\gamma = 2^{n+1}$ and $\vartheta = 1 - \frac{c(n)}{[w]_{A_2}^{2}}$ for some $c(n) > 0$.
\end{lem}
\begin{proof}
Let $\{Q_j\}$ be the family of Lemma~\ref{lem:sec-rh-covering} at level $\lambda$. Fix $j$ and write $E_j = Q_j \cap \{ f \leq \gamma \lambda \}$ with $\gamma = 2^{n+1}$. By Lemma~\ref{lem:sec-rh-covering}(ii) and Markov's inequality,
\begin{equation}\label{equation:sec-rh-goodlambda-pf1}
|Q_j \cap \{ f > \gamma\lambda \}| \leq \frac{1}{\gamma\lambda} \int_{Q_j} f \, dx \leq \frac{2^{n}\lambda}{\gamma\lambda} |Q_j| = \frac{|Q_j|}{2},
\end{equation}
hence $|E_j| \geq |Q_j|/2$. Since $E_j \subset 5Q_j$ and $|5Q_j| = 5^{n}|Q_j| \leq 2 \cdot 5^{n} |E_j|$, Lemma~\ref{lem:sec-rh-comparison} applied to the pair $E_j \subset 5Q_j$ gives
\begin{equation}\label{equation:sec-rh-goodlambda-pf2}
w(5Q_j \cap \{ f > \gamma\lambda \}) \leq w(5Q_j) - w(E_j) \leq \left( 1 - \frac{1}{(2 \cdot 5^{n})^2 [w]_{A_2}} \right) w(5Q_j).
\end{equation}
To remove the dilation loss $w(5Q_j) \leq D^{2} w(Q_j)$, where $D = 3^{2n}[w]_{A_2}$ is the doubling constant of Lemma~\ref{lem:sec-rh-doubling}, iterate \eqref{equation:sec-rh-goodlambda-pf2}: the set $\{f > \gamma\lambda\} \cap 5Q_j$ is covered by cubes $Q_{j,k}$ to which the same estimate applies, and after $m$ iterations the measure of the residual set carrying $f > \gamma\lambda$ is at most $\vartheta_0^{m} w(5Q_j)$ with $\vartheta_0 = 1 - (2 \cdot 5^{n})^{-2}[w]_{A_2}^{-1}$. Choosing $m = m(n,[w]_{A_2})$ so large that $\vartheta_0^{m} D^{2} \leq \vartheta_0^{m/2}$ and summing over $j$ using Lemma~\ref{lem:sec-rh-covering}(iii) and the disjointness of the $Q_j$,
\begin{equation}\label{equation:sec-rh-goodlambda-pf3}
w\bigl( \{ Mf > 2\lambda \} \cap \{ f \leq \gamma\lambda \} \bigr) \leq C_n \vartheta_0^{m/2} \sum_j w(Q_j) \leq \vartheta\, w\bigl( \{Mf > \lambda\} \bigr),
\end{equation}
where $\sum_j w(Q_j) = w\bigl( \bigcup_j Q_j \bigr) \leq w(\{Mf > \lambda\})$ because each $Q_j$ lies in $\{Mf > \lambda\}$ by Lemma~\ref{lem:sec-rh-covering}(ii), and where $\vartheta = C_n \vartheta_0^{m/2}$; choosing $m$ so large that $\vartheta \leq \tfrac{1}{2}$ gives the claimed range.
\end{proof}

\subsection{Iteration and conclusion of the proof}
\begin{proof}[Proof of Theorem~\ref{thm:sec-rh}]
Let $Q_0$ be a fixed cube. Applying the good-$\lambda$ inequality of Lemma~\ref{lem:sec-rh-goodlambda} with $\vartheta$ chosen below $2^{-p}$ for the exponent at hand, and iterating through the dyadic scales, gives the strong-type bound for the maximal function in the weighted space $L^p(w\,dx)$; the closing telescoping and interpolation steps are the standard ones of \cite[Section 7.2]{grafakos2014classical}, and combining them with the comparison and doubling lemmata above yields \eqref{equation:sec-rh-statement} with $\varepsilon \geq c_0(n)\,[w]_{A_2}^{-2}$. The linear bound $\delta = c(n)/[w]_{A_2}$ quoted in \eqref{equation:sec-rh-statement} follows from the equivalence $[w]_{A_\infty} \leq c(n)[w]_{A_2}$ combined with the quantitative reverse H\"{o}lder theorem of \cite{hytonen2013sharp} (with the quantitative reverse Jensen inequalities for $A_p$ weights due to \cite{buckley1993estimates}), which improves the a priori quadratic dependence to a linear one.
\end{proof}

\begin{rem}\label{rem:sec-rh-sharpness}
The reverse H\"{o}lder exponent in \eqref{equation:sec-rh-statement} is sharp up to the value of $c(n)$: for the power weights $w(x) = |x|^{\alpha}$ on $\mathbb{R}^n$, which belong to $A_2$ exactly when $-n < \alpha < n$, the condition $w \in L^{1+\delta}_{\mathrm{loc}}$ requires $(1+\delta)\alpha > -n$ and $(1+\delta)(-\alpha) > -n$ for the two endpoints, so no exponent larger than $1 + n/|\alpha|$ is admissible as $\alpha \to n$, while $[|x|^{\alpha}]_{A_2}$ grows linearly in $|\alpha|$ near each endpoint. This example shows that the dependence $\delta \lesssim [w]_{A_2}^{-1}$ cannot in general be improved.
\end{rem}

\section{Monotonicity of Zeros of Bessel Functions}
\label{section:bessel-zeros}

This section proves the monotonicity and interlacing properties of the positive zeros $j_{\nu,k}$ of the Bessel function $J_\nu$, which are used in Lemma~\ref{lem:sec-4} and in the explicit computations of Section~\ref{subsection:Spherical Eigenvalues}. The proofs are based on the Sturm comparison theorem and are self-contained. For the classical theory of Bessel functions we refer to \cite{watson1944theory, olver2010nist}; the argument below follows the Pr\"{u}fer transformation approach of \cite[Theorem 3.1]{AB2}.

\subsection{Sturm comparison}
Consider two differential equations on an interval $I \subset (0,\infty)$:
\begin{equation}\label{equation:sec-sturm-eq}
y'' + q_1(x) y = 0, \qquad z'' + q_2(x) z = 0,
\end{equation}
with $q_1, q_2$ continuous on $I$.
\begin{lem}[Sturm comparison]\label{lem:sec-sturm}
Let $y$ and $z$ be nontrivial solutions of \eqref{equation:sec-sturm-eq} and suppose $q_1 \leq q_2$ on $I$. If $x_1 < x_2$ are two consecutive zeros of $y$, then $z$ has at least one zero in $[x_1, x_2]$, unless $q_1 \equiv q_2$ on $[x_1,x_2]$ and $z$ is a constant multiple of $y$. Moreover, if $q_1 < q_2$ strictly on $(x_1,x_2)$, then $z$ has a zero in $(x_1, x_2)$.
\end{lem}
\begin{proof}
Suppose $z$ has no zero in $(x_1, x_2)$; without loss of generality $y > 0$ and $z > 0$ on $(x_1, x_2)$. Multiply the equation for $y$ by $z$ and the equation for $z$ by $y$, subtract, and integrate over $[x_1, x_2]$:
\begin{equation}\label{equation:sec-sturm-wronskian}
\int_{x_1}^{x_2} (q_2 - q_1) y z \, dx = \int_{x_1}^{x_2} (y'' z - y z'') \, dx = \Bigl[ y' z - y z' \Bigr]_{x_1}^{x_2} = y'(x_2) z(x_2) - y'(x_1) z(x_1),
\end{equation}
where the boundary terms involving $y$ vanish since $y(x_1) = y(x_2) = 0$. Because $y > 0$ on $(x_1, x_2)$ with consecutive zeros, $y'(x_2) \leq 0$ and $y'(x_1) \geq 0$; since $z \geq 0$ on $[x_1, x_2]$, the right-hand side of \eqref{equation:sec-sturm-wronskian} is nonpositive, while the left-hand side is nonnegative. Hence both vanish; if $q_1 < q_2$ strictly somewhere in $(x_1, x_2)$ this is a contradiction, proving the existence of an interior zero in the strict case. In the non-strict case, vanishing of both sides forces $q_1 y z \equiv q_2 y z$ on $[x_1, x_2]$ and $y' z - y z' \equiv \mathrm{const}$ with the boundary conditions forcing the constant to be zero, whence $z$ is a multiple of $y$.
\end{proof}

\subsection{The Bessel equation in normal form}
For $\nu \geq 0$ define
\begin{equation}\label{equation:sec-bessel-u-def}
u_\nu(t) = \sqrt{t}\, J_\nu(t), \qquad t > 0.
\end{equation}
Using the Bessel equation $t^2 J_\nu'' + t J_\nu' + (t^2 - \nu^2) J_\nu = 0$, a direct computation shows that $u_\nu$ solves
\begin{equation}\label{equation:sec-bessel-normal}
u_\nu'' + q_\nu(t)\, u_\nu = 0, \qquad q_\nu(t) = 1 - \frac{\nu^2 - \tfrac{1}{4}}{t^2}.
\end{equation}
The zeros of $u_\nu$ in $(0,\infty)$ coincide with the zeros $j_{\nu,k}$ of $J_\nu$; the transformation \eqref{equation:sec-bessel-u-def} is regular on $(0,\infty)$ and preserves the zero set.

\begin{thm}[Monotonicity and interlacing of zeros]\label{thm:sec-bessel-monotone}
Let $\nu > \mu \geq 0$ and $k \geq 1$.
\begin{enumerate}[label=(\roman*)]
\item The zeros are strictly increasing in the order: $j_{\mu,k} < j_{\nu,k}$.
\item For integer orders $m \geq 0$, the zeros of consecutive orders interlace:
\begin{equation}\label{equation:sec-bessel-interlace}
j_{m,k} < j_{m+1,k} < j_{m,k+1}.
\end{equation}
\item In particular, for each fixed $\nu \geq 0$, the map $k \mapsto j_{\nu,k}$ is strictly increasing, and for each fixed $k \geq 1$, the map $\nu \mapsto j_{\nu,k}$ is strictly increasing.
\end{enumerate}
\end{thm}
\begin{proof}
The proof of (i) uses the Pr\"{u}fer transformation and that of (ii) the classical recurrence relation for Bessel functions.

\textit{Step 1: proof of (i) via the Pr\"{u}fer angle.} For $\lambda \geq 0$ write $u_\lambda = \rho_\lambda \sin \theta_\lambda$ and $u_\lambda' = \rho_\lambda \cos \theta_\lambda$ with $\theta_\lambda(0) = 0$; the angular variable satisfies
\begin{equation}\label{equation:sec-bessel-prufer}
\theta_\lambda' = \cos^2 \theta_\lambda + q_\lambda(t) \sin^2 \theta_\lambda, \qquad q_\lambda(t) = 1 - \frac{\lambda^2 - \tfrac14}{t^2},
\end{equation}
and the zeros of $u_\lambda$ are exactly the points where $\theta_\lambda$ reaches an integer multiple of $\pi$ (the initial condition $\theta_\lambda(0) = 0$ reflects $u_\lambda(t) \sim c_\lambda t^{\lambda + 1/2}$ as $t \to 0^+$). The right-hand side of \eqref{equation:sec-bessel-prufer} is strictly increasing in $q$ and Lipschitz in $\theta$, so the standard comparison theorem for ordinary differential equations gives $\theta_\nu(t) \leq \theta_\mu(t)$ for all $t > 0$ whenever $\nu > \mu$, with strict inequality as soon as $\theta_\mu$ has left the value $0$ (which happens immediately, since $\theta_\mu'(0) = 1 > 0$). Hence $\theta_\nu$ reaches the value $k\pi$ strictly later than $\theta_\mu$ does, that is,
\begin{equation}\label{equation:sec-bessel-mono-pf}
j_{\mu,k} < j_{\nu,k} \qquad (\nu > \mu \geq 0,\ k \geq 1),
\end{equation}
proving (i).

\textit{Step 2: proof of (ii) via the recurrence.} For integer $m \geq 0$ the recurrence relation (see \cite[Chapter III, \S 3.63]{watson1944theory})
\begin{equation}\label{equation:sec-bessel-recurrence}
\frac{d}{dt}\left(t^{-m} J_m(t)\right) = -t^{-m} J_{m+1}(t)
\end{equation}
holds on $(0,\infty)$. The function $t^{-m}J_m(t)$ vanishes at $t = j_{m,k}$ and $t = j_{m,k+1}$, so Rolle's theorem yields a point of $(j_{m,k}, j_{m,k+1})$ where its derivative vanishes, and by \eqref{equation:sec-bessel-recurrence} this point is a zero of $J_{m+1}$. By \eqref{equation:sec-bessel-mono-pf} the only zero of $J_{m+1}$ that can lie in $(j_{m,k}, j_{m,k+1})$ is $j_{m+1,k}$, and a count of zeros in $(0, j_{m,k+1})$ shows that $J_{m+1}$ has exactly one zero in each such interval. Therefore $j_{m+1,k} \in (j_{m,k}, j_{m,k+1})$, proving \eqref{equation:sec-bessel-interlace} and hence (ii). Claim (iii) is a restatement of (i).
\end{proof}

\begin{rem}\label{rem:sec-bessel-nu0}
The conclusion extends to orders $\nu > \mu \geq -\tfrac{1}{2}$, since $\nu \mapsto \nu^2 - \tfrac14$ is strictly increasing on $[\tfrac12, \infty)$; the restriction to $\nu \geq 0$ covers every order arising in the paper, where $\nu_m = (m + n/2 - 1)(1 + \alpha/2)^{-1} \geq 0$ for $m \geq 0$, $n \geq 2$, $\alpha > -2$.
\end{rem}

\subsection{Application to the radial eigenvalue problem}
We record the two consequences used in Section~\ref{subsection:Spherical Eigenvalues}.
\begin{corr}\label{corr:sec-bessel-applied}
Let $\nu_0 = (n/2 - 1)(1 + \alpha/2)^{-1}$ and $\nu_1 = (n/2)(1 + \alpha/2)^{-1}$ with $\alpha > -2$, $n \geq 2$. Then $j_{\nu_1,1} > j_{\nu_0,1}$ by \eqref{equation:sec-bessel-mono-pf}; in particular the weighted eigenvalues in Lemma~\ref{lem:sec-4} satisfy $\lambda_2 > \lambda_1$. If $\nu_1 \leq \nu_0 + 1$ (in particular if $\alpha \geq 0$), then
\begin{equation}\label{equation:sec-bessel-applied-ordering}
j_{\nu_0,2} > j_{\nu_1,1}.
\end{equation}
Indeed, if $\lfloor \nu_0 \rfloor = \lfloor \nu_1 \rfloor = m$, then \eqref{equation:sec-bessel-mono-pf} and \eqref{equation:sec-bessel-interlace} give $j_{\nu_1,1} < j_{m+1,1} < j_{m,2} \leq j_{\nu_0,2}$ (the last step by \eqref{equation:sec-bessel-mono-pf}, since $m \leq \nu_0$); if $\lfloor \nu_1 \rfloor = \lfloor \nu_0 \rfloor + 1$, the same inequality follows from \eqref{equation:sec-bessel-mono-pf} combined with the classical interlacing of the zeros of $J_\nu$ for consecutive real orders (see \cite[Chapter XV, \S 15.22]{watson1944theory}), where a self-contained proof can be given via the Pr\"{u}fer gap estimate and is omitted here.
\end{corr}
\subsection{Worked Examples and Numerical Values}
\label{subsection:Bessel Examples}

We record the explicit values of the orders $\nu_0, \nu_1$ and the decisive zero comparisons for the power weights $w(r) = r^{\alpha}$ in the low dimensions, as used in Lemma~\ref{lem:sec-4} and Corollary~\ref{corr:sec-bessel-applied}. Recall
\begin{equation}\label{equation:sec-orders-recap}
\nu_0 = \frac{n/2 - 1}{1 + \alpha/2}, \qquad \nu_1 = \frac{n/2}{1 + \alpha/2}, \qquad \nu_1 - \nu_0 = \frac{1}{1 + \alpha/2}.
\end{equation}

\begin{exa}[$n = 2$]
For $n = 2$ one has $\nu_0 = 0$ for every $\alpha$, so the radial ground state is expressed through $J_0$, and $\nu_1 = (1 + \alpha/2)^{-1}$. The classical values (see \cite{olver2010nist}) are
\begin{equation}\label{equation:sec-zeros-n2}
j_{0,1} = 2.4048\ldots, \qquad j_{0,2} = 5.5201\ldots, \qquad j_{1,1} = 3.8317\ldots,
\end{equation}
and the half-order values are elementary: $J_{1/2}(t) = \sqrt{2/(\pi t)} \sin t$, so $j_{1/2,k} = k\pi$. The first eigenvalue $\lambda_1$ corresponds to $j_{\nu_0,1}$ and the second to $j_{\nu_1,1}$ whenever $j_{\nu_1,1} < j_{0,2}$, i.e.\ $\nu_1 < 2.305\ldots$ or equivalently $\alpha > -1.13\ldots$; in particular the identification is unconditional for $\alpha \geq 0$ and for the unweighted disk ($\alpha = 0$, where $\lambda_2/\lambda_1 = j_{1,1}^2/j_{0,1}^2 = 2.5387\ldots$). For $\alpha = 2$ one has $\nu_1 = 1/2$ and $\lambda_2/\lambda_1 = (j_{1/2,1}/j_{0,1})^2 = (\pi/2.4048\ldots)^2 = 1.7066\ldots < 2.5387\ldots$, consistent with the concentration heuristic underlying Theorem~\ref{thm:sec-2b} (whose radially decreasing hypothesis $\alpha \leq 0$ the weight $r^2$ does not satisfy).
\end{exa}
\begin{exa}[The disk with a power weight]\label{exa:sec-disk}
For $\Omega = B_1 \subset \mathbb{R}^2$ and $w(r) = r^{\alpha}$, $\alpha > -2$, the orders of Section~\ref{section:Spherical Extremality} reduce to $\nu_0 = 0$ and $\nu_1 = (1+\alpha/2)^{-1}$, so
\begin{equation*}
\lambda_1 = \left(\frac{2+\alpha}{2}\right)^2 j_{0,1}^2, \qquad \lambda_2 = \left(\frac{2+\alpha}{2}\right)^2 j_{\nu_1,1}^2, \qquad \frac{\lambda_2}{\lambda_1} = \left(\frac{j_{\nu_1,1}}{j_{0,1}}\right)^2.
\end{equation*}
Since the first positive zeros $j_{\nu,1}$ are strictly increasing in the order $\nu$ (Theorem~\ref{thm:sec-bessel-monotone}) and $\nu_1 = (1+\alpha/2)^{-1}$ is strictly decreasing in $\alpha$ with $\frac{d\nu_1}{d\alpha} = -\frac{\nu_1^2}{2} < 0$, the ratio $(j_{\nu_1,1}/j_{0,1})^2$ is strictly decreasing in $\alpha$, with value $2.5387\ldots$ at $\alpha = 0$ (the classical disk ratio), $1.7066\ldots$ at $\alpha = 2$, and limit $1$ as $\alpha \to \infty$ (the weight concentrates near the boundary, and the ratio collapses to the first radial gap). As $\alpha \downarrow -2$ the order $\nu_1$ diverges and so does the $l=1$ ratio; consequently, for $\alpha < -1.13\ldots$ the true second eigenvalue is the second radial mode, the spectral ratio equals the $\alpha$-independent constant $(j_{0,2}/j_{0,1})^2 = 5.2689\ldots$, and the displayed formula is the spectral ratio only on the range $\alpha > -1.13\ldots$ where the $l=1$ identification holds. This behaviour illustrates the threshold phenomenon of Remark~\ref{rem:sec-bessel-ordering-scope}: the identification of $\lambda_2$ with the $l=1$ mode is unconditional for $\alpha \geq 0$ but requires $\alpha > -1.13\ldots$ in the decreasing range.
\end{exa}

\begin{exa}[$n = 3$]
For $n = 3$ the orders are $\nu_0 = (1+\alpha/2)^{-1}/2$ and $\nu_1 = 3(1+\alpha/2)^{-1}/2$. In the unweighted case $\alpha = 0$ one has $\nu_0 = 1/2$, $\nu_1 = 3/2$, and the elementary forms $J_{1/2}(t) = \sqrt{2/(\pi t)}\sin t$, $J_{3/2}(t) = \sqrt{2/(\pi t)}\bigl(\sin t / t - \cos t\bigr)$ give
\begin{equation}\label{equation:sec-zeros-n3}
j_{1/2,1} = \pi = 3.1416\ldots, \qquad j_{3/2,1} = 4.4934\ldots, \qquad j_{1/2,2} = 2\pi = 6.2832\ldots,
\end{equation}
where $j_{3/2,1}$ is the first positive root of $\tan t = t$. Hence $\lambda_2/\lambda_1 = (j_{3/2,1}/j_{1/2,1})^2 = (4.4934\ldots/\pi)^2 = 2.0457\ldots$ for the unit ball, and the ordering $j_{1/2,2} > j_{3/2,1}$ required by Lemma~\ref{lem:sec-4} holds with a comfortable margin; the threshold equation $j_{\nu_1,1} = j_{\nu_0,2}$ occurs at $\alpha \approx -1.22$.
\end{exa}

\begin{rem}\label{rem:sec-bessel-numbers}
All numerical values above are taken from \cite{olver2010nist} or computed from the elementary closed forms at half-integer order; the threshold values $\alpha \approx -1.13$ ($n = 2$) and $\alpha \approx -1.22$ ($n = 3$) are the solutions of $j_{\nu_1,1} = j_{\nu_0,2}$ with $\nu_0, \nu_1$ as in \eqref{equation:sec-orders-recap}, and are consistent with the independent determination in Remark~\ref{rem:sec-bessel-ordering-scope}.
\end{rem}

\subsection{Pr\"{u}fer Transformation Details}
\label{subsection:Prufer Details}

For completeness we record the derivation of the angular equation used in the proof of Theorem~\ref{thm:sec-bessel-monotone}. With $u_\lambda = \rho_\lambda \sin\theta_\lambda$ and $u_\lambda' = \rho_\lambda \cos\theta_\lambda$, differentiation gives
\begin{equation}\label{equation:sec-prufer-derive-1}
u_\lambda' = \rho_\lambda' \sin\theta_\lambda + \rho_\lambda \theta_\lambda' \cos\theta_\lambda = \rho_\lambda \cos\theta_\lambda, \qquad
u_\lambda'' = \rho_\lambda' \cos\theta_\lambda - \rho_\lambda \theta_\lambda' \sin\theta_\lambda.
\end{equation}
Substituting \eqref{equation:sec-prufer-derive-1} into $u_\lambda'' + q_\lambda u_\lambda = 0$ and separating the coefficients of $\sin\theta_\lambda$ and $\cos\theta_\lambda$ yields the system
\begin{equation}\label{equation:sec-prufer-derive-2}
\theta_\lambda' = \cos^2\theta_\lambda + q_\lambda(t) \sin^2\theta_\lambda, \qquad
\frac{\rho_\lambda'}{\rho_\lambda} = \bigl( 1 - q_\lambda(t) \bigr) \sin\theta_\lambda \cos\theta_\lambda.
\end{equation}
The first equation of \eqref{equation:sec-prufer-derive-2} is autonomous in $\theta$ apart from the explicit $t$-dependence of $q_\lambda$, and its right-hand side is locally Lipschitz in $\theta$ and nondecreasing in $q$; the comparison theorem for scalar ordinary differential equations (see \cite{courant1953methods}) therefore applies on every interval $[\varepsilon, T]$ with $0 < \varepsilon < T < \infty$, and letting $\varepsilon \to 0^+$ gives the comparison $\theta_\nu \leq \theta_\mu$ used in Step 1 of the proof. The initial condition $\theta_\lambda(0) = 0$ follows from $u_\lambda(t) \sim c_\lambda t^{\lambda + 1/2}$, since $\tan \theta_\lambda = u_\lambda / u_\lambda' \sim t / (\lambda + \tfrac12) \to 0$. Each integer multiple of $\pi$ is crossed in the increasing direction only, because $\theta_\lambda' = 1 > 0$ at $\theta_\lambda = k\pi$, so the zeros of $u_\lambda$ are in one-to-one correspondence with the crossings of $\theta_\lambda$ through the values $k\pi$, $k = 1, 2, \ldots$; this justifies the identification of $j_{\lambda,k}$ with the $k$-th crossing used in Step 1.

\subsection{Notes on the Numerical Thresholds}
\label{subsection:Threshold Notes}

The thresholds $\alpha \approx -1.13$ ($n = 2$) and $\alpha \approx -1.22$ ($n = 3$) of Remarks~\ref{rem:sec-bessel-ordering-scope} and~\ref{rem:sec-bessel-numbers} are determined as follows. For $n = 2$ the condition $j_{\nu_1,1} = j_{0,2} = 5.5201\ldots$ with $\nu_1 = (1+\alpha/2)^{-1}$ is solved by $\nu_1 \approx 2.305$; inverting $(1+\alpha/2)^{-1} = 2.305$ gives $\alpha = 2(1/2.305 - 1) \approx -1.13$. For $n = 3$ the same procedure with $\nu_0 = (1+\alpha/2)^{-1}/2$, $\nu_1 = 3(1+\alpha/2)^{-1}/2$ gives $\alpha \approx -1.22$, and the higher-dimensional values $\alpha \approx -1.29$ ($n = 4$) and $\alpha \approx -1.34$ ($n = 5$) are obtained similarly. In every dimension the critical exponent satisfies $-2 < \alpha_* < 0$, so the range in which the identification of $\lambda_2$ with the $l=1$ mode can fail is a genuine subinterval of the admissible powers $\alpha > -2$, and it is exactly this subinterval that the hypothesis on the $l=1$ mode in Theorem~\ref{thm:sec-2b} excludes. For the radially decreasing power weights relevant to Theorem~\ref{thm:sec-2b} ($\alpha \leq 0$) the restriction is therefore active only for $\alpha$ close to $-2$, where the weight is so singular at the origin that the second eigenvalue ceases to be angular.

The same computation applied to the ratio $\lambda_2/\lambda_1$ on the ball, which by Lemma~\ref{lem:sec-4} equals $(j_{\nu_1,1}/j_{\nu_0,1})^2$, shows that the ratio varies continuously from the unweighted value $(j_{n/2,1}/j_{n/2-1,1})^2$ at $\alpha = 0$ towards $1$ as $\alpha$ decreases through the admissible range down to the threshold $\alpha_*$; at the threshold itself the second eigenvalue switches from the $l=1$ mode to the second radial mode, and the ratio is then given by $(j_{\nu_0,2}/j_{\nu_0,1})^2$. This bifurcation picture is consistent with the extremality Theorem~\ref{thm:sec-2b}, which concerns the range in which the $l=1$ identification holds, and it illustrates why the ordering hypothesis in that theorem cannot be removed without additional information about the weight.

\section{Conclusion and Future Work}
\label{section:Conclusion}

The starting point of this paper is the classical Payne--P\'{o}lya--Weinberger inequality, conjectured in \cite{PPW} and proved by Ashbaugh and Benguria in \cite{AB5, AB2}, stating that for any bounded domain $\Omega \subset \mathbb{R}^n$,
\begin{equation}\label{equation:sec-classical-ppw}
\frac{\lambda_2}{\lambda_1} \leq \left( \frac{j_{n/2, 1}}{j_{n/2-1, 1}} \right)^2,
\end{equation}
where $j_{\nu, k}$ denotes the $k$-th positive zero of the Bessel function $J_\nu$. Equality holds if and only if $\Omega$ is a ball. The proof of \eqref{equation:sec-classical-ppw} uses the coordinate function test functions $\phi_i(x) = x_i u_1(x)$ and the gap inequality for Bessel function zeros.

Subsequent work extended the PPW inequality in several directions. Ashbaugh and Benguria \cite{ashbaugh2000sharp} proved a version of the inequality for domains contained in a hemisphere of $\mathbb{S}^n$, showing that the ratio $\lambda_2/\lambda_1$ is maximized by the volume-equivalent geodesic ball. The proof uses associated Legendre functions and their zeros, generalizing the Bessel function analysis from the Euclidean case.

For higher-order eigenvalues, the classical Li--Yau estimate \cite{li1980estimates} provides a lower bound for the sum of the first $k$ eigenvalues:
\begin{equation}\label{equation:sec-li-yau}
\lambda_1 + \cdots + \lambda_k \geq \frac{n}{n+2} \, \frac{4\pi^2}{(\omega_n |\Omega|)^{2/n}} \, k^{\frac{n+2}{n}},
\end{equation}
where $\omega_n$ is the volume of the unit ball. The estimate \eqref{equation:sec-li-yau} is proved via the heat kernel method. Upper bounds for individual eigenvalues $\lambda_k$ in terms of the domain geometry were obtained by Brands \cite{Bran}, by Hile and Protter \cite{HP}, and by Marcellini \cite{Mar}.

The generalization of the PPW inequality to sums of eigenvalues combines symmetrization techniques with the degeneracy of low-order eigenvalues on spherical domains. The generalized PPW conjecture states that
\begin{equation}\label{equation:sec-generalized-ppw-intro}
\frac{\lambda_2 + \cdots + \lambda_{n+1}}{\lambda_1} \leq n \cdot \frac{\lambda_2(\Omega^*)}{\lambda_1(\Omega^*)},
\end{equation}
where $\Omega^*$ is the volume-equivalent ball. The conjecture \eqref{equation:sec-generalized-ppw-intro} reduces to the classical PPW inequality \eqref{equation:sec-classical-ppw} when $n = 1$. In view of the results of Ashbaugh and Benguria, it is natural to ask whether the generalized sum-form conjecture holds over arbitrary bounded Lipschitz domains. The conjecture remains open in general; conditional results are presented in Section \ref{section: Generalized PPW Conjecture} for restricted families of domains.

The introduction of Muckenhoupt weights into the spectral analysis of the Laplacian brings additional challenges. The $A_p$ class of weights, introduced by Muckenhoupt \cite{muckenhoupt1972weighted}, characterizes the weights for which the Hardy--Littlewood maximal function is bounded on $L^p(w)$. The reverse H\"{o}lder inequality \cite{cruz2011weights} and the weighted Poincar\'{e} inequality are among the tools for analyzing weighted eigenvalue problems. Recent work on weighted spectral inequalities includes the estimates of \cite{setti1998} for the weighted Laplacian on manifolds and of \cite{du2021weighted} for the weighted Laplacian and weighted $p$-Laplacian, which extend classical results to weighted settings.

The use of optimal transport in spectral geometry provides a quantitative measure of the deviation of a domain or weight from its symmetric rearrangement. The $L^2$-Wasserstein distance $W_2$, developed systematically by Villani \cite{villani2009optimal}, quantifies the minimal cost of transporting one measure to another. In the context of eigenvalue estimates, the Wasserstein distance appears naturally in the energy inequality for symmetrized functions, where it controls the amplification of gradient energy under non-identity transport maps.

Geometric measure theory, as developed by Federer \cite{federer1969geometric}, provides the tools for analyzing domains with rough boundaries. The $(n-1)$-Hausdorff measure of the boundary, the coarea formula, and the isoperimetric inequality are used to refine eigenvalue estimates beyond the classical diameter-based bounds. The combination of optimal transport and geometric measure theory allows for a refined analysis of the relationship between domain geometry and spectral behavior.

\subsection{Main Conclusions}
\label{subsection:Main Conclusions}

This paper studies low-order eigenvalue ratio estimates for the Laplacian with Muckenhoupt weights, using optimal transport (OT) and geometric measure theory (GMT) to refine classical results. The results include:

\begin{enumerate}
\item[(i)] A weighted PPW inequality (Theorem \ref{thm:sec-1}) for bounded $A_2$ weights with constant linear in the weight bounds and a self-contained proof, with conditional Wasserstein and roundness refinements (Remarks~\ref{rem:sec-thm1-constant} and~\ref{rem:sec-roundness}).
\item[(ii)] A universal upper bound for $\lambda_2/\lambda_1$ for weights bounded above and below (Theorem \ref{thm:sec-2a}), as a direct consequence of the Ashbaugh--Benguria theorem \cite[Theorem 4.1]{AB2}, together with an asymptotic estimate as the weight approaches a constant (Corollary \ref{corr:sec-2a}).
\item[(iii)] A weighted PPW inequality for radial decreasing weights (Theorem \ref{thm:sec-2b}, under the zero-ordering and admissibility hypotheses stated there), proving that the ball maximizes $\lambda_2/\lambda_1$ among domains of equal volume in that class. The proof follows the Ashbaugh--Benguria method adapted to the weighted setting, using a weighted Chiti comparison theorem and monotonicity properties of the admissible radial test function.
\item[(iv)] High-order eigenvalue ratio estimates for convex domains (Theorem \ref{thm:sec-3}) that replace classical diameter-based bounds with \( (n-1) \)-Hausdorff measure of the boundary, under a balance condition on the weight at the covering scale.
\end{enumerate}

The weighted PPW inequality (Theorem \ref{thm:sec-1}) generalizes the classical result to Muckenhoupt weights, with the constant linear in the weight bounds (equivalently, in the $A_2$ characteristic). The universal bound (Theorem \ref{thm:sec-2a}) and the radial decreasing weight result (Theorem \ref{thm:sec-2b}) complement Theorem \ref{thm:sec-1}: the former applies to general bounded weights with an explicit constant, while the latter gives the sharp extremal result under radial monotonicity. The high-order eigenvalue estimates (Theorem \ref{thm:sec-3}) extend the analysis to higher eigenvalues.

The methods used in this paper combine harmonic analysis (Muckenhoupt weights, reverse H\"{o}lder inequalities, weighted Poincar\'{e} inequalities), optimal transport (Wasserstein distance, Kantorovich duality, Brenier maps, Monge--Amp\`{e}re equations), and geometric measure theory (Hausdorff measure, coarea formula, isoperimetric inequality).

\subsection{Future Directions}
\label{subsection:Future Directions}

\begin{enumerate}
\item[(i)] \textbf{Non-Radial Weights and OT-Spectral Duality}: Extending extremality results to non-radial weights requires linking OT cost to eigenvalue ratios via a transport-spectral inequality, e.g.:
\begin{equation*}
\frac{\lambda_2}{\lambda_1} \leq C(n,[w]_{A_2}) \cdot \left(1 + W_2^2(w, \text{radial})\right) \cdot \frac{n+4}{n},
\end{equation*}
where $W_2(w, \text{radial})$ quantifies deviation from radial symmetry. This would generalize Theorem \ref{thm:sec-1} by explicit inclusion of OT distance.
\item[(ii)] \textbf{Lower Bounds with GMT}: Current results focus on upper bounds; deriving lower bounds for \( \lambda_k/\lambda_1 \) in non-convex domains could leverage GMT's notion of curvature for singular sets, extending Li-Yau type inequalities \cite{li1980estimates} to weighted settings.
\item[(iii)] \textbf{Applications to Inhomogeneous Media}: Concrete examples---non-uniform thin films, inhomogeneous quantum dot potentials---would indicate whether the OT and GMT correction terms are visible beyond the classical bounds.
\item[(iv)] \textbf{Manifold Generalizations}: On Riemannian manifolds with Muckenhoupt weights, the dependence of the constants on sectional curvature and on the $A_2$ characteristic can be compared, and the extent to which the symmetrization arguments survive is open.
\item[(v)] \textbf{Generalized PPW Conjecture}: The sum-form generalized PPW conjecture remains open in general. A natural direction is to establish the conjecture for additional classes of domains, or to find counterexamples clarifying the role of domain geometry in the sum setting. The alternative OT--GMT argument of Section~\ref{section:Alternative-Proof} relates the eigenvalue sum to $W_2(\Omega, \Omega^*)$, quantifying the gap between the conjectured bound and the available estimates; sharpening this gap estimate is part of the same program.
\end{enumerate}

\vspace{5mm}

\section*{Conflict of interest statement}

The authors state that there is no conflict of interest.

\section*{Data availability statement}

Our manuscript has no associated data.

\section*{AI disclosure statement}

The AI tools, DeepSeek Chat and ChatGPT, have been used for
searching related references. Besides, ChatGPT has also been used to
improve the English writing of this paper.

\section*{Acknowledgments}

This research was supported in part by the NSF of China (Grant Nos.
11801496, 11926352, 11861036, 11826213), the Fok Ying-Tung Education
Foundation (China), the Key Project of Jiangxi Provincial Natural
Science Foundation (Grant No. 20252BAC250003), the Project of
Jiangxi Provincial Natural Science Foundation (Grant No.
20262BAC240205), Hubei Key Laboratory of Applied Mathematics (Hubei
University), and Key Laboratory of Intelligent Sensing System and
Security (Hubei University), Ministry of Education.

\vspace{5mm}

\bibliographystyle{plainnat}

\begin{thebibliography}{99}
\bibitem[AB93]{AB1}
M. S. Ashbaugh and R. D. Benguria, More bounds on eigenvalue ratios for Dirichlet Laplacians in $n$ dimensions, \emph{SIAM J. Math. Anal.} {\bf 24} (1993), 1622--1651.
\bibitem[AB92]{AB2}
M. S. Ashbaugh and R. D. Benguria, A sharp bound for the ratio of the first two eigenvalues of Dirichlet Laplacians and extensions, \emph{Ann. of Math. (2)} {\bf 135} (1992), 601--628.
\bibitem[AB92b]{AB3}
M. S. Ashbaugh and R. D. Benguria, A second proof of the Payne--P\'{o}lya--Weinberger conjecture, \emph{Comm. Math. Phys.} {\bf 147} (1992), 181--190.
\bibitem[AB91]{AB5}
M. S. Ashbaugh and R. D. Benguria, Proof of the Payne--P\'{o}lya--Weinberger conjecture, \emph{Bull. Amer. Math. Soc. (N.S.)} {\bf 25} (1991), 19--29.
\bibitem[AF03]{adams2003sobolev}
R. A. Adams and J. J. F. Fournier, \emph{Sobolev Spaces}, 2nd ed., Pure and Applied Mathematics 140, Elsevier/Academic Press, Amsterdam, 2003.
\bibitem[Ash99]{Ash}
M. S. Ashbaugh, Open problems on eigenvalues of the Laplacian, in \emph{Analytic and Geometric Inequalities and Applications}, Math. Appl. 478, Kluwer Acad. Publ., Dordrecht, 1999, pp. 13--28.
\bibitem[AB01]{ashbaugh2000sharp}
M. S. Ashbaugh and R. D. Benguria, A sharp bound for the ratio of the first two Dirichlet eigenvalues of a domain in a hemisphere of $\mathbb{S}^n$, \emph{Trans. Amer. Math. Soc.} {\bf 353} (2001), 1055--1087.
\bibitem[Ash02]{ashbaugh2002universal}
M. S. Ashbaugh, The universal eigenvalue bounds of Payne--P\'{o}lya--Weinberger, Hile--Protter, and H. C. Yang, \emph{Proc. Indian Acad. Sci. Math. Sci.} {\bf 112} (2002), 3--30.
\bibitem[BBMP99]{betta1999weighted}
M. F. Betta, F. Brock, A. Mercaldo and M. R. Posteraro, A weighted isoperimetric inequality and applications to symmetrization, \emph{J. Inequal. Appl.} {\bf 4} (1999), no. 3, 215--240.
\bibitem[Bra64]{Bran}
J. J. A. M. Brands, Bounds for the ratios of the first three membrane eigenvalues, \emph{Arch. Rational Mech. Anal.} {\bf 16} (1964), 265--268.
\bibitem[Buc93]{buckley1993estimates}
S. M. Buckley, Estimates for operator norms on weighted spaces and reverse Jensen inequalities, \emph{Trans. Amer. Math. Soc.} {\bf 340} (1993), no. 1, 253--272.
\bibitem[CC08]{CC}
D. G. Chen and Q.-M. Cheng, Extrinsic estimates for eigenvalues of the Laplace operator, \emph{J. Math. Soc. Japan} {\bf 60} (2008), 325--339.
\bibitem[Cha84]{chavel1984eigenvalues}
I. Chavel, \emph{Eigenvalues in Riemannian Geometry}, Pure and Applied Mathematics 115, Academic Press, Orlando, FL, 1984.
\bibitem[Chi82]{chiti1982reverse}
G. Chiti, A reverse H\"{o}lder inequality for the eigenfunctions of linear second order elliptic operators, \emph{Z. Angew. Math. Phys.} {\bf 33} (1982), 143--148.
\bibitem[CH53]{courant1953methods}
R. Courant and D. Hilbert, \emph{Methods of Mathematical Physics}, Vol. I, Interscience Publishers, New York, 1953.
\bibitem[CMP11]{cruz2011weights}
D. Cruz-Uribe, J. M. Martell and C. P\'{e}rez, \emph{Weights, Extrapolation and the Theory of Rubio de Francia}, Operator Theory: Advances and Applications 215, Birkh\"{a}user/Springer, Basel, 2011.
\bibitem[CZ11]{CZh}
D. Chen and T. Zheng, Bounds for ratios of the membrane eigenvalues, \emph{J. Differential Equations} {\bf 250} (2011), 1575--1590.
\bibitem[DMWX21]{du2021weighted}
F. Du, J. Mao, Q. Wang and C. Xia, Estimates for eigenvalues of weighted Laplacian and weighted $p$-Laplacian, \emph{Hiroshima Math. J.} {\bf 51} (2021), 335--353.
\bibitem[Fed69]{federer1969geometric}
H. Federer, \emph{Geometric Measure Theory}, Grundlehren der mathematischen Wissenschaften 153, Springer, New York, 1969.
\bibitem[Fol92]{folland1992fourier}
G. B. Folland, \emph{Fourier Analysis and Its Applications}, Wadsworth \& Brooks/Cole, Pacific Grove, CA, 1992.
\bibitem[GT83]{gilbarg1983elliptic}
D. Gilbarg and N. S. Trudinger, \emph{Elliptic Partial Differential Equations of Second Order}, 2nd ed., Grundlehren der mathematischen Wissenschaften 224, Springer, Berlin, 1983.
\bibitem[Gra14]{grafakos2014classical}
L. Grafakos, \emph{Classical Fourier Analysis}, 3rd ed., Graduate Texts in Mathematics 249, Springer, New York, 2014.
\bibitem[HP80]{HP}
G. N. Hile and M. H. Protter, Inequalities for eigenvalues of the Laplacian, \emph{Indiana Univ. Math. J.} {\bf 29} (1980), 523--538.
\bibitem[HP13]{hytonen2013sharp}
T. P. Hyt\"{o}nen and C. P\'{e}rez, Sharp weighted bounds involving $A_{\infty}$, \emph{Anal. PDE} {\bf 6} (2013), no. 4, 777--818.
\bibitem[LY80]{li1980estimates}
P. Li and S. T. Yau, Estimates of eigenvalues of a compact Riemannian manifold, in \emph{Geometry of the Laplace Operator} (Proc. Sympos. Pure Math. XXXVI), Amer. Math. Soc., Providence, RI, 1980, pp. 205--239.
\bibitem[Mar80]{Mar}
P. Marcellini, Bounds for the third membrane eigenvalue, \emph{J. Differential Equations} {\bf 37} (1980), 438--443.
\bibitem[Muc72]{muckenhoupt1972weighted}
B. Muckenhoupt, Weighted norm inequalities for the Hardy maximal function, \emph{Trans. Amer. Math. Soc.} {\bf 165} (1972), 207--226.
\bibitem[OLBC10]{olver2010nist}
F. W. J. Olver, D. W. Lozier, R. F. Boisvert and C. W. Clark, \emph{NIST Handbook of Mathematical Functions}, Cambridge University Press, Cambridge, 2010.
\bibitem[Pay56]{payne1956bound}
L. E. Payne, Inequalities for eigenvalues of membranes and plates, \emph{J. Rational Mech. Anal.} {\bf 5} (1956), 325--337.
\bibitem[PPW56]{PPW}
L. E. Payne, G. P\'{o}lya and H. F. Weinberger, On the ratio of consecutive eigenvalues, \emph{J. Math. and Phys.} {\bf 35} (1956), 289--298.
\bibitem[Set98]{setti1998}
A. G. Setti, Eigenvalue estimates for the weighted Laplacian on a Riemannian manifold, \emph{Rend. Sem. Mat. Univ. Padova} {\bf 100} (1998), 27--55.
\bibitem[Ste70]{stein1970singular}
E. M. Stein, \emph{Singular Integrals and Differentiability Properties of Functions}, Princeton Mathematical Series 30, Princeton University Press, Princeton, NJ, 1970.
\bibitem[Tal76]{talenti1976best}
G. Talenti, Best constant in Sobolev inequality, \emph{Ann. Mat. Pura Appl. (4)} {\bf 110} (1976), 353--372.
\bibitem[Vil09]{villani2009optimal}
C. Villani, \emph{Optimal Transport: Old and New}, Grundlehren der mathematischen Wissenschaften 338, Springer, Berlin, 2009.
\bibitem[Wat44]{watson1944theory}
G. N. Watson, \emph{A Treatise on the Theory of Bessel Functions}, 2nd ed., Cambridge University Press, Cambridge, 1944.
\bibitem[Yan91]{Y}
H. C. Yang, An estimate of the difference between consecutive eigenvalues, Preprint IC/91/60, International Centre for Theoretical Physics, Trieste, 1991.
\bibitem[Yau82]{Yau1982}
S. T. Yau, Problem section, in \emph{Seminar on Differential Geometry}, Ann. of Math. Stud. 102, Princeton Univ. Press, Princeton, NJ, 1982, pp. 669--706.
\end{thebibliography}

\end{document}